\documentclass[english,12pt,leqno]{article}
\usepackage[latin1]{inputenc}

\usepackage{fullpage}

\usepackage{tipa}
\usepackage{dsfont}
\usepackage{tikz-cd}
\usepackage{float}
\usepackage{graphicx}

\usepackage{amsxtra}
\usepackage{amsmath}
\usepackage{amssymb}
\usepackage{amsfonts}
\usepackage[all,2cell]{xy}
\usepackage{mathrsfs}
\usepackage{amsthm}
\usepackage{enumitem}
\usepackage{hyperref}
\usepackage{subfigure}
\usepackage{mathabx}
\usepackage{euscript}
\usepackage[bbgreekl]{mathbbol}
\usepackage{pgfplots}
 
\usetikzlibrary{patterns}
\usetikzlibrary{calc}

\usepackage[OT2,T1]{fontenc}
\DeclareSymbolFont{cyrletters}{OT2}{wncyr}{m}{n}
\DeclareMathSymbol{\Sha}{\mathalpha}{cyrletters}{"58}
 
\makeatletter
\hypersetup{
    unicode=true,           
    pdftoolbar=true,        
    pdfmenubar=true,        
    pdffitwindow=false,     
    pdfstartview={FitH},    
    pdftitle={\@title},     
    pdfauthor={\@author},   
    pdfsubject={},          
    pdfcreator={},          
    pdfproducer={},         
    pdfkeywords={},         
    pdfnewwindow=true,      
    colorlinks,             
    linkcolor=black,        
    citecolor=black,        
    filecolor=black,        
    urlcolor=black          
}
\makeatother

\newtheorem{thm}[equation]{Theorem}
\newtheorem{cor}[equation]{Corollary}
\newtheorem{lem}[equation]{Lemma}
\newtheorem{prop}[equation]{Proposition}

\newtheoremstyle{example}{\topsep}{\topsep}%
     {}
     {}
     {\bfseries}
     {.}
     {2pt}
     {\thmname{#1}\thmnumber{ #2}\thmnote{ #3}}

   \theoremstyle{example}
   
   \newtheorem{Defi}[equation]{Definition}
   \newtheorem{rem}[equation]{Remark}
   \newtheorem{rems}[equation]{Remarks}
   \newtheorem{prob}[equation]{Problem}
   \newtheorem{exas}[equation]{Examples}
   \newtheorem{ex}[equation]{Example}

  \newtheorem{observ}[equation]{Observation}
    \newtheorem{conj}[equation]{Conjecture}

\numberwithin{equation}{subsection}

\renewcommand{\theparagraph}{\Alph{paragraph}.} 
\def\AAA{\mathbb{A}}

\def\BB{\mathbb{B}}
\def\CC{\mathbb{C}}
\def\DD{\mathbb{D}}

 \def\HH{\mathbb{H}}

\def\RR{\mathbb{R}}
\def\ZZ{\mathbb{Z}}

\def\gen{\mathfrak{g}}

\def\len{\mathfrak{l}}

\def\pen{\mathfrak{p}}

\def\Ien{\mathfrak{I}}
\def\Pen{\mathfrak{P}}
\def\Ren{\mathfrak{R}}

\def\Ac{\mathcal{A}}

\def\Cc{\mathcal{C}}

\def\Dc{\mathcal{D}}
\def\Ec{\mathcal{E}}
\def\Fc{\mathcal{F}}
\def\Gc{\mathcal{G}}

\def\Lc{\mathcal{L}}

\def\Hc{\mathcal{H}}

\def\Pc{\mathcal{P}}

\def\Rc{\mathcal{R}}
\def\Sc{\mathcal{S}}
\def\Tc{\mathcal{T}}

\def\Vc{\mathcal{V}}

\def\Yc{\mathcal{Y}}

 \def\Aut{{\on{Aut}}}
 
 \def\BBB{{\BB\BB}}
 \def\Br{{\on{Br}}}
  \def\bw{{{\bf w}}}
 \def\bz{{{\bf z}}}

 \def\c{{\bf{c}}}

 \def\constr{{\on{constr}}}
 \def\Conv{{\on{Conv}}}

 \def\eps{{\varepsilon}}

  \def\Exit{{\on{Exit}}}
\def\Ext{{\on{Ext}}}

 \def\fd{{\on{fd}}}
 \def\fine{{\on{fine}}}
 \def\FPS{{\on{FPS}}}
 \def\Fun{{\on{Fun}}}

  \def\Hom{\operatorname{Hom}\nolimits}
 \def\Hur{{\on{Hur}}}
 
 \def\Id{\operatorname{Id}\nolimits}
 \newcommand{\iso}{\overset{\sim}{\longrightarrow}}
 
  \def\Im{{\on{Im}}}

 \def\k{\mathbf k}

 \def\lex{{\on{lex}}}
  \def\Lex{{\on{Lex}}}
 
   \def\lra{\longrightarrow}

\def\Ob{\operatorname{Ob}\nolimits}
\def\on{\operatorname}
\def\ol{\overline}
\def\op{{\on{op}}}
\def\OR{{\on{OR}}}

 \def\PB{{\on{PB}}}
 \def\Perv{{\on{Perv}}}

\def\Re{{\on{Re}}}
 \def\Rep{\on{{Rep}}}

 \def\sa{{\on{sa}}}
  
  \def\sgn{{\on{sgn}}}
 \def\Sh{{\on{Sh}}}
 \def\ShD{{{\Sha} \on{D}}}

   \def\Sym{{\on{Sym}}}
   
\newcommand{\tb}{\tilde b}
\newcommand{\tc}{\tilde c}

\def\Top{{\Tc op}}
\def\Tor{{\on{Tor}}}
\def\tot{{\on{tot}}}

\def\ul{\underline}

\def\Vect{\on{Vect}}
\def\vzeta{{\vec{\zeta}}}
 
\def\wt{\widetilde}

\def\XO{{\on{XO}}}

\def\1{{\mathbf{1}}}
\def\2{{\mathbf{2}}}
\def\(({(\hskip -1mm (}
\def\)){)\hskip -1mm )}
 \def\be{\begin{equation}}
\def\ee{\end{equation}}
\def\ed{\end{document}}

\title{ Shuffle algebras and perverse sheaves}
\author{  Mikhail Kapranov, Vadim Schechtman}
 
\begin{document}

\maketitle

\thanks{\em to Kyoji Saito for his 75th birthday}

\begin{abstract}
 We relate shuffle algebras, as defined by Nichols, Feigin-Odesskii and Rosso, to perverse sheaves on symmetric
 products of the complex line (i.e., on the spaces of monic polynomials stratified by multiplicities of roots). More precisely,
 we construct an equivalence between:
 
 \begin{itemize}
 \item[(i)] Braided Hopf algebras of a certain type.
 
 \item[(ii)] Factorizable collections of perverse sheaves on all the symmetric products. 
 \end{itemize}
 Under this eqiuvalence, the Nichols algebra associated to an object $V$ corresponds to the collection of
 the intersection cohomology extensions of the local systems on the open configuration spaces associated to
 the tensor powers of $V$. 
  Our approach is based on using  real skeleta of  complex configuration spaces. 
 
\end{abstract}
 
\tableofcontents


\section*{Introduction}

\paragraph{First goal of the paper. }\label{par:first-goal}

Shuffle algebras were introduced by Feigin--Odesskii \cite{feigin-odesskii}  and Rosso \cite{rosso} but a closely related construction
goes back to the 1978 paper of Nichols \cite{nichols}. In modern terms,    an object $V$ of
a braided monoidal abelian category $(\Vc, \otimes, R)$  gives rise to three bialgebras in $\Vc$:
\begin{itemize}
\item[(1)] The {\em tensor algebra}  $T_!(V) = \bigoplus_{n\geq 0} V^{\otimes n}$ with tensor multiplication
and with  the comultipliction making $V$ primitive. 

\item[(2)] The {\em cotensor}, or ``big shuffle''  algebra $T_*(V)$ which is the same object $ \bigoplus_{n\geq 0} V^{\otimes n}$
but with the shuffle multiplication involving the $R$-matrices,
 see  \eqref {eq:shuffle-prod} below. We have a canonical morphism $\alpha: T_!(V)\to T_*(V)$.

\item[(3)]   The {\em Nichols algebra} (or  the ``true shuffle algebra'')  $T_{!*}(V)$ defined as the image of $\alpha$.
\end{itemize}

For instance, the positive part $U_q^+(\gen)$ of the quantum group with Cartan matrix $\|a_{ij}\|_{i,j=1}^r$ is obtained
as $T_{!*}(V)$ for $V$ an object of the category of $\ZZ^r$-graded vector spaces with braiding given by the factors $q^{a_{ij}}$,
see \cite{feigin-odesskii}, \cite{rosso} and Example \ref {ex:abel-mon} below  for  a precise formulation.

 \vskip .2cm
 
 The first goal of this paper is to relate the above structures with the familiar phenomenon in the theory of perverse sheaves.
 If $X$ is a complex algebraic variety and $j: U\hookrightarrow X$ is an affine open subset, then a local system
 $\Lc$ on $U$ gives three perverse sheaves on $X$:
 
 \begin{itemize}
 \item[($1'$)] The zero extension $j_!\Lc$.
 
 \item[($2'$)]The full direct image $j_*\Lc$.  We have a canonical morphism $\alpha': j_!\Lc\to j_*\Lc$.
 
 \item[($3'$)] The {\em intermediate} (or {\em intersection cohomology}) extension $j_{!*}\Lc$ which is the image of $\alpha'$. 
 \end{itemize}
 
 It turns out that ($1$)-($3$) can be interpreted as a particular case of ($1'$)-($3'$), if we take for  $X$  the symmetric products
 $\Sym^n(\CC)$, put 
 $U= \Sym^n_\neq (\CC)$ to be the generic stratum and take $\Lc$ to be the local system on $U$ corresponding to the
 object $V^{\otimes n}$ with its natural action of the braid group. Here we consider all $n$ at once and work with local systems
 and perverse sheaves with values in $\Vc$. 
 
 \vskip .2cm
 
 Our interpretation can be seen as a nonlinear analog of several  ``localization''  results,  starting from the Beilinson-Bernstein  picture  \cite{BB}, which deal 
 with {\em modules}
 over  a given (Lie)  algebra, not algebras themselves. In such
 results, ($1'$) and ($2'$) correspond to the induced and co-induced Verma
 modules, $\alpha'$ to the Shapovalov form and ($3'$) to the irreducible quotient of the Verma.

 \paragraph{Second goal: skeletology.} 
 Our second goal is 
   to  develop, a little further, the method of {\em totally real skeleta} in the study of perverse sheaves
   i.e., in describing the category $\Perv(X,S)$ for a stratified complex manifold $(X,S)$.     While the full range of
   this method is still not clear to us, in good cases   \cite{galligo-GM} 
   \cite{KS-hyp-arr} \cite{KS-graphs}  
 it works like this.  We choose
   an appropriate ``skeleton'' (totally real subanalytic subset) $i: K\hookrightarrow X$ with the properties:
   \begin{itemize}
   \item For any $\Fc\in\Perv(X,S)$, the complex $i^!\Fc$ is quasi-isomorphic to a single sheaf $R_K(\Fc)$. 
   Therefore, taking the stalk of $R_K(\Fc)$ at any given $x\in K$ gives an exact functor $R_{K,x}: \Perv(X,S)\to\Vect$.
   
   \item The functors $R_{K,x}$ take Verdier duality (denoted $\Fc\mapsto \DD \Fc$) to vector space duality. Therefore
   the sheaf structure on $R_K(\DD\Fc)$ defines an extra {\em cosheaf structure} on $R_K(\Fc)$
   with the same (co)stalks. 
   
   \item The sheaf and cosheaf structures on $R_K(\Fc)$ are compatible in a certain sense. Spelling out these
   compatibility conditions amounts to defining the notion of a {\em bisheaf} on $K$. The category $\Perv(X,S)$
   is then identified with the category of such bisheaves. 
     \end{itemize}
      In the above situation of  $X=\Sym^n(\CC)$ with its natural  ``diagonal'' stratification , we can take $K=\Sym^n(\RR)$.
     In this case a source of ``bisheaves'' can be found in  graded {\em bialgebras}: such a
      bialgebra gives a bisheaf on $K$ encoding the degree $n$ parts of the two Hochschild complexes:
    for
     the multiplication and comultiplication. 
     
     \paragraph { Description of results.} As far as the first goal goes, we establish, in Theorem  \ref{thm:localize}, an
      equivalence of categories $L$ between:
     \begin{itemize}
     \item [(1)] A certain class of graded Hopf algebras (called {\em primitive bialgebras}, see Definition \ref{def:prim-gen-bialg})
     in an abelian braided category $\Vc$.
     
     \item [(2)] Perverse sheaves on $\Sym(\CC) = \bigsqcup_{n\geq 0} \Sym^n(\CC)$
     (i.e., collections of perverse sheaves $\Fc_n$ on $\Sym^n(\CC)$ for all $n\geq 0$)
      with values in $\Vc$
     which are:
     
     \begin{itemize}
     \item[(2a)]  Smooth with respect to the   diagonal stratification,  and
     
     \item[(2b)]  {\em Factorizable}, that is, are equipped with a natural multiplicative structure with respect to the natural
     semigroup structure on $\Sym(\CC)$, see Definition \ref{def:FPS}. 
     
     \end{itemize}
     \end{itemize}
     We further prove, in Theorem  \ref{thm:shuffle},  that the functor $L$ takes the tensor, cotensor and Nichols algebras
     to the perverse sheaves obtained as $j_!$, $j_*$ and $j_{!*}$ as explained above in \S \ref {par:first-goal}
     
     \vskip.2cm
     
     As far as the second goal goes, we use the skeletal philosophy to connect, directly,  bialgebras with ``bisheaves''
     (and to motivate the idea of bisheaves in the first place). 
     We also use it  to give a very simple proof of the result of Takeuchi (quoted in Proposition \ref {prop:takeuchi}) on the dimension
    of the linear algebra data (cohomology with support in $\RR^n$) of a perverse sheaf on $\CC^n$ smooth with respect
    to the complexification of a real hyperplane arrangement. 
     
     \paragraph{Summary of the paper.} 
     In Ch.  \ref{sec:skeleta}  we discuss, at some length, our philosophy of skeleta and bisheaves (\S\S  \ref{subsec:skel-pure}
    and \ref{subsec:sh-bish}).  In \S \ref{subsec:takeuchi}, we consider examples of representation of skeleta
    as limit positions of deformations of ``real fronts'' $\{\Re(f)\geq 0\}$ of holomorphic functions $f$
    and give a  simple proof of    Takeuchi's theorem (Prop. \ref {prop:takeuchi}).  
    In \S  \ref{perv-ab-cat}  we give  
         a (rather routine) formal treatment of perverse sheaves with values in an abelian category, the concept which will be used later in the paper. 
         Here it is convenient
         to use the sub-analytic topology of Kashiwara-Schapira \cite {KS-ind}.

    \vskip.2cm
    
    In Ch. \ref {subsec:bish-bial} we introduce our main object of study: the symmetric product $\Sym^n(\CC)$
    together with its natural diagonal stratification.  In particular,  we ``break the symmetry''  of purely complex geometry,
    by introducing, in \S \ref{subsec:im-fine}, the {\em imaginary stratification} and the  decomposition into {\em fine cells}
    (extension of the Fox-Neuwirth cells  decomposition of $\Sym^n_\neq(\CC)$). 
     This allows us to construct, for any perverse sheaf $\Fc$ smooth with respect to the diagonal stratification,
     a natural Cousin resolution, an explicir complex of sheaves $\Ec^\bullet(\Fc)$ quasi-isomorphic to $\Fc$,
     see Proposition \ref{prop:cousin}. 
     In \S \ref{subsec:bar-bialg} we introduce a natural class of examples of ``bisheaves'' on the real skeleton:
     the data encoding both the bar and cobar complex of a braided bialgebra. In particular, we introduce,
     in Definition \ref{def:prim-gen-bialg}
     a concept of   primitive bialgebras.   Such bialgebras will be later (Theorem  \ref{thm:localize})
     identified with factorizable perverse sheaves. 
     
     \vskip .2cm
     
     In Ch. \ref{sec:fact-perp-bi} we, first of all, recall the formalism of shuffle algebras in the generality we need
     (\S \ref{subsec:shuf-braid}). In \S \ref{subsec:fact-perv} we introduce the concept of factorizable
     perverse sheaves with values in a braided category. For this it is convenient to use the definition
     of a braided tensor structure based on the $E_2$-operad.  In \S \ref{par:loc-fun} we formulate our main
     results: the Localization Theorem \ref{thm:localize} on equivalence between primitve bialgebras and
     factorizable perverse sheaves and Theorem \ref{thm:shuffle} identifying perverse sheaves corresponding
     to shuffle and Nichols algebras. 
     
     \vskip .2cm
     
     Ch. \ref{sec:cousin} is devoted to the proof of Theorems  \ref{thm:localize} and  \ref{thm:shuffle} 
     by a very direct method. Starting from a primitive bialgebra $A$, we produce, for each $n$, an explicit complex
     of sheaves $\Ec_n^\bullet(A)$ on $\Sym^n(\CC)$. A priori, it is only $\RR$-constructible  with respect to
     the stratification into fine cells, but we prove that it is in fact, a factorizable perverse sheaf with respect
     to the complex (diagonal) stratification.
     This remarkable phenomenon (the stalks of the complex $\Ec_n^\bullet(A)$ jump at real ``walls'' while their cohomology does not) comes from contractibility of some combinatorial complexes which appear as  multiplicities and
     which are treated systematically in Appendix \ref {app:shuff-com}. 
     
     \vskip .2cm
     
      In Ch. \ref {sec:examples}, which does not depend on the rest of the paper, we illustrate Theorem  \ref{thm:localize}
      using well known descriptions of perverse sheaves on $\Sym^n(\CC)$ for small values of $n$:
      the classical $(\Phi,\Psi)$-description for $n=2$ and the Granger-Maisonobe description \cite {GM-cusp} for $n=3$. 
      
         \vskip .2cm
         
         Finally,   Appendix  \ref{app:shuff-com}
         discusses various incarnations of a classical combinatorial object, the {\em shuffle complex},
         whose vertices are shuffles of some type $(d_1,\cdots, d_m)$, edges correspond to elementary moves etc. 
         Complexes of this nature appear as multiplicities in the stalks of the Cousin complexes in \S  \ref{sec:cousin}.

\paragraph{Relation to previous work.}\label{par:rel-prev}

 In \cite{BFS} (inspired in turn by \cite{SV}), 
intermediate extensions of local systems on  configuration spaces have been used 
 to construct the Hochschild complex for the positive 
part of a small quantum group  with coefficients in a tensor product of irreducible modules. 
In particular, one gets in this way the weight components of the irreducible modules themselves.  

To achieve this goal one works at first on {\em Cartesian} configuration spaces
(i.e., open parts of Cartesian, not symmetric products of $\CC$) and then 
passes to invariants with respect to actions of the appropriate symmetric groups.

\vskip .2cm

Unlike \cite{BFS} (and earlier precursors), we work directly on the symmetric products $\Sym^n(\CC)$
rather than Cartesian powers $\CC^n$ with diagonal stratifications.  The advantage of this approach is that
the linear algebra data have much smaller size (for $n=2$ we get the classical $(\Phi,\Psi)$-description)
 and match the features of a graded Hopf algebra very directly.  

\vskip .2cm

 We have been very much influenced by the unpublished preprint of D. Gaitsgory 
\cite{gaitsgory-fact}  on factorizable algebras
on configuration spaces.   The method of  \cite{gaitsgory-fact} (see also \cite{lysenko}, 
Theorem 29.2.3)  is based on the Koszul duality relation
between bialgebras and $E_2$-algebras, a purely ``derived'' concept which goes back to the works of Tamarkin and Kadeishvili, see    \cite{tamarkin-form},   
\cite{tamarkin-quant}, \cite{kadeishvili}, \cite{costello-yang},  
 and on the correspondence between $E_2$-algebras and factorizable constructible complexes on Ran spaces due to Lurie
 \cite{lurie:ha}.

  \vskip .2cm
  
 Unlike  \cite{gaitsgory-fact}, our approach does not use $E_2$-algebras (even though the $E_2$-operad
 provides a natural language for defining braided categories). Instead,  we
 produce the factorizable perverse sheaf $L(A)$ associated to
  a bialgebra $A$
 in a completely direct and explicit way:  as the sequence of  Cousin complexes,
 using the $2n$-dimensional real geometry of the  ``imaginary strata''
  in $\Sym^n(\CC)$. In fact, this geometry may be useful for constructing   manageable
  models for the dg-operad ${\bf e}_2$ (singular chains of the topological operad  $E_2$).  
  
    \vskip .2cm
    
    More recently, Ellenberg, Tran and Westerland \cite{westerland}, working in a general braided context,
    identified the Ext-cohomology of $T_*(V)$ with the homology of the braid group $\Br_n$ with coefficients
    in $V^{\otimes n}$ (twisted by the sign character). This result (Corollary 4.5 of \cite{westerland}) is recovered,
    in an equivalent form, as our Corollary \ref {cor:braid-cohom}.  An example of particular interest for \cite{westerland}
    is $\Vc$ being the category of Yetter-Drinfeld modules for a finite group $G$ and $V$ being associated to
    a conjugation invariant subset $C\subset G$. In this case $V^{\otimes n}$ is the permutation representation
    associated to the natural action of $\Br_n$ on $C^n$, and its homology is identified with the topological
    homology of the {\em Hurwitz space} $\on{Hur}^C_{G,n}$, the unramified cover of $\Sym^n_\neq (\CC)$
    associated to the $\Br_n$-set $C^n$, see \cite {ellenberg-venkatesh, westerland} for background. 
    Our approach puts these results into perspective and suggests a similar relation  of the Nichols
    algebra $T_{!*}(V)$ (for $V$ associated to $G,C$) to  intersection homology of some compactifications of
     $\on{Hur}^C_{G,n}$,

 \paragraph{Future directions.} Our approach, being completely elementary, suggests various generalizations,
 of which we mention two.
 
 \vskip .2cm

 (1) In this paper we consider only non-linear objects:
  graded Hopf algebras $A$. One can
 similarly describe, in terms of perverse sheaves, the natural linear objects:  {\em  Yetter-Drinfeld modules} over $A$
 (which form an abelian category) \cite{yetter} \cite{bespalov}. For this, we need to use
 the same symmetric products $\Sym^n(\CC)$ but with different stratifications:  choose a
 distinguished point $x_0\in\CC$, say $x_0=0$, and stratify $\Sym^n(\CC)$ according to the multiplicity of $x_0$
 (a separate non-negative integer) and the collection of multiplicities at other points (an un-ordered partition).
 The natural real skeleton corresponding to this stratification is $\Sym^n(\RR_{\geq 0})$.  
 After passing, as in \cite{gaitsgory-fact},  to $E_2$-algebras Koszul dual to bialgebras,  this fits into the package of 
 factorization algebras on $\CC$ stratified by several points,
 see \cite {ayala, CG, ginot}.  To  further extend from $\CC$ to more general Riemann surfaces, one needs to
 work in a {\em ribbon} braided category \cite{bakalov}.

 \vskip .2cm
 
 (2) The main context
  of the Feigin-Odesskii definition \cite{feigin-odesskii} of shuffle algebras (and of their many applications ever since)
 was  that of  $R$-matrices which  depend meromorphically on a parameter $t\in C$ (with $C$ typically an algebraic
 curve).   Conceptually,
 this corresponds to working in a {\em meromorphic braided category} \cite{soibelman}, where one has
 some algebro-geometric structure on the ``moduli space'' of objects.  It was pointed out by 
 Kontsevich and Soibelman  \cite{KoSo}  that the parametric shuffle algebra gives rise to some factorization data
 on the cartesian products of $C$.  The definition of these factorization data in \cite{KoSo} 
is similar to the  definition from \cite{BFS}, but 
the data in \cite{KoSo}  
 are of  ``coherent'' nature (i.e., involve coherent sheaves). In the case of constant $R$-matrices,
  considered in the present paper,  we deal with  ``constructible'' data (perverse sheaves). 
  It seems very plausible that these two directions: constructible and coherent, can be unified to
  give a conceptual interpretation of parametric shuffle algebras.

\paragraph{Acknowledgements.} It is a pleasure to dedicate this work to Kyoji Saito. 
His proposal was at the origin of the Granger - Maisonobe example \cite{GM-cusp}, which, in turn, 
was at the origin of our reflections on these topics. 

\bigskip

We are grateful to D. Gaitsgory, V. Hinich,   P. Safronov, K. Saito and C. Westerland for useful discussions. 
The research of M.K. was supported by World Premier International Research Center Initiative (WPI Initiative), 
 MEXT, Japan and by the IAS School of Mathematics. V.S. is grateful to Kavli IPMU 
for the hospitality during the conference on vertex algebras in the summer 2018, where 
a part of the present work has been reported. 

 
 \section{Perverse sheaves and  ``bisheaves'' on skeleta}\label{sec:skeleta}
 
 \subsection{Skeleta pure with respect to a complex stratification}\label{subsec:skel-pure}
 
 \paragraph{Pure skeleta. }\label{par:skel-pure}

 We fix, once and for all, a base field $\k$ of characteristic $0$. We denote $\Vect_\k$, resp. $\Vect_\k^\fd$
 the category of all, resp. of finite-dimensional, $\k$-vector spaces. 
 
  \vskip .2cm

  Let $X$ be a connected complex manifold, $\dim_\CC(X)=n$. Let $S=\{X_\alpha\}_{\alpha\in A}$ be a 
 complex analytic Whitney stratification of $X$. We denote by $i_\alpha: X_\alpha\to X$ the
 embedding of the stratum $X_\alpha$. By $X_0$ we denote the unique open stratum of $S$, so $0\in A$
 is a distinguished element.

 \vskip .2cm

 Let $\Sh(X,S)$ be the abelian category of sheaves of $\k$-vector spaces $\Fc$ on $X$ constructible with respect to $S$.
 This means $\Fc$  is locally constant on each $X_\alpha$, with finite-dimensional stalks.
 Let further
 $D^b(X,S)$ be the bounded derived category formed by  complexes of sheaves $\Fc$ of $\CC$-vector spaces on
 $X$   constructible with respect to $S$. This means 
  each cohomology sheaf $\ul H^i(\Fc)$  belongs to $\Sh(X,S)$.    
  
    \vskip .2cm

 Let $\Perv(X,S)\subset D^b(X,S)$ be the abelian category formed by  perverse sheaves  
 with respect to the middle perversity. Explicitly, see, e.g., \cite{KS-sheaves}  a complex $\Fc$ is perverse, if:
 \begin{itemize}
 \item[($P^+$)] For every $\alpha\in A $, we have $\ul H^n(i_\alpha^*\Fc) = 0$ for $n > -\dim_\CC(X_\alpha)$. 
	\item[($P^-$)] For every $\alpha\in A $, we have $\ul H^n(i_\alpha^! \Fc) = 0$ for $n < -\dim_\CC (X_\alpha)$. 
 \end{itemize}
  In particular, the conditions imply that the restriction of a perverse $\Fc$ to the open stratum $X_0$ is
 a local sysytem in degree $(-n)$. 
  The category $\Perv(X,S)$ carries the Verdier duality $\Fc\mapsto\Fc^*$. 
  
  We will also denote 
  \[
  D^b_\constr(X) = \bigcup_S D^b(X,S),\quad \Perv(X)=\bigcup_S \Perv(X,S)\subset D^b_\constr(X)
  \]
 the derived category  
  of complexes of sheaves on $X$ constructible with respect
  to some (indeterminate) stratification  and the  subcategory of perverse sheaves there.

    \vskip .2cm
    
    Let $K\subset X$ be a closed  real subanalytic set of pure real dimension $n$.  We denote $i_K: K\to X$
    the embedding. 
    The  Whitney stratification $S$ induces
    a stratification (disjoint decomposition) $S_K$ of $K$ into strata $K_\alpha=K\cap X_\alpha$. 
    
    By a {\em topological ball}  we will mean a subanalytic subset isomorphic, as such, to $\RR^m$ for some $m$.
    
    \begin{Defi}
    (a) $K$ is called a {\em skeleton} for $(X,S)$, if each $K_\alpha$ is a disjoint union of topological balls
    of the same dimension.
    
    \vskip .2cm
    
    (b) A skeleton $K$ is called {\em pure}, if for any $\Fc\in \Perv(X,S)$ the complex $i_K^!\Fc$
    on $K$ is quasi-isomorphic to a single sheaf in degree $0$, which we denote $R_K(\Fc) = \ul\HH_K^0 (\Fc)$.  
    \end{Defi}
    
    We refer to the connected components of $K_\alpha$ as {\em cells} of $K$ and denote by $(\Cc, \leq)$
    the set of cells of $K$, ordered by inclusion of closures. 
    
       \vskip .2cm
    
   Let $\Sh(K, \Cc)$ be the abelian category of sheaves of $\k$-vector spaces on $K$ which
   are locally constant (hence constant) on each cell of $\Cc$, with finite-dimensional stalks.
 
      \vskip .2cm

   If $K$ is a pure skeleton, we have an exact functor of abelian categories
   \[
   R_K = i_K^!: \Perv(X,S)\lra \Sh(K,\Cc). 
   \]
   In particular, for any cell $C\in\Cc$ we have an exact functor
   \[
   R_{K,C}: \Perv(X,S)\lra\Vect_\k^\fd, \quad \Fc\mapsto R_K(\Fc)_C := H^0(C, R_K(\Fc))
   \]
   obtained by taking the stalk of $R_K(\Fc)$ at $C$ or, what is the same, at 
   any point $x\in C$ (these stalks are canonically
   identified for different $x\in C$). 
   
   \paragraph{Pure skeleta and quiver descriptions of $\Perv(X,S)$.}
   In studying perverse sheaves, considerable effort had been spend on ``elementary"
   descriprions of the categories $\Perv(X,S)$ for some particular stratified spaces $(X,S)$, see 
   \cite {beil-gluing, GM-cusp, galligo-GM, MV}  .
   Typically, such a description has the form of an identification
   \be\label{eq:perv-rep}
   \Perv(X,S)\buildrel \simeq \over\lra \Rep(Q,R)
   \ee
   with 
  $\Rep(Q, R)$, the category of
   representations of some quiver (oriented graph) $Q$ subject to a set of relations $R$,
   a category which is manifestly abelian.  An object of $\Rep(Q,R)$
   is a family of vector spaces $V_i\in \Vect$ (corresponding to the vertices of  $Q$) and
   linear operators $a_e: V_{\on{source}(e)} \to V_{\on{target}(e)}$ for any oriented edge $e$ of $Q$,
   subject to the relations from the set $R$. 
   
   An identification \eqref{eq:perv-rep} contains, in particular, a datum of several exact functors
   \[
   \Perv(X,S)\lra\Vect_\k^\fd,\quad \Fc\mapsto V_i(\Fc), 
   \]
    where the  $V_i(\Fc)$ are the vector spaces of the quiver
   representation corresponding to $\Fc$. Now, a priori it is not easy to construct such exact functors
   in the first place. But if $K$ is a pure skeleton for $(X,S)$, then we have a natural
   supply of such functors, namely the $R_{K,C}, C\in\Cc$.

   \begin{observ}
   In most examples, the $V_i(\Fc)$ can indeed be found as $R_{K,C}(\Fc)$ for an appropriate skeleton $K$.
   Further, the natural maps between the $R_{K,C}(\Fc)$ coming  from the structure of  a sheaf on $R_K(\Fc)$,
   provide {\em precisely one-half} of the arrows $a_e$ in the quiver description, the other half being 
  the  ``formal adjoints" of these. The key to a quiver description consists therefore in finding a sufficiently
  rich pure skeleton for $(X,S)$. 
   \end{observ} 
   
   \begin{ex}[(Graphs on surfaces)]
   Let $n=1$, i.e., $X$ be a Riemann surface. A stratification $S$ consists then of a discrete set $N\subset X$
   so that the strata are $X\setminus N$ and the points of $N$. A skeleton for $(X,S)$ is nothing but
   a (subanalytic) graph $K$ embedded into $X$, regardless of its position with respect to $N$. 
   In this case it is not hard to prove that every graph is in fact a pure skeleton, see \cite{KS-graphs}, Prop. 2.2. 
   By taking $K$ sufficiently rich (i.e., passing through all the points of $N$ and being a spanning graph for $X$),
   one can obtain from this approach a quiver description of $\Perv(X,N)$, see \cite{KS-graphs}. 
   This contains the following more well known examples.
   \end{ex}
   
   \begin{ex}[(The $(\Phi, \Psi)$ description)]\label{ex:phi-psi}
   Let $X=\CC$ and let $N=\{0\}$ consist of one point. 
   Denote the corresponding category of perverse sheaves $\Perv(\CC,0)$. 
   Take for $K$ the closed half-line $[0,\infty)$ through $0$. The cells of $K$ are $\{0\}$ and $(0,\infty)$.
    Given $\Fc\in\Perv(\CC,0)$, the  stalks of the sheaf $R_{[0,\infty)}(\Fc)=\ul\HH^0_{[0,\infty)}(\Fc)$ at these cells
    are $\Phi(\Fc)$ and $\Psi(\Fc)$, the classical spaces of {\em vanishing and nearby cycles} of $\Fc$. 
    
    Further, a constructible sheaf in $\Sh([0,\infty) ,\Cc)$ is  the same as a diagram $v: E_0 \to E_1$, where $E_0$
   is the stalk at $0$, $E_1$ the stalk at $(0,\infty)$ and $v$ is the generalization map. So the sheaf structure
   of $R_{[0,\infty)}(\Fc)$ associates to $\Fc$ the arrow $\Phi(\Fc)\buildrel v\over\to\Psi(\Fc)$ which is
   precisely one half of the two arrows in the classical description \cite{galligo-GM}
\cite{beil-gluing} of $\Perv(\CC,0)$ as the category of diagrams
      \[
      \xymatrix{
  \Phi
 \ar@<.7ex>[r]^{v}&\Psi
\ar@<.7ex>[l]^{u}   
}, \quad T_\Psi:=  \Id_\Psi-vu , \,\, T_\Phi= \Id_\Phi-uv \,\,\text{ invertible}. 
\]
We note that invertibility of $T_\Psi$ implies invertibility of $T_\Phi$ and vice versa (even without assuming finite dimensionality of
$\Phi$ and $\Psi$), because of the identity\be\label{eq:(1-uv)}
(1 - uv)^{-1} = 1 + u(1 - vu)^{-1}v,
\ee
 familiar in the theory of the Jacobson radical, cf. \cite{lambek} \S 3.2, proof of Proposition 4. 
 We are grateful to V. Drinfeld for pointing it  out  to us. 
  \end{ex}
   
 \begin{ex}[(The ``Dirac'' description)]\label{ex:dirac}
  Keeping the same $X=\CC$ and $S= \bigl\{ \{0\}, \CC\setminus\{0\}\bigr\}$ as above, we take  $K=\RR=(-\infty, \infty)$ the real line through $0$. This $K$ has three cells:
 $\RR_- = (-\infty, 0)$, $\{0\}$ and $\RR_+=(0,\infty)$. 
 A constructible sheaf in $\Sh(\RR, \Cc)$ is the same as a diagram of stalks at these cells and the generalization maps
  \[
      \xymatrix{
E_- & E_0
\ar[l]_{\gamma_-}
\ar[r]^{\gamma_+}& E_+
 }. 
\]
This corresponds to the so-caled ``Dirac description" of $\Perv(\CC,\Fc)$, see \cite{KS-hyp-arr} \S 9A,  in terms of the diagrams  
\[
\begin{gathered}
      \xymatrix{
E_-\ar@<-.7ex>[r]_{\delta_-}& E_0
\ar@<-.7ex>[l]_{\gamma_-}
\ar@<.7ex>[r]^{\gamma_+}& E_+, 
\ar@<.7ex>[l]^{\delta_+}   
}
\\
\gamma_-\delta_- = \on{Id}_{E_-}, \,\,\,\gamma_+
\delta_+ = \on{Id}_{E_+}, \quad 
 \gamma_-\delta_+: E_+\to E_-, \,\,\, \gamma_+\delta_-: E_-\to E_+ 
 \text{ 
 invertible.}
\end{gathered}
\]
Again, the structure of the sheaf $R_{\RR}(\Fc)$ for $\Fc\in\Perv(\CC,0)$ provides precisely one-half
of the arrows in the diagram. 
  \end{ex}
  
     \begin{ex}[(Hyperplane arrangements)]\label{ex:hyp-arr}
   More generally, let $X=\CC^n$ and let $\Hc = \{H_i\}_{i\in I}$ be an arrangement of linear hyperplanes in $\RR^n$. 
   Then $\Hc_\CC = \{H_i^\CC = H_i\otimes\CC)$ is an arrangement of hyperplanes in $\CC^n$ and it gives a Whitney
   stratifications $S=S_\Hc$  into the strata
   \[
   L^\circ \, = \, L\, \setminus\,   \bigcup_{H_i\not\supset L} (L\cap H_i^\CC).
   \]
   Here $L$ runs over  {\em flats} of $\Hc_\CC$, i.e., linear subspaces obtained as intersections
   $\bigcap_{j\in J} H_j^\CC$ for any subset $J\subset I$  (including $J=\emptyset$ which corresponds to
   $L=\CC^n$). 
   
   Denote the corresponding category of perverse sheaves $\Perv(\CC^n, \Hc)$. It was studied
   in \cite{KS-hyp-arr}. In particular, it was proved there that the real space $K=\RR^n\subset\CC^n=X$
   is a pure skeleton for this stratification. 
   The set $\Cc$ consists of {\em faces} (of all dimensions)  of the real arrangement $\Hc$. 
  \end{ex}

   We recall, see, e.g.,  \cite{KS-hyp-arr} \S 1D:
   
   \begin{Defi}\label{def:quasi-reg}
   A {\em regular cell decomposition} of a space $X$ is a stratification
   into open topological balls (cells) such that the closure of any cell is a closed topological ball. 
   A {\em quasi-regular cell
   complex} is a difference $X\setminus X'$   of two regular cell complexes with $X'$ a  closed cellular subcomplex of $X$. 
     \end{Defi}
   
   In the examples above (in fact, in most of the examples of skeleta we know),
   the stratification $S_K$ is a  quasi-regular cell decomposition.
     For a quasi-regular cell complex,  it is standard that taking the stalks at various
  cells and generalization maps among them gives an identification
    \be
   \Sh(K,\Cc) \,\simeq \,  \Rep(\Cc) \, := \, \Fun((\Cc, \leq), \Vect_\k^\fd). 
   \ee
  Here $\Rep(\Cc)$ is the category of representations of the poset $(\Cc, \leq)$, i.e., functors from
  $(\Cc, \leq)$ considered as a category, to $\Vect_\k^\fd$. 
   
   For more general cell decompositions the role of $(\Cc, \leq)$ is played by
   $\Exit(K,\Cc)$, the  {\em  category of exit paths} of the stratification $\Cc$, see \cite{treumann}. 
   
   \vfill\eject
   
   \subsection{Sheaves and bisheaves}\label{subsec:sh-bish}

  \paragraph{The idea of bisheaves. }
   In the examples we have seen, the linear algebra data describing $\Fc\in\Perv(X,S)$
   give, in particular,  an actual constructible sheaf on the skeleton $K$, in terms of its stalks and generalization maps.
   But they also contain companion maps in the opposite direction which, similarly, give a {\em cosheaf}
   on $K$, that is, a sheaf with values in $\Vect_\k^\op$, the opposite to the category of vector spaces. 
   
   It is suggestive, therefore, to introduce the term  {\em bisheaf}
  to describe the entire set of linear algebra data representing $\Fc$. Informally, a bisheaf should consist
  of an $S_K$-constructible  sheaf and cosheaf on $K$, with the same stalks at the cells, satisfying some
  compatibility conditions between the two structures. 
  
  \begin{ex} For instance, in the situation of Example
  \ref {ex:dirac} we have $K=\RR$ with stratification $K_S$ into $\RR_-$, $\{0\}$ and $\RR_+$
   A  $S_K$-constructible sheaf on $K$ is the same as a diagram
   \[
      \xymatrix{
E_- & E_0
\ar[l]_{\gamma_-}
\ar[r]^{\gamma_+}& E_+,
 }
\]
a cosheaf is the same as a diagram
\[
 \xymatrix{
E_-\ar [r]^{\delta_-}& E_0
 & E_+, 
\ar [l]_{\delta_+}   
}
 \]  
and a bisheaf should be a diagram 
   \[
  \xymatrix{
E_-\ar@<-.7ex>[r]_{\delta_-}& E_0
\ar@<-.7ex>[l]_{\gamma_-}
\ar@<.7ex>[r]^{\gamma_+}& E_+, 
\ar@<.7ex>[l]^{\delta_+}   
}
 \]  
 satisfying the conditions $\gamma_\pm\delta_\pm = \Id$ as well as the invertibility of $\gamma_+\delta_-$ and
 $\gamma_-\delta_+$. 
 \end{ex}
 
 We do not know a general definition of a bisheaf. It seems that it should not be a purely topological
 concept involving  the stratified space $(K, S_K)$ but should involve information about the
 complexification $(X,S)$ as well. 
 
 \paragraph{Self-dual skeleta.}  The skeleta that have been used, have, in addition to purity, the following property:
 the functors $R_{K,C}: \Perv(X,S)\to\Vect_\k^\fd$, take Verdier duality to the usual duality fo vector spaces, i.e., we have
 functorial isomorphisms
 \[
 R_{K}(\DD\Fc )_C \,\simeq \,  (R_K(\Fc)_C)^*.
 \]
This duality then interchanges the sheaf and cosheaf structures. 
That is, the additional ``cosheaf" structure on  the sheaf $R_K(\Fc)$
comes from the structure of sheaf on $R_K(\Fc^*)$. 

We will call this property of a pure skeleton {\em self-duality}. Using the standard properties of Verdier duality
(interchanging  $i^*$ and $i^!$),  for a self-dual skeleton $K$ we have
natural isomorphisms
\[
i_C^* i_K^! \Fc \,\simeq i_C^! i_K^*\Fc, \quad \Fc\in \Perv(S,X), \,\, C\in\Cc. 
\]

\vfill\eject

\subsection {Pure skeleta as degenerations of real fronts. Takeuchi's theorem}\label{subsec:takeuchi}

Let $f$ be a holomorphic function on $X$. We then have the hypersurface $X_f=\{f=0\}$
and the closed subanalytic set $X^+_f = \{\Re(f)\geq 0\}$ which we call the {\em real front} of $f$.  Denote
\[
X_f\buildrel i_f\over\hookrightarrow X^+_f \buildrel j_f\over\hookrightarrow X
\]
 the embeddings. It is known  \cite{KS-sheaves} that the functor of vanishing cycles
 \[
 \Phi_f: D^b_\constr(X) \lra D^b_\constr(X_f)
 \]
 perserves perverse sheaves and commutes with Verdier duality. It is also known ({\em ibid.})
 that we have a canonical quasi-isomorphism
 \[
\Phi_f(\Fc) \,\simeq \,  i_f^* j_f^!(\Fc). 
 \]
 On the other hand, for a skeleton $K$ and its cell $C$,  the stalk at $C$ of the complex $R_K(\Fc)$ can be
 written as $i_C^* i_K^! (\Fc)$. This suggests a way to construct pure and self-dual skeleta as limits of
 isotopic deformations of regions of the form $X_f^+$. 
 
 \begin{ex}\label{ex;phi-psi-def}
 (a) In the situation of Example \ref{ex:phi-psi}, the skeleton $[0,\infty)$ can be seen as a limit of
 isotopic deformations (compatible with our stratification) of the real front $\Re(z)\geq 0$.
 
  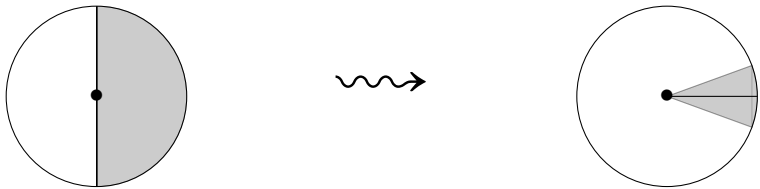
\begin{figure}[ht]
   \centering
    
\begin{tikzpicture}[scale=0.3]

\def\centerarc[#1](#2)(#3:#4:#5)
    { \draw[#1] ($(#2)+({#5*cos(#3)},{#5*sin(#3)})$) arc (#3:#4:#5); }
    
\node at (0,0){$\bullet$};

\draw (0,0) circle (4); 
  \clip (0,-4) rectangle (4,4);
    \filldraw[opacity=0.2] (0,0) circle(4);
    \draw [line width = 0.9] (0,-4) -- (0,4);

\end{tikzpicture}
\quad\quad\quad\quad
\begin{tikzpicture}[scale=0.3]
\node at (0,0){\huge$\rightsquigarrow$};
\node at (0,-4){};
\end{tikzpicture}
\quad\quad\quad\quad
 \begin{tikzpicture}[scale=0.3]
 
 \node at (0,0){$\bullet$};

\draw (0,0) circle (4); 
\draw (0,0) -- (4,0); 

\node (a) at (340:4){};
\node (b) at (20:4){};
\filldraw[opacity=0.2]  (0,0) -- (a.center)
  (a.center) arc(340:380:4)
  (b.center) -- (0,0); 
  \filldraw [opacity=0.2] (0,0) -- (a.center) -- (b.center) -- (0,0); 
  
 \end{tikzpicture}
 \caption{Deforming $\{\Re(z)\geq 0\}$  into the half-line. }\label{fig:half-line}
\end{figure}

\end{ex}

\begin{ex}\label{ex:dirac-def}
 In the situation of Example \ref{ex:dirac}, the skeleton $\RR=(-\infty, \infty)$ can be seen as a limit of isotopic
deformations of the real front $\{\Re(z^2)\geq 0\}$. The deformation is achieved by the transformations
$T_\eps: \CC\to \CC$, $0\leq \eps\leq 1$,  taking $x+iy$ to $x+\eps iy$, see Fig. \ref{fig:full-line}. Note that such transformations
cannot be made holomorphic, since we need to decrease the angle of aperture of the sectors. 
Nevetheless, they preserve the topology of the situation. 

 \begin{figure}[ht]
   \centering
    
\begin{tikzpicture}[scale=0.3]

\node at (0,0){$\bullet$};

\node (a-) at (225:4){}; 
\node (a+) at (45:4){}; 
\node (b-) at (315:4){}; 
\node (b+) at (135:4){}; 

\draw (0,0) circle (4); 

\filldraw[opacity=0.2]  (0,0) -- (b-.center)
  (b-.center) arc(315:405:4)
  (a+.center) -- (0,0); 
  \filldraw [opacity=0.2] (0,0) -- (b-.center) -- (a+.center) -- (0,0); 
  
  \filldraw[opacity=0.2]  (0,0) -- (b+.center)
  (b+.center) arc (135:225:4)
  (a-.center) -- (0,0); 
  
    \filldraw[opacity=0.2] (0,0) -- (b+.center) -- (a-.center) -- (0,0);

\end{tikzpicture}
\quad\quad\quad\quad
\begin{tikzpicture}[scale=0.3]
\node at (0,0){\huge$\rightsquigarrow$};
\node at (0,-4){};
\end{tikzpicture}
\quad\quad\quad\quad
 \begin{tikzpicture}[scale=0.3]
 
 \node at (0,0){$\bullet$};

\draw (0,0) circle (4); 
\draw (-4,0) -- (4,0); 

\node (b-) at (340:4){};
\node (a+) at (20:4){};
\filldraw[opacity=0.2]  (0,0) -- (a.center)
  (a.center) arc(340:380:4)
  (b.center) -- (0,0); 
  \filldraw [opacity=0.2] (0,0) -- (b-.center) -- (a+.center) -- (0,0);

  \node (b+) at (160:4){}; 
  \node (a-) at (200:4){}; 
  
  \filldraw[opacity=0.2] (0,0) -- (b+.center)
  (b+.center) arc (160:200:4)
  (a-.center) -- (0,0); 
  
    \filldraw[opacity=0.2] 
    (0,0) -- (b+.center) -- (a-.center) -- (0,0);

 \end{tikzpicture}
 \caption{Deforming $\{\Re(z^2)\geq 0\}$  into the real line. }\label{fig:full-line}
\end{figure}
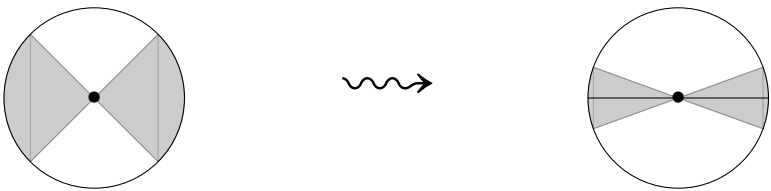
  \end{ex}
 
Let us  now consider the situation of Example \ref{ex:hyp-arr}, so we have  an arrangement $\Hc$ of hyperplanes in
 $\RR^n$ and  the corresponding 
  category $\Perv(\CC^n,\Hc)$. In the following proposition both parts are known but rather non-trivial.
  Part (a), already mentioned in Example \ref{ex:hyp-arr}, was proved in \cite{KS-hyp-arr} but  an equivalent
  statement goes back
  to \cite{BFS}. Part (b) is a result of Takeuchi \cite{takeuchi-Dmod}.

  \begin{prop}\label{prop:takeuchi}
  (a) The skeleton $\RR^n\subset\CC^n$ is pure and self-dual.
  
  \vskip .2cm
  
  (b) For $\Fc\in\Perv(\CC^n, \Hc)$ the dimension of the stalk $R_{\RR^n}(\Fc)_0$ is equal to
  \[
  \sum_{L} \on{mult}_{T^*_L\CC^n} (\Fc),
  \]
  where $L$ runs over all complex flats of $\Hc$ and $\on{mult}_{T^*_L\CC^n} (\Fc)$
  is the multiplicity of the Lagrangian variety $T^*_L\CC^n$ in the characteristic cycle of $\Fc$. 
  \end{prop}
  
  We now present an extremely simple proof of both statements by the method of degenerations,
   developing a higher-dimensional version of Example \ref{ex:dirac-def}. That is, we take a quadratic form 
  $Q(z)$ on $\CC^n$ with real coefficients such that the restriction of $Q$ to $\RR^n$ is positive definite. 
  Explicitly,  we choose a linear coordinate system $x_1,\cdots, x_n$ on $\RR^n$ and define 
  \[ 
  Q: \CC^n\to\CC, 
  \quad Q(z_1,\cdots, z_n) \,=\,\sum z_i^2.
  \] 
  The real part $\Re(Q): \CC^n\to\RR$ is a nondegenerate quadratic form of signature $(n,n)$,
  and the transformations
  \[
  T_\eps:\CC^n\to\CC^n, \quad {\bf x} + i {\bf y} \,\mapsto {\bf x} + i\eps{\bf y}, \quad 0\leq \eps\leq 1,
  \]
 ``flatten''  the region $\Re(Q) \geq 0$ into $\RR^n$ in a way compatible with the stratification
  $S_\Hc$, much like in Fig. 
 \ref{fig:full-line}.  Therefore
 \be\label{eq:hyperb-van}
 i_0^* i_{\RR^n}^! (\Fc) \,\simeq \, \Phi_Q(\Fc)_0
 \ee
  is the stalk at $0$ of the perverse sheaf of vanishing cycles $\Phi_Q(\Fc)$.
  If we now assume that $\Hc$ is a central arrangement i.e., that $\{0\}$ is the minimal flat
  (we can always reduce to this case by quotienting by the minimal flat), then $\Phi_Q(\Fc)$
  is supported at 0 and so reduces to a single vector space in degree $0$. This shows
  that the stalk of $i_{\RR^n}^!(\Fc)$ at $0$, the minimal celll of the real arrangement,
  is concentrated in degree $0$. The stalk at any other cell $C$ is analyzed similarly, by taking
  the transversal slice to $C$ and using a positive definite quadratic form on that slice. This proves (a).
  
  \vskip .2cm
  
  Further,  (b) follows at once \eqref{eq:hyperb-van} and the standard microlocal interpretation of
  vanishing cycles \cite{KS-sheaves}. That is, the cotangent space of $\CC^n$ is $\CC^n\times\CC^{n*}$.
  The stratification $S=\{L^\circ\}$ gives rise to
  the set of Langangian varieties $T^*_L\CC^n = L\times L^\perp$
   which can contribute to the characteristic cycle of $\Fc$.
  The graph of $dQ$ is the Lagrangian variety $\Gamma_{dQ}$
  which is the graph of the isomorphism $\CC^n\to\CC^{n*}$
  given by the bilinear form corresponding to $Q$. It meets each $T^*_L\CC^n$ exactly once, transversely, 
  at the point $(0,0)$.  Let us now deform $\Gamma_{dQ}$ a little,    replacing $Q$ with $Q'=Q+l$
  where $l$ is a generic linear function on $\CC^n$. Then
  \[
  \dim\Phi_Q(\Fc)_0 \,=\,\sum_{x\in\CC^n} \dim\Phi_{Q'}(\Fc)_x.
  \]
  But $\Gamma_dQ'$ is a generic shift of 
    $\Gamma'$ which will meet each $T^*_L\CC^n$ exactly once, transversely, at a point $(x_L,\xi_L)$  which does not
    belong to any other $T^*_{M}\CC^n$, $M\neq L$. 
  This implies that
 \[
 \dim \Phi_{Q'}(\Fc)_{x_L} \,=  \on{mult}_{T^*_L\CC^n}(\Fc)
 \]
 and the claim follows. \qed


 \subsection{Perverse sheaves with values in an abelian category}\label{perv-ab-cat}
 
 In the remainder of the paper it will be convenient to work with sheaves and perverse sheaves formed not by
 vector spaces but by objects of some abelian category $\Vc$. In this section we formally sketch this concept.
 The reader can assume that $\Vc$ consists of vector spaces with some extra structure (e.g., of modules over a $\k$-algebra $\Lambda$). 
 In such a case the generalization we describe is really trivial (perverse sheaves with  the same extra structure). 
 
 \paragraph{Sheaves with values in a Grothendieck category.}
  Let $\Vc$ be an $\k$-linear abelian category.  Note that there is an intrinsic concept of tensoring an object $V$ of $\Vc$ with
  a  finite-dimensional $\k$-vector space $E$, the result of which will be denoted $V\otimes_\k E$.

  \vskip .2cm 
  
  Let us first assume that
 $\Vc$ is a {\em Grothendieck category} \cite{KS-categories},
 in particular, that  it has arbitrary direct sums and products.  
 For example, one can take for $\Vc$ the category of {\em all}  modules over a $\k$-algebra $\Lambda$. 
 
 \vskip .2cm
  
 Let $X$ be a topological space. We can speak about {\em sheaves on $X$ with values in} $\Vc$. These are contravariant
 functors $\Fc$ from the poset of open sets in $X$ to $\Vc$  satisfying the descent axiom: for any open covering 
 $U=\bigcup_{i\in I} U_i$ of an open set $U$, the sequence
 \be\label{eq:descent}
 0\to\Fc(U) \lra \prod_{i\in I} \Fc(U_i) \lra \prod_{i,j\in I} \Fc(U_i\cap U_j) 
 \ee
 is exact.
 
  The operation of passing from a presheaf with values in $\Vc$ to the associated sheaf 
is  defined as usual. A constant sheaf with values in $\Vc$ is the sheaf 
associated with a constant presheaf $(U\subset X)\mapsto $ a fixed $V\in \Vc$. 
A {\it  local system on $X$ with values in $\Vc$} is a sheaf locally isomorphic 
to a constant sheaf. It is the same as a functor from the fundamental groupoid
$
\Pi_1(X) \to \Vc. 
$

\vskip .2cm
We denote by $\Sh(X,\Vc)$ the category of sheaves on $X$ with values in $\Vc$. It is again a Grothendieck abelian category. 
By $D^b(\Sh(X))$ we denote the bounded derived category of $\Sh(X,\Vc)$.

Further, let $(X,S)$ be a stratified complex manifold. Then, as in \S \ref{subsec:skel-pure}, we have:
\begin{itemize}
\item  The abelian category
$\Sh(X,S;\Vc)$ of $S$-constructible sheaves on $X$ with values in $\Vc$.

\item  The triangulated category
 $D^b(X,S;\Vc)$
of $S$-constructible complexes on $X$ with values in $\Vc$.

\item  The abelian subcategory 
$\Perv(X,S;\Vc ) \subset D^b(X,S; \Vc)$. 
\end{itemize}

\paragraph{Finiteness issues: subanalytic topology.}

 The requirement for $\Vc$ to be a Grothendieck category is, in general, forced by the
necessity to have infinite products as in \eqref{eq:descent}. This excludes, for instance, the case $\Vc=\Vect_\k^\fd$
or any other abelian category with a perfect duality. 

However, if we are interested only in constructible sheaves
and complexes, then the following more restricted framework is sufficient.

\vskip .2cm

Let $X$ be a real analytic manifold. We recall \cite{KS-sheaves}, \cite{KS-ind} the
concept of a {\em subanalytic set} in $X$. In particular, we have the 
 {\em subanalytic site} $X_\sa$ of $X$. 
This is the poset of  relatively compact  subanalytic open subsets $U\subset X$ with coverings being finite coverings in the usual sense. 
Since all coverings are finite, we can speak about sheaves on $X_\sa$ with values in any abelian category $\Vc$, not necessarily
Grothendieck one. 

 In particular, if $(X,S)$ is a stratified  complex manifold, then we have the categories
\[
\Sh(X_\sa,S;\Vc)\,\subset \, D^b(X_\sa ,S;\Vc)\, \supset \, \Perv(X_\sa,S;\Vc). 
\]
If $\Vc$ has a perfect duality $V\mapsto V^*$, then $D^b(X_\sa,S;\Vc)$ and $\Perv(X_\sa,S;\Vc )$  also have
a perfect duality which we call the  Verdier duality and 
denote $\Fc\mapsto \DD\Fc$.  For an arbitrary $\Vc$ we have the Verdier duality in the form of anti-equivalences
\be\label{eq:verdier-gener}
\DD: \Perv(X_\sa, S,\Vc) \to  \Perv(X_\sa, S, \Vc^\op), \quad \DD: D^b(X_\sa ,S;\Vc)\to D^b(X_\sa ,S;\Vc^\op). 
\ee

\begin{prop}
If $\Vc$ is a Grothendieck abelian category, then the restriction functors
\[
\Sh(X ,S;\Vc)\to \Sh(X_\sa,S;\Vc), \quad D^b(X ,S;\Vc)\to D^b(X_\sa ,S;\Vc), \quad  \Perv(X,S;\Vc)\to  \Perv(X_\sa,S;\Vc)
\]
are equivalences of categories. 
\end{prop}

\noindent{\sl Proof:} Same argument as in \cite{KS-ind}, \S 7.1. \qed

\paragraph{Conventions for the rest of the paper.}

(1) In the sequel we will drop the subscript ``$\sa$'' and simply speak about $S$-constructible sheaves
and complexes on $X$ with values in $\Vc$, while having in mind the subanalytic topology.
 In particular, we will use freely the standard functorialities  such as $f_*, f^*, f_!, f_!$
on constructible complexes on $X_\sa$, see \cite {KS-ind}, \S 7.1 for their construction. 
For a holomorphic function $f$, we  will also use the vanishing cycle functor $\Phi_f$ as well as the fact
that it preserves $\Vc$-valued  perverse sheaves and commutes with Verdier duality in the form 
\eqref{eq:verdier-gener}. 

\vskip .2cm

(2) We will freely use the (obvious) analogs of the concepts and results of \S\S \ref{subsec:skel-pure} -\ref{subsec:takeuchi} for
perverse sheaves with values in $\Vc$. In particular, we will use the concept of a pure skeleton, of a self-dual
skeleton (either when $\Vc$ has a perfect duality or with respect to \eqref{eq:verdier-gener}), as well as both the $(\Phi, \Psi)$-description and the Dirac descriprion
 of  $\Vc$-valued perverse sheaves on a disk in terms of diagrams formed by objects and morphisms of $\Vc$. 
 In the same vein, we will use Proposition \ref{prop:takeuchi}(a) in the $\Vc$-valued context.


 \section{Bisheaves and bialgebras: perverse sheaves on symmetric products} \label{subsec:bish-bial}

\subsection {The symmetric product $\Sym^n(\CC)$,  its complex strata and real skeleton}
\label{subsec:symprod-basics}

\paragraph{The complex strata.} We denote by $\Sigma_n$ the symmetric group of order $n$. 
Fix $n\geq 1$ and 
let $X=\Sym^n(\CC)=\CC^n/\Sigma_n$ be the $n$th symmetric product of $\CC$. As
 well known, it can be seen in either of two ways:
\begin{itemize}
\item[(1)] As the
space of effective divisors of degree $n$, i.e.,  of formal sums
\[
\sum \lambda_i x_i,\,\,\, \lambda_i\in\ZZ_{>0}, \,\,\, x_i\in\CC, \,\,\, \sum \lambda_i=n.
\]
\item[(2)] As the space of monic polynomials of degree $n$
\[
f(x) = x^n + a_{n-1} x^{n-1} +\cdots + a_0, \,\,\, a_i\in\CC. 
\]
\end{itemize}
Explicitly, $x_i$  in (1) is a root of $f(x)$ with multiplicity $\lambda_i$. 
The description (2) makes it clear that $\Sym^n(\CC)$ is isomorphic to $\CC^n$. 
The description (1) gives a natural Whitney stratification of $\Sym^n(\CC)$ by the type of divisors.
It is labelled by (unordered) {\em partitions} of $n$, i.e., by sequences
\[
\lambda = (\lambda_1 \geq \cdots \geq \lambda_p>0), \quad |\lambda| :=\sum \lambda_i = n. 
\]
Here we can think of $\lambda$ as a Young diagram with $n$ cells. 
We denote the set of such partitions (or Young diagrams) by $\Pc(n)$. 
For $\lambda\in\Pc(n)$ we
 denote $l(\lambda)=p$ the
number of parts of $\lambda$, i.e., the number of rows in the Young diagram. 

\vskip .1cm

We denote by $\leq$ the standard partial ordering on $\Pc(n)$ whereby $\lambda\leq\mu$ if and only if 
$\mu$ is obtained from $\lambda$ by moving some cells to the right and up. 
With respect to this order, the $1$-part partition $(n)$ is the minimal one while the partition $(1,\cdots, 1)$
is the maximal one. 

\vskip .2cm

By definition, the stratum $X_\lambda^\CC$
associated to $\lambda$ consists of divisors that can be represented as $\sum\lambda_i x_i$
with all $x_i$ distinct. For example, $X_{1,\cdots, 1}^\CC$ is the open stratum
(all points distinct), while $X_{n}^\CC\simeq \CC$
is the minimal stratum (all points coincide). The following is straightforward.

\begin{prop}
(a)  $X_\lambda^\CC$ is a smooth locally closed subvariety in $\Sym^n(\CC)$
of dimension  $l(\lambda)$. The collection $S=\{X_\lambda^\CC\}$ forms a Whitney stratification of $X=\Sym^n(\CC)$.

(b) Fo  $\lambda, \mu\in\Pc(n)$ we have $X_\lambda^\CC\subset \overline X_\mu^\CC$
if and only if $\lambda\leq\mu$. 
\qed
\end{prop}

 We will refer to $S$ as the {\em diagonal stratification} and use the notation $S_n$ to empahize the dependence on $n$,
 if needed. 
  We will be interested in the category $\Perv(\Sym^n(\CC), S)$. 
  
  \vskip .2cm

Note that the group structure on $\CC$ gives a morphism
\be\label{eq:sum-divisor}
\sigma: \Sym^n(\CC)\lra \CC, \quad \sum\lambda_i x_i\,\,\, \mapsto\,\,\,  \sum\nolimits^\CC \lambda_i x_i, 
\ee
where $\sum^\CC$ is the operation of addition in $\CC$. 
We denote $\Sym^n_0(\CC)$ the preimage of 0
under $\sigma$. It is identified with the space of monic polynomials $f(x)$ with $a_{n-1}=0$. 
Note that we have an identification
\[
\Sym^n(\CC) \,\,\simeq \,\,\CC\times \Sym^n_0(\CC). 
\]
It is given by the action of $\CC$ on $\Sym^n(\CC)$ by additive shifts
  $f(x) \mapsto f(x+c)$.  Each stratum $X_\lambda$ is the product of $\CC$ and the induced
  stratum in $\Sym^n_0(\CC)$.

 \paragraph{The real skeleton.} 
  Let $K=\Sym^n(\RR)\subset\Sym^n(\CC)$ be the $n$th symmetric product of $\RR$, i.e.,
the subspace formed by effective divisors consisting of real points. 
As in the complex case, we have the action of $\RR$ on $\Sym^n(\RR)$ given by adding the
same number to all the points of a divisor. This action induces identifications
  \be\label{eq:sym-r-n-0}
\begin{gathered}
\Sym^n(\RR) \,\,\simeq \,\,\RR\times \Sym^n_0(\RR),
\\
\Sym^n(\RR)/\RR  \,\,\simeq \,\,  \Sym^n_0(\RR). 
\end{gathered}
\ee
Alternatively, $\Sym^n(\RR)$ can be seen as the set of those monic polynomials whose all roots are real,
and $\Sym^n_0(\RR)$ as the set of monic polynomial with all roots real and summing to 0. 
So $\Sym^n(\RR)$ is a closed  ``curvilinear cone" in the Euclidean space $\RR^n$ formed by
 all monic polynomials with real coefficients. 
 
 \vskip .2cm
 
 Similarly, let $\RR^n_0\simeq \RR^{n-1}$ be the Euclidean space of all 
   of monic polynomials with real coefficients and with $a_1=0$. Then  $\Sym^n_0(\RR)$ is
   a  ``curvilinear cone" in $\RR^n_0$. In examples, it is convenient to visualize $\Sym^n_0(\RR)$
   instead of $\Sym^n(\RR)$,
   factoring out the translation symmetry.

\begin{figure}[H]
\centering
\begin{tikzpicture}[scale=0.5]
\centering
\draw plot  [ line width = 1, domain=-2:2, variable=\x] ({\x*\x}, {0.5*\x*\x*\x}); 

\node at (0,0){$\bullet$};

\fill [gray, opacity=0.3,  domain=-2:2, variable=\x]
  (4,9) -- (4,-9)
  -- plot  [ domain=-2:2, variable=\x] ({\x*\x}, {0.5*\x*\x*\x})
 -- cycle;
 
 \draw  
 (4,4) -- (4,-4) -- (-4,-4) -- (-4, 4) -- (4, 4); 
 
 \node at (-0.7,-0.7){\small$(3)$}; 
 
 \node at (1,2){\small$(2,1)$};
 \node at (1,-2){\small$(1,2)$}; 
 \node at (2.5,0){\small$(1,1,1)$};
 \node at (3.2,2){\small$K$}; 
 
\end{tikzpicture}
\caption{ The real skeleton $K=K_3$ in $\Sym^3_0(\CC)$ and its strata labelled by $\2^{3-1}$. }\label{fig:semicubic}
\end{figure}
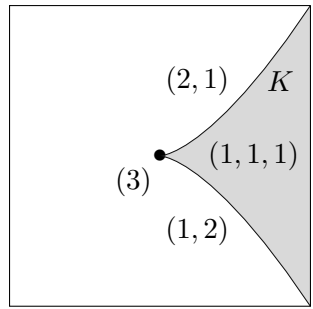

\begin{exas} (a) Let $n=2$. Then $\RR^2_0=\RR$ is the real line formed by polynomials
$x^2+a=0$, and $\Sym^2_0(\RR) = \{ a\leq 0\}$ is the 
negative real line. 

\vskip .2cm

(b) Let $n=3$. Then $\RR^3_0$ consists of real cubic polynomials  of the form $x^3+ax + b$, and 

\[
\Sym^3_0(\RR)\,\,= \,\,\bigl\{ (a,b) | \,\, 4a^3 + 27 b^2 \leq 0\bigr\}
\]
 is the interior of the  real  semi-cubical parabola, see Fig. \ref{fig:semicubic}. 
\end{exas}

We note that $K=\Sym^n(\RR)$ has
a natural stratification $S_K$  labelled by {\em  ordered partitions}    of $n$, 
i.e., sequences 
\be\label{eq:order-part} 
\alpha = (\alpha_1, \cdots, \alpha_p),\,\,\,\alpha_i\in \ZZ_{>0}, \,\,\,\sum \alpha_i=n. 
\ee
We denote $l(\alpha)=p$ the length of $\alpha$, and by $\ol\alpha = (\alpha_{\sigma(1)} \geq 
\cdots \geq \alpha_{\sigma(p)})$  (for appropriate $\sigma\in S_p$) the ordered partition corresponding to $\alpha$.

 The stratum
 $K_\alpha = \Sym^n_\alpha(\RR) $ corresponding to $\alpha$,  consists
of divisors of the form
\[
\alpha_1 x_1 +\cdots + \alpha_p x_p, \quad x_i\in\RR, \,\,\,  x_1 <\cdots < x_p.
\]
The number of ordered partitions of $n$ is $2^{n-1}$, they are in bijection with arbitrary subsets of the
set of  plus signs in the equality
\[
n= 1+1+1+\cdots + 1.
\]
Conceptually, these  plus signs correspond to simple roots of the root system $A_{n-1}$. 
To get $\alpha$ we simply perform the additions corresponding to the chosen subset of the plus signs. 
We therefore denote the set of ordered partitions of $n$ by $\2^{n-1}$
and equip it with the partial order of reverse inclusion of subsets. Geometrically, we view $\2^{n-1}$
as the set of vertices of an $(n-1)$-dimensional cube. 

Further, given $\alpha\in\2^{n-1}$, we denote by $\alpha^\vee$ the {\em complementary ordered partition},
for which the corresponding subset in $\{1,\cdots, n-1\}$ are complementary. Thus, for example,
$(1,1,1)^\vee= (3)$ etc. Note that the un-ordered partition $\ol{\alpha^\vee}$ corresponds to the transpose
of the Young diagram corresponding to the un-ordered partition $\alpha$. 

\begin{ex}
Let $n=3$. The poset $\2^{3-1}$, with order indicated by arrows,  has the form
 
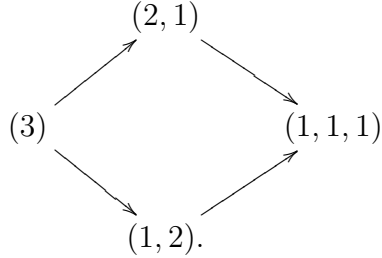
\[
\xymatrix{
&(2,1)\ar[dr]&
\\
(3)\ar[ur] \ar[dr] && (1,1,1)\\
&(1,2). \ar[ur]&
}
\]
 
\end{ex}

\begin{prop}\label{prop:K-strata}
For each $\alpha\in\2^{n-1}$, the stratum $K_\alpha$ is a cell of real dimension $l(\alpha)$.
 We have $K_\alpha \subset \ol K_\beta$ if and only
if $\alpha\leq\beta$ in $\2^{n-1}$. The image of $K_\alpha$ under the embedding
 $\Sym^n(\RR)\subset \Sym^n(\CC)$ lies in the stratum $X_{\ol\alpha}$. \qed
\end{prop}


\subsection{Imaginary strata and fine cells}\label{subsec:im-fine}

\paragraph{The imaginary strata of $\Sym^n(\CC)$.} 
 \label{par:im-cells}
The imaginary part map $\Im: \CC\to\RR$ gives a mapping
\[
\Ien: \Sym^n(\CC)\to\Sym^n(\RR) = K, \quad {\bf z} =  \sum \lambda_i z_i \mapsto \sum \lambda_i \, \Im(z_i). 
\]
The fibers of $\Ien$ are totally real subanalytic subspaces (not necessarily submanifolds) of $\Sym^n(\CC)$. 

We call the {\em imaginary stratification} of $\Sym^n(\CC)$ the pullback, under $\Ien$, of the stratification $(K_\alpha)$ of $K$ and
denote the corresponding strata, as well as their embeddings, by
\[
X_\alpha^\Ien \,=\, \Ien^{-1}(K_\alpha),\quad j_\alpha: X_\alpha^\Im \hookrightarrow X=\Sym^n(\CC),   \quad \alpha\in\2^{n-1}. 
\]

Explicitly, a point  (written as a diivisor) ${\bf z} \in \Sym^n(\CC)$ 
     lying in $X_\alpha^\Im$  splits (in a unique way) into a sum of
   effective divisors ${\bf z}_1 + \cdots + {\bf z}_p$ such that 
    \be\label{eq:S-a-split}
   \Ien({\bf z}_1) = \alpha_1 \cdot y_1, \,\,\,\cdots,\,\,\, \Ien({\bf z} _p) = \alpha_p \cdot y_p, \quad y_1 < \cdots < y_p, 
   \ee
   see Fig. \ref{fig:S-a-Im}.

 \begin{figure}[H]
 \centering
 \begin{tikzpicture}[scale=.4, baseline=(current bounding box.center)]

 \draw[->] (-2,0) -- (13,0); 
 
 \draw[->] (0,-2) -- (0,11); 
 
 \draw[dashed] (0,2) -- (11,2); 
 
  \draw[dashed] (0,4) -- (11,4); 
  
   \draw[dashed] (0,7) -- (11,7); 
   
   \node at (2,2){$\bullet$}; 
    \node at (5,2){$\bullet$};  \node at (8,2){$\bullet$}; 
 \node at (3,4){$\bullet$};  \node at (7,4){$\bullet$}; 
 
  \node at (1,7){$\bullet$};  \node at (4,7){$\bullet$}; 
   \node at (9,7){$\bullet$};  \node at (11,7){$\bullet$}; 
   
   \node at (6,5.5){$\cdots\cdots$}; 
   
   \node at (-1,2){$y_1$};   \node at (-1,4){$y_2$};   \node at (-1,7){$y_p$}; 
   
   \node at (-1, 5.8){$\vdots$};  

\node at (13,2){${\bf z} _1$}; \node at (13,4){${\bf z}_2$}; \node at (13,7){${\bf z}_p$}; 

 \end{tikzpicture}
 \caption{A point $\bf z$ of $X_\alpha^\Im$.}\label{fig:S-a-Im}
 \end{figure}
 
 We also have the similar real part map
 \[
 \Ren: \Sym^m\CC)\to\Sym^m(\RR), \quad m\geq 0. 
 \]
 Applying it  to  the divisor ${\bf z}_\nu$ corresponding to $\bf z$ as above gives a map 
  $\rho_{\alpha, \nu}: X_\alpha^\Ien \to\Sym^{\alpha_\nu}(\RR)$, $\nu=1,\cdots, p$. These maps combine into a projection
 \be
\rho_\alpha = (\rho_{\alpha,\nu})_{\nu=1}^p : X_\alpha^\Ien \lra \prod_{\nu=1}^p \Sym^{\alpha_\nu}(\RR),\quad \rho_{\alpha,\nu}({\bf z}) 
= \Ren({\bf z}_\nu),
\ee
and the following is clear.

\bigskip


\begin{prop}\label{prop:rho-fib}
 $\rho_\alpha$
 is a locally trivial fibration with fibers homeomorphic to $\RR^{l(\alpha)}$. 
 In particular, $X_\alpha^\Ien$ is a contractible CW-complex  of  pure real dimension   $n+l(\alpha)$. \qed
\end{prop}

\begin{exas}
(a)  If $\alpha = (n)$ is the partition with one part, then $X_{(n)}^\Im$ consists of $\{z_1, \cdots, z_n\}$ with
$\Im(z_1) = \cdots = \Im(z_n)$ and so is identified with $\RR\times \Sym^n(\RR)$. 

 \vskip .2cm
 
 (b) Note that $\Ien$ preserves the total sum of the points of a divisor so induces a map
\[
\Ien_0: \Sym^n_0(\CC) \lra\Sym^n_0(\RR). 
\]
Consider $n=2$. 
After undertifying $\Sym^2_0(\CC)$ with $\CC$ and $\Sym^2_0(\RR)$ with $\RR_{\geq 0}$ we get the map
 \[
\CC \lra \RR_{\geq 0}, \quad z\mapsto  (\Im \sqrt{z})^2. 
\]
The preimage of 0 under this map is the positive line $\RR_{\geq 0}$ and the other preimages
are parabolas positioned horizontally (so the positive line is the ``infinitely sharp parabola'').  
the induced imaginary strata on $\Sym^2_0(\CC)=\CC$ are the positive real line
  and its complement:
  \[
  S^\Im_{2}\cap\Sym^2_0(\CC) = \RR_{\geq 0}, \quad S_{1,1}^\Im\cap\Sym^2_0(\CC) = \CC - \RR_{\geq 0}.
  \]
\end{exas}

\paragraph{ The fine cells.}
We will further consider the preimages, under $\rho_\alpha$, of the  product strata of $\prod_{\nu=1}^p \Sym^{\alpha_\nu}(\RR)$
and call these preimages the {\em fine cells}. From this point of view,  fine cells are labelled by pairs $(\alpha, \Gamma)$, where
$\alpha = (\alpha_1, \cdots, \alpha_p)\in\2^{n-1}$ and $\Gamma = (\gamma^{(1)}, \cdots, \gamma^{(p)})$ 
is a sequence of ordered partitions, with $\gamma^{(\nu)}\in\2^{\alpha_\nu-1}$.   

\vskip .2cm

Note that for a given $\alpha$,  a datum of $\Gamma$  as above is equivalent to a datum of a single partition
$\gamma\in \2^{n-1}$ refining $\alpha$. This $\gamma$ is obtained simply by writing all the parts of all the $\gamma^{(i)}$
together in a single sequence. Thus, a fine cell is labelled by a pair $\alpha\leq \gamma$ of ordered partitions of $n$,\
and will be denoted $X_{[\alpha:\gamma]}$.  Thus
\[
X_\alpha^\Ien \,=\,\bigsqcup_{\gamma\geq \alpha} X_{[\alpha:\gamma]}. 
\]

\begin{prop}\label{prop:fine-cells}
(a) $X_{[\alpha:\gamma]}$ is a topological cell (space homeomorphic to an open ball) of dimension $l(\alpha) + l(\gamma)$. 

\vskip .2cm

(b) The collection of the $X_{[\alpha:\gamma]}$, $\alpha\leq\gamma$, 
 forms a cell decomposition (stratification into cells) of $\Sym^n(\CC)$
refining the complex stratification $S$. More precisely,  let  $\lambda\in\Pc(n)$ is an unordered partition of $n$. 
 The complex stratum $X_\lambda^\CC$ is the following union of fine cells:
 \[
 X_\lambda^\CC \,=\,\bigsqcup_{\alpha\leq \gamma\atop \ol\gamma=\lambda} X_{[\alpha:\gamma]}. 
 \]
\vskip .2cm

(c) We have $X_{[\alpha:\gamma]}\subset \ol{X_{[\alpha':\gamma']}}$
if and only if $\alpha\leq\alpha'$ and $\gamma\leq\gamma'$ in $\2^{n-1}$. 

\vskip .2cm

\vskip .2cm

(d) In addition, for any $\alpha\in\2^{n-1}$ and $\lambda\in\Pc(n)$ we have
\[
X_\alpha^\Ien \cap X_\lambda^\CC\, =\,\bigsqcup_{\ol\gamma=\lambda\atop \alpha \leq\gamma} X_{[\alpha:\gamma]}.
\]

\end{prop}

\noindent {\sl Proof:} Proposition \ref{prop:rho-fib} implies part (a). The other parts are straightforward. \qed

\vskip .2cm

We denote the stratification into fine cells by $S^{\fine}$. 

   \begin{exas}\label{ex:fuchs-1}
   
   (a) For $n=2$ the  stratification induced by $S^\fine$ on $\Sym^2_0(\CC)=\CC$ is the cell decomposition
   consisting of the point $0$, the open half line $\RR_{>0}$ and the complement $\CC\setminus \RR_{\geq 0}$.
   This example illustrates the fact that the cell decomposition  $S^{\fine}$ is not regular: the closure of one cell can approach a lower
   dimensional cell from more than one direction.
   
   \vskip .2cm

 (b)   The open complex stratum $\Sym^n_\neq(\CC) = X^\CC_{1^n}\subset\Sym^n(\CC)$ is the union of fine 
   cells $X_{[\alpha: 1^n]}$ for all $\alpha=(\alpha_1, \cdots, \alpha_p)\in\2^{n-1}$.
   This is the classical Fox-Neuwirth  cell decomposition of $\Sym^n_\neq(\CC)$, see \cite{fox-neuwirth} \cite{fuchs} used 
   for computing the cohomology of the braid
   group $\on{Br}_n = \pi_1(\Sym^n_\neq(\CC))$. 
   \end{exas}


\subsection{Perverse sheaves on $\Sym^n(\CC)$: generalities}

Let $\Vc$ be an abelian category. For each $n\geq 1$ we
consider the category $\Perv(\Sym^n(\CC), S,\Vc)$ of $\Vc$-valued perverse sheaves on $\Sym^n(\CC)$ smooth with respect to the diagklnal stratification.
An elementary description of this category
is known for $n\leq 3$. The case $n=1$ is trivial (the objects  are constant sheaves), and the cases $n=2,3$
are recalled in Ch. \ref{sec:examples}.  In this section we collect several general remarks.

\paragraph{Purity of the skeleton.}

 \begin{prop}\label{prop:K-pure}
 The skeleton $K=\Sym^n(\RR)$ is pure and self-dual for $(X=\Sym^n(\CC), S)$ with respect to perverse sheaves with values in
  any category $\Vc$. 
 \end{prop}
 
 To prove this,  
note that we have the diagram
  \be
    \xymatrix{
  \CC^n\ar[r]^{{\Ien}_{\CC^n}}
  \ar[d]_{p=p_\CC}& \RR^n\ar[d]^{p_\RR}
  \\
  \Sym^n(\CC) \ar[r]^{\Ien} & \Sym^n(\RR),
  }
  \ee
  where the vertical arrows are the  maps of taking the quotient by $\Sigma_n$.
  Inside $\CC^n$, we have the {\em root arrangement} (for the Lie algebra $\gen\len_n(\CC)$)
\[
\Hc_\CC \,= \bigcup_{0\leq i<j\leq n} H_{ij}^\CC, \quad H_{ij}^\CC \,=\,\bigl\{ (z_1, \cdots, z_n) \, \bigl| \, z_i=z_j \bigr\}. 
\]
 It is the complexification of the similar arrangement $\Hc$ in $\RR^n$.
 
  \begin{prop}\label{prop:p*F-perv}
 Suppose $\Fc\in\Perv(\Sym^n(\CC), S, \Vc)$. Then:
 \begin{itemize}
 \item[(a)] $p^*\Fc \in\Perv(\CC^n, \Hc, \Vc)$.
 
 \item[(b)] $\Fc= (p_* p^*\Fc)^{\Sigma_n}$. 
 \end{itemize}
 \end{prop}
 
 Proposition \ref{prop:p*F-perv} implies  that we can embed $\Perv(\Sym^n(\CC), S, \Vc)$ into the
 category of $\Sigma_n$-eqiivariant perverse sheaves on $(\CC^n, \Hc)$. Objects of the latter category can,
 in principle, be analyzed by the methods of \cite {KS-hyp-arr} \cite{weissman}.  In this paper we use
 a different approach, intrinsic to $\Sym^n(\CC)$. It has the advantage of involving linear algebra data
 of much smaller size and can be seen as  a generalization of the most classical $(\Phi,\Psi)$-description
 of $\Perv(\CC,0)$ (Example \ref{ex:phi-psi}). 
 
 \vskip .2cm

\noindent\underbar {\sl Proof of Proposition \ref {prop:p*F-perv}:} (a) Because the strata of $S$ are the images of the strata of $S_\Hc$, we have that $p*\Fc$ is an $S_\Hc$-constructible complex.
Let us verify the perversity condition $(P^+)$ for $p^*\Fc$. Since $p$ preserves the dimensions of the strata, the condition follows from the
same condition for $\Fc$. The condition $(P^-)$ for $p^*\Fc$ is equivalent to $(P^+)$ for the Verdier dual $\DD(p^*\Fc)$, which is identified with
$p^*(\DD\Fc)$. So $(P^-)$ for $p^*\Fc$ follows from $(P^+)$ for $\DD\Fc$, i.e., from $(P^-)$ for $\Fc$. Here $\DD$ is understood in the
sense of \eqref{eq:verdier-gener}. 

Part (b) is clear. \qed
 
 \vskip .2cm
 
 \noindent \underbar{\sl Proof of Proposition  \ref{prop:K-pure}.} Proposition \ref{prop:p*F-perv} implies that
\be\label{eq:RK-cover}
\ul{R\Gamma}_K \Fc\,=\, \bigl( p_* \ul{R\Gamma}_{\RR^n}(p^*\Fc)\bigr)^{\Sigma_n}, 
\ee
and so our statement follows from Proposition \ref{prop:takeuchi}(a). \qed

\paragraph{The total vanishing cycle functor.} For $\Fc\in\Perv(\Sym^n(\CC), S, \Vc)$ we denote
\be\label{eq:total-van}
\Phi_\tot(\Fc) \,=\, R_K(\Fc)_0\,\in \, \Ob(\Vc)
\ee
the stalk at $0\in K$ of the sheaf $R_K(\Fc)$ and call it the {\em total vanishing cycle space} of $\Fc$. 

Recall that the coefficients $a_1,\cdots, a_n$ of an indeterminate monic polynomial form a coordinate system on $\Sym^n(\CC)$,
and $\Sym^n_0(\CC)$ is defined by $a_1=0$. 

\begin{prop}
\[
\Phi_\tot(\Fc) \,\simeq \, \Phi_{a_2}\bigl(\Fc|_{\Sym^n_0(\CC)}\bigr)_0
\]
is identified with the stalk at $0$ of the perverse sheaf of vanishing cycles of the restriction $\Fc|_{\Sym^n_0(\CC)}$
with respect to the linear function $a_2$. 
\end{prop}

\noindent {\sl Proof:} 
Under the projection $p: \CC^n\to \Sym^n(\CC)$ the coordinates $a_m$ are the elementary symmetric functions of the standard
coordinates $x_1,\cdots, x_n$ on $\CC^n$. In particular, the quadratic form $Q(x) = x_1^2+\cdots + x_n^2$
descends to the function $a_1^2-2a_2$ on $\Sym^n(\CC)$. So by \eqref{eq:RK-cover} and   \eqref{eq:hyperb-van},
\[
\Phi_\tot(\Fc) \,\simeq \, \Phi_Q(p^*\Fc)^{\Sigma_n}_0 \,=\, \Phi_{a_1^2-2a_2}(\Fc)_0 \,=\,\Phi_{-2a_2}\bigl(\Fc|_{\Sym^n_0(\CC)}\bigr)_0,
\]
as claimed.  \qed

 \paragraph{ The Cousin resolution.}
 
 We now formulate a general scheme of analyzing objects of $\Perv(\Sym^n(\CC), S,\Vc)$
   which is similar in spirit
 to the one used for hyperplane arrangements in \cite{KS-hyp-arr} but is intrinsic to $\Sym^n(\CC)$.
   We fix a $\k$-linear abelian category $\Vc$. 
 Recall that $j_\alpha: X_\alpha^\Ien\to \Sym^n(\CC)$ is the embedding of the imaginary stratum corresponding
 to $\alpha\in\2^{n-1}$. 
 
\begin{prop}\label{prop:cousin}
Let $\Fc\in \Perv(\Sym^n(\CC), S, \Vc)$.

\vskip .2cm

(a) For any $\alpha\in\2^{n-1}$ the complex $j_\alpha^! \Fc$  on $X_\alpha^\Im$ is quasi-isomorphic to a single sheaf in degree 
equal to $-l(\alpha)$. 

\vskip .2cm

(b) Denoting the sheaf in (a) by $\wt\Ec_\alpha = \wt\Ec_\alpha(\Fc)$, the complex 
$ Rj_{\alpha*}\wt\Ec_\alpha$
reduces to a single sheaf $\Ec_\alpha=\Ec_\alpha(\Fc) = R^0j_{\alpha *} j_\alpha ^!\Fc$. This sheaf is constant on each fine cell. 

\vskip .2cm

(c)  $\Fc$ has an explicit representative (the {\em Cousin resoluton}) of the form
\[
\Ec^\bullet(\Fc) \,=\,\biggl\{ \Ec_{(1,\cdots, 1)} (\Fc) \buildrel\ul \delta\over\to \bigoplus_{l(\alpha) = n-1} \Ec_\alpha(\Fc) \buildrel
\ul \delta\over\to \cdots \buildrel \ul\delta\over\to \Ec_{(n)}(\Fc)\biggr\}.
\]
 The grading of this complex is such that $\Ec_n(\Fc)$ is situated in degree $0$. 
\end{prop}

\noindent{\sl Proof:} (a) Suppose first that $\alpha = (n)$, so $X_\alpha^\Ien$ consists of divisors whose
all components have the same imaginary part. In this case our statement is essentially equivalent to
 Proposition   \ref{prop:K-pure}: indeed, $K$ consists of divisors whose all components have the same
 imaginary part which is, moreover, equal to $0$.  So we establish the claim for $\alpha=(n)$. 
 Let now 
  $\alpha = (\alpha_1,\cdots, \alpha_p)$ be arbitrary  and suppose $\bz\in X_\alpha^\Ien$ is any point.
Let us write
 $\bz=\bz_1 +\cdots + \bz_p$ as in \eqref{eq:S-a-split}.  Then, near $\bz$, the pair $(X_\alpha^\Ien, \Sym^n(\CC))$
 is homeomorphic to the product of the pairs $(X_{(\alpha_i)}^\Ien,  \Sym^{\alpha_i}(\CC))$, and so
 our claim follows. 
 
 \vskip .2cm
 
 (b)
  Let $\bz\in \Sym^n(\CC)$. The stalk  at $\bz$ of $R^q j_* \wt\Ec^\alpha(\Fc)$ is $H^q(U\cap X_\alpha^\Ien, \wt\Ec^\alpha(\Fc))$, where 
   $U$ is a small ball around $\bz$ in $\Sym^n(\CC)$. 
   We note that  the imaginary strata $X_\alpha^\Ien$ are the images, under $p$, of the {\em tube cells}
  $\RR^n+i C$ where $C$ runs over the cells of the real arrangement $\Hc_\RR$, see \cite {KS-hyp-arr} , \S 3C. 
   Therefore 
    $U\cap X_\alpha^\Ien$ is the disjoint union of the images, under $p$, of the intersections of balls in $\CC^n$
  with the tube cells and $H^q$ of (the pullback of) $\wt\Ec^\alpha(\Fc)$ over such intersections vanish for $q>0$. 
  
  \vskip .2cm
  
  (c) This  is a general formal consequence  of (a) and (b) and  of  a spectral sequence  (Postnikov system) 
  corresponding to any complex of sheaves
on any filtered topological space, see, e.g.,   \cite{KS-hyp-arr} \S 1B.

\subsection{Expected answer:  bar-construction for bialgebras}\label{subsec:bar-bialg}

\paragraph{Representations and double representations of $\2^{n-1}$.}

We retain the notation of \S \ref{subsec:symprod-basics}.

\begin{Defi} Let $\Vc$ be an abelian category.

(a) By a {\em representation} of $\2^{n-1}$ in $\Vc$  we mean a covariant functor $\gamma:\2^{n-1}\to \Vc$, 
 i.e., a commutative cube of objects and morphisms
\[
\begin{gathered}
\bigl( E_\alpha, \alpha  \in\2^{n-1};  \,\,\, \gamma_{\alpha\alpha'}: E_\alpha \lra E_{\alpha'}, \,\,\, \alpha\leq\alpha'\bigr),\\
\gamma_{\alpha\alpha''}=\gamma_{\alpha'\alpha''}\circ\gamma_{\alpha\alpha'} \text{ for }
\alpha\leq\alpha'\leq\alpha''. 
\end{gathered}
\]
 We denote by $\Rep(\2^{n-1}, \Vc)$ the category of representations of $\2^{n-1}$ in $\Vc$
 
(b)  By an {\em anti-representation} of $\2^{n-1}$ in $\Vc$  we mean a contravariant functor $\delta: \2^{n-1}\to\Ac$,
 i.e., a commutative cube of objects and morphisns of the form
\[
  \bigl( E_\alpha, \alpha  \in\2^{n-1};  \,\,\, \delta_{\alpha'\alpha}: E_{\alpha'} \lra E_{\alpha}, \,\,\, \alpha\leq\alpha'\bigr).
  \]
  
(c)  By a {\em double representation} of  $\2^{n-1}$  in $\Vc$ we mean a diagram in $\Vc$ 
 \[
 \bigl( E_\alpha, \alpha  \in\2^{n-1};  \,\,\, \gamma_{\alpha\alpha'}: E_\alpha \lra E_{\alpha'},\,\,
  \delta_{\alpha'\alpha}: E_{\alpha'} \lra E_{\alpha}, 
  \,\,\, \alpha\leq\alpha'\bigr)
 \]
 such that $(E_\alpha, \gamma_{\alpha\alpha'})$ is a representation of $\2^{n-1}$, and
 $(E_\alpha, \delta_{\alpha'\alpha})$ is an anti-representation. We denote by $\Rep^{(2)}(\2^{n-1},\Vc)$
 the category of double representations of $\2^{n-1}$ in $\Vc$.
 \end{Defi}
 
  Proposition \ref{prop:K-strata} implies at once: 
 
\begin{cor}
The category $\Sh(K, S_K, \Vc)$  is identified with  $\Rep(\2^{n-1}, \Vc)$.   \qed
 \end{cor}
 
\begin{rems}
(a)  Note that a representation of $\2^{n-1}$ is determined by its {\em elementary arrows} 
 \[
 \gamma_{\alpha\alpha'}, \,\,\,\alpha \leq \alpha', \,\, l(\alpha')=l(\alpha) +1
 \]
 which correspond to the edges of the cube. The existence of compatible $\gamma_{\alpha\alpha'}$
 for all $\alpha\leq\alpha'$ means that the cube is commutative in the usual sense (each 2-dimensional
 face is a commutative square). Similarly for anti-representations. 
 \vskip .2cm
 
 (b) Any representation $E = (E_\alpha, \gamma_{\alpha\alpha'})$
  of $\2^{n-1}$ (a commutative cube) in $\Vc$ can be converted into a cochain
 complex in $\Vc$ by putting alternating $\pm$ signs on the elementary
arrows and summing over $\alpha$ with fixed $l(\alpha)$. We denote the
complex thus obtained by
 \[
 C^\bullet(E) \,\,=\,\,\biggl\{ E_{(n)} \buildrel\gamma\over\to \bigoplus_{\stackrel{\alpha_1+\alpha_2=n}{ a_i>0} } E_{(\alpha_1, \alpha_2)} \buildrel\gamma\over
\to  \bigoplus_{\stackrel{\alpha_1+\alpha_2+\alpha_3=n}{ a_i>0} } E_{(\alpha_1, \alpha_2, \alpha_3)}\buildrel\gamma\over \to\cdots \buildrel\gamma\over\to  
E_{(1,1,\cdots, 1)} \biggr\}. 
 \]
 We fix its grading so that $E_{(n)}$ has degree $1$. 
 
 Similarly, any anti-representation $E= (E_\alpha, \delta_{\alpha'\alpha})$ gives rise to a chain complex
 $C_\bullet(E)$. We can view  $C_\bullet(E)$  as a {\em homological complex} (differental $\delta$ lowering the degree by $1$)
 with the same grading as $C^\bullet(E)$. Alternatively, we can convert it into a cohomological complex by reversing the
 degree, that is, by putting
  \[
 C_\bullet(E) \,\,=\,\,\biggl\{ E_{(1,1,\cdots, 1)} \buildrel\delta\over
\to \cdots \buildrel\delta\over
\to 
 \bigoplus_{\stackrel{\alpha_1+\alpha_2+\alpha_3=n}{ a_i>0} } E_{(\alpha_1, \alpha_2, \alpha_3)}
 \buildrel\delta\over
\to 
\bigoplus_{\stackrel{\alpha_1+\alpha_2=n}{ a_i>0} } E_{(\alpha_1, \alpha_2)}
\buildrel\delta\over
\to E_{(n)}  \biggr\}. 
 \] 
 with  $E_{(n)}$ 
 in degree $(-1)$.  
 
 A double representation $E=(E_\alpha, \gamma_{\alpha\alpha'}, \delta_{\alpha'\alpha})$ gives therefore a complex with
 two differentials $\gamma$ and $\delta$, one raising, the other lowering the degrees.

  \end{rems}

\paragraph{ (Anti-)representations from (co)algebras.} \label{par:birep-bial}
Let $(\Vc, \otimes,\1)$ be a $\k$-linear monoidal abelian category. We assume that $\otimes$ is exact in each variable.

\vskip .2cm

By a {\em connected graded algebra} in $\Vc$ we mean a $\ZZ_{\geq 0}$-graded
associative algebra $A=\bigoplus_{n\geq 0} A_n$ in $\Vc$
with $A_0=\1$ being the unit. We denote by $\mu_{p,q}: A_p\otimes A_q\to A_{p+q}$ the components of the
multiplication $\mu$ in $A$.

\vskip .2cm

Given such an $A$, for each $n\geq 0$ we have the anti-representation $\BB_n(A)$ of $\2^{n-1}$ called the
$n$th {\em bar-cube} of $A$. Explicitly, 

\begin{enumerate}
 \item [(1)] 
$
\BB_n(A)_{(\alpha_1, \cdots, \alpha_p)}= A_{\alpha_1} \otimes \cdots \otimes A_{\alpha_p}
$
 \item[(2)]  The  elementary maps $\delta_{\alpha'\alpha'}$   are  given by the 
appropriate components  $\mu_{pq}$ tensored by the identities.
\end{enumerate}

The commutativity of the cube follows from the associativity of $A$. For example, for $n=3$  the diagram $\BB_3(A)$ has the form
\[
 \xymatrix{
&A_2\otimes A_1
 \ar@<.5ex>[dl]_{\mu_{2,1}}
&
\\
A_3
 && A_1\otimes A_1\otimes A_1. 
\ar@<-.5ex>[dl]^{\Id\otimes\mu_{1,1}}
\ar@<.5ex>[ul]_{\mu_{1,1}\otimes\Id}
 \\
&A_1\otimes A_2. 
\ar@<-.5ex>[ul]^{\mu_{1,2}}
&
}
\]
The corresponding chain complex, written as a cohomological one, will be called the $n$th {\em bar-complex} of $A$ and
denoted
\be
B_n(A)\,=\, 
\biggl\{ A_1^{\otimes n} \buildrel\delta\over
\to \cdots \buildrel\delta\over
\to 
 \bigoplus_{\stackrel{\alpha_1+\alpha_2+\alpha_3=n}{ a_i>0} } A_{\alpha_1}\otimes A_{\alpha_2}\otimes A_{\alpha_3} 
 \buildrel\delta\over
\to 
\bigoplus_{\stackrel{\alpha_1+\alpha_2=n}{ a_i>0} }  A_{\alpha_1}\otimes A_{\alpha_2} 
\buildrel\delta\over
\to A_n \biggr\},
\ee
its grading normalized so that $A_n$ is in degree $(-1)$. We note that
\be\label{eq:bar-tor}
H^j B_n(A) \,=\,\Tor_{-j,n}^A(\1, \1)
\ee
is the $n$th graded component of $\Tor_{-j}^A(\1,\1)$ with respect to the additional grading coming from the grading of $A$. 

\vskip .2cm

Dually, by a {\em graded coconnected coalgebra}  in $\Vc$ we mean a  $\ZZ_{\geq 0}$-graded
coassociative coalgebra $A=\bigoplus_{n\geq 0} A_n$ in $\Vc$
with $A_0=\1$ and the counit being the projection to $A_0=\1$. We denote by 
$\Delta_{p,q}: A_{p+q}\to A_p\otimes A_q$ the components of the comultiplication $\Delta$  in $A$. 
We note that 
\[
\Delta_{0,n}: A_n\lra A_0\otimes A_n= \1\otimes A_n=A_n
\]
 is the identity, and similarly for $\Delta_{n,0}$. In particular, $A_1$ ``consists of primitive elements", i.e.,
 $\Delta|_{A_1} = \1\otimes\Id + \Id\otimes \1$. 
 
 \vskip .2cm
 
 A graded coconnected coalgebra $A$ gives a representation $\BB_n^\dagger(A)$ of $\2^{n-1}$ in $\Vc$
 called the $n$th {\em cobar-cube}. Explicitly, 
 \begin{enumerate}
 \item [(1)] 
$
\BB^\dagger_n(A)_{(\alpha_1, \cdots, \alpha_p)}= A_{\alpha_1} \otimes \cdots \otimes A_{\alpha_p}
$
 \item[(2)]  The  elementary maps $\gamma_{\alpha'\alpha}$   are  given by the 
appropriate components  $\Delta_{pq}$ tensored by the identities.
\end{enumerate}
The cochain complex of $\BB^\dagger_n(A)$ will be called the $n$th {\em cobar-complex} of $A$ and denoted
\be
B^\dagger_n(A)\,=\, 
\biggl\{ 
 A_n
 \buildrel\gamma\over
\to 
\bigoplus_{\stackrel{\alpha_1+\alpha_2=n}{ a_i>0} }  A_{\alpha_1}\otimes A_{\alpha_2} 
 \buildrel\gamma\over
\to 
 \bigoplus_{\stackrel{\alpha_1+\alpha_2+\alpha_3=n}{ a_i>0} } A_{\alpha_1}\otimes A_{\alpha_2}\otimes A_{\alpha_3} 
\buildrel\gamma\over
\to 
\cdots \buildrel\gamma\over
\to 
A_1^{\otimes n}
  \biggr\},
\ee
its grading normalized so that $A_n$ is in degree $1$. As in the algebra case,
\[
H^j B_n^\dagger(A)  \,=\,\on{Cotor}^{j,n}_A(\1, \1).
\]

\paragraph{Bialgebras and Hopf algebras in a braided monoidal category.}\label{par:braid-bial}
Let $(\Vc, \otimes, \1, R)$ be an   $\k$-linear braided monoidal category. 
So for any $V_1, V_2\in \Vc$ one has the $R$-matrix
$$
R_{V_1,V_2}:\ V_1\otimes V_2 \iso V_2\otimes V_1
$$
satisfying the braiding axioms  \cite{joyal-street}. These axioms imply, in particular, the Yang-Baxter equation for any
$R_{V,V}: V\otimes V\to V\otimes V$. (Some physicists would call a braiding structure 
{\em a statistics.}). 

\vskip .2cm

  Recall \cite{majid}  \cite {takeuchi-hopf} that  we can speak about bialgebras in $\Vc$. More precisely, if $A$ is an associative algebra in 
  $\Vc$ with
 multiplication $\mu: A\otimes A\to A$, then $A\otimes A$ becomes an associative algebra with multiplication being the
 composition
 \be\label{eq:AA-alg}
 \begin{gathered}
 (A\otimes A)\otimes (A\otimes A) \, \simeq \,  A\otimes (A\otimes A)\otimes A
 \buildrel \Id_A\otimes R_{A,A}\otimes \Id_A\over\lra A\otimes (A\otimes A)\otimes A \,  \simeq
 \\
 \simeq \,  (A\otimes A)\otimes (A\otimes A) \buildrel \mu\otimes\mu  \over\lra A\otimes A. 
 \end{gathered}
 \ee
 
 \begin{Defi}\label{def:br-bialg}
 (a) A {\em bialgebra} in $\Vc$ is an object $A$ made into an associative algebra with  unit $u: \1\to A$  via $\mu: A\otimes A\to A$ 
 and into a coassociative coalgebra  with counit $\eps: A\to\1$
 via $\Delta: A\to A\otimes A$
  so that $\Delta$ is a morphism of algebras where $A\otimes A$ has the algebra structure 
 \eqref{eq:AA-alg}. 
 
 \vskip .2cm
 
 (b) An {\em antipode} in a bialgebra $A$ in $\Vc$ is a morphism $\sigma: A\to A$ such that the two 
 compositions
 \[
 \mu\circ (\sigma \otimes\Id)\circ\Delta, \,\,\, \mu\circ(\Id\otimes\sigma)\circ\Delta: A \lra A
 \]
 are equal to $\eta\circ\eps$. A bialgebra equipped with antipode is called a {\em Hopf algebra}. 
 \end{Defi} 
 
 We further recall the concept of the {\em duals}, see, e.g., \cite{bakalov} \S 2.1 or \cite{takeuchi-hopf} \S 2. 
  Let $\Vc$ be any monoidal category and $V,W$ be two objects.
 A {\em duality datum} between $V$ and $W$ is a pair of morphisms
 \[
 u: \1\lra V\otimes W, \quad c: V\otimes W\lra \1
 \]
 such that the following diagrams commute:
 \[
 \xymatrix{
 V\ar[r]^{\hskip -0.7cm V\otimes u} \ar[dr]_{\Id}
 & V\otimes W\otimes V\ar[d]^{c\otimes V}
 \\
  & V
 }
 \quad\quad\quad
 \xymatrix{
 W \ar[r]^{\hskip -0.7cm  u\otimes W} \ar[dr]_{\Id}&V\otimes W\otimes V
 \ar[d]^{W\otimes c}
 \\
 & V
 }
 \]
 In this case $W$ is defined by $V$ uniquely up to a unique isomorphism, is called the {\em right dual} of $V$  and denoted $W=V^*$.
 Similarly $V$ is defined by $W$ uniquely up to a unique isomorphism, is called the {\em left dual} of $W$  and denoted $V={} ^*W$.
 An object having a left (resp. right) dual is called {\em left} (resp, {\em right}) {\em dualizable}. 
 
 \vskip .2cm

  Any monoidal category $\Vc$ can be considered as a 2-category
 $\Omega^{-1}\Vc$ with one object, so  objects $V\in \Vc$ become 1-morphisms of $\Omega^{-1}\Vc$.  
 The concept of  left and right duals of objects in $\Vc$ becomes a particular case of the concept of the left and right adjoint of
 a 1-morphism of a 2-category (e.g., of the left and right adjoint of a functor between usual categories).  The duality datum
 corresponds to the unit and counit of an adjunction. 
 
 \vskip .2cm
 
 Let  us now return to the case when $\Vc$ is a braided monoidal category. In this case a left dual is also a right dual and vice versa,
 so we use the notation $V^*$ and 
 speak simply about {\em dualizable objects}. A morphism $f: V\to W$ between dualizable objects gives the
 adjoint morphism $f^t: W^*\to V^*$. 
 If every object of $\Vc$ is dualizable, then $\Vc$ is called {\em rigid}. In this case we can think of $V\mapsto V^*$ as
 a perfect duality on $\Vc$. 
 
 \begin{prop}\label{prop:dual-antipode}
 Let $\Vc$ be a  braided monoidal $\k$-linear abelian category. 
 Let $A$ be a bialgebra  $\Vc$ which is dualizable as an object of $\Vc$. Then:
 \vskip .2cm
 
 (a)  $A^*$ is
 also a bialgebra with $\mu_{A^*}=\Delta_A^t$ and $\Delta_{A^*} = \mu_A^t$. Further, if $A$ is a Hopf algebra,
 then $A^*$ is a Hopf algebra with antipode $\sigma_{A^*} = \sigma_A^t$. 
 
 \vskip .2cm
 
 (b) If $A$  is a Hopf algebra, then the antipode in $A$ is invertible. 
 \end{prop}
 
 \noindent{\sl Proof:} (a)  is Theorem 2.16  and (b) is a particular case of Theorem 4.1 of \cite{takeuchi-hopf}.\qed

 \paragraph{ Primitive bialgebras. Relation to $\Perv(\Sym^n(\CC), S)$.} Let $\Vc$ be a  braided monoidal $\k$-linear abelian category
 as above.

 \begin{Defi}\label{def:prim-gen-bialg}
 A {\em primitive  bialgebra} in $\Vc$ is
 a graded bialgebra $A = \bigoplus_{n\geq 0} A_n$  which is connected and coconnected. 
  \end{Defi}
 
 For any partition $(\alpha_1,\cdots, \alpha_k)\in\2^{n-1}$ of $n$, $\alpha_j>0$ 
  and 
   a primitive bialgebra $A$,   we denote
 \[
 \Delta_{\alpha_1,\cdots, \alpha_k}: A_n\lra A_{\alpha_1}\otimes\cdots\otimes A_{\alpha_k}, \quad 
 \mu_{\alpha_1,\cdots, \alpha_k}: A_{\alpha_1}\otimes\cdots\otimes A_{\alpha_k}
 \lra A_n
 \]
 the corresponding components of the iterated (co)multiplication. For the trivial partition $(n)$  of $n$ we put
 $\Delta_{(n)} = \mu_{(n)} = \Id_{A_n}$. 
 
 \begin{prop}\label{prop:antipode}
  Let $A$ be a primitive bialgebra. For $n\geq 0$ define the morphism $\sigma_n: A_n\to A_n$  by
  \[
  \sigma_n \,\, =\,\,   \sum_{(\alpha_1, \cdots, \alpha_k)\in \2^{n-1}} (-1)^{k-1}  \mu_{\alpha_1,\cdots, \alpha_k}\circ 
   \Delta_{\alpha_1,\cdots, \alpha_k}. 
  \]
 Then:
 
 \vskip .2cm
 
 (a)  $\sigma = \bigoplus_n \sigma_n: A\to A $ is an antipode for $A$. 
 
 \vskip .2cm
 
(b) If each $A_n$ is dualzable, then $\sigma$ is invertible. 
 \end{prop}
 
 \noindent {\sl Proof:}  Part  (a) is verified by a direct check,
 cf. \cite{loday} \S 1.2 for a similar formula in the non-braided case.  
  Part (b) follows from Proposition \ref{prop:dual-antipode}, if we consider $A$ itself as a dualizable object in the
 braided category formed by graded objects of $\Vc$. \qed
 
 \vskip .2cm
 
 Given a primitive bialgebra $A$ and $n\geq 0$, the co-representation $\BB_n(A)$ and the representation $\BB_n^\dagger(A)$
 of $\2^{n-1}$ combine into a double representation which we denote $\BBB_n(A)$ and call the
 $n$th {\em bibar-cube} of $A$. For example, for $n=3$  the diagram $\BBB_3(A)$ has the form
\be\label{eq:bibar-3}
 \xymatrix{
&A_2\otimes A_1\ar@<.5ex>[dr]^{\Delta_{1,1}\otimes\Id}
\ar@<.5ex>[dl]^{\mu_{2,1}}
&
\\
A_3\ar@<.5ex>[ur]^{\Delta_{2,1}} \ar@<-.5ex>[dr]_{\Delta_{1,2}} && A_1\otimes A_1\otimes A_1
\ar@<-.5ex>[dl]_{\Id\otimes\mu_{1,1}}
\ar@<.5ex>[ul]^{\mu_{1,1}\otimes\Id}
 \\
&A_1\otimes A_2. \ar@<-.5ex>[ur]_{\Id\otimes\Delta_{1,1}}
\ar@<-.5ex>[ul]_{\mu_{1,2}}
&
}
\ee

It is natural to expect that the double representation $\BBB_n(A)$ actually comes from a perverse sheaf on $\Sym^n(\CC)$
smooth with respect to the diagonal factorization $S$. In the next chapter we give a precise result to this effect.

\vfill\eject

\section{Factorizable perverse sheaves and bialgebras}\label{sec:fact-perp-bi}

\subsection{Shuffle algebras in braided categories}\label{subsec:shuf-braid}

 \paragraph{The tensor algebra.}
 We return to the situation of \S \ref{subsec:bar-bialg}\ref{par:braid-bial}, that is, of a braided $\k$-linear
monoidal category $(\Vc, \otimes,  \1, R)$. 
Let $V\in \Vc$; consider the {\em tensor algebra} 
\[
T_!(V) = \bigoplus_{n=0}^\infty\ V^{\otimes n}
\]
 with  the multiplication 
\[
V^{\otimes m}\otimes V^{\otimes n} \lra V^{\otimes m + n}
\]
begin the identity map. So as an algebra, $T_!(V)$ is free. 
It has a unique structure of a bialgebra in $\Vc$ for which 
 the comultiplication
\[
\Delta^! = \Delta^{!,R}:\ T_!(V) \lra T_!(V)\otimes T_!(V)
\]
is defined as the unique algebra map for which  $T_!^1(V) = V$ is primitive, that is, 
\[
\Delta^!|_{T_!^1(V)} = \1\otimes\Id + \Id\otimes \1: V\lra V\otimes V. 
\]
If we think of $V$ as consisting of ``vectors'' $x\in V$, then
 \[
\Delta^!(x) = x\otimes 1 + 1 \otimes x, \quad x\in V. 
\]
Since the algebra structure on  
  $T_!(V)\otimes T_!(V)$ is defined using the $R$ matrix, see \eqref{eq:AA-alg}, the comultiplication $\Delta^!$ depends on $R$. 
  We have put an index $R$ at $\Delta^!$ to stress its  
dependence on $R$.

\paragraph{The cotensor (big shuffle) algebra.} 
Dually, we have the {\em cotensor algebra} $T_*(V)$ which is a bialgebra in $\Vc$ 
with the same underlying object:
\[
T_*(V) = \bigoplus_{n=0}^\infty\ V^{\otimes n}. 
\]
Its comultiplications 
\[
\Delta_{m,n} = \Delta^{(*)}_{m,n}:\ V^{\otimes m + n} \lra V^{\otimes m}\otimes V^{\otimes n}
\]
are the identity maps, whereas 
the multiplication $\bigstar_R$ on $T_*(V)$  is the  so-called {\em shuffle product},  depending on $R$. That is, 
\[
\bigstar_R^{m,n}:\ V^{\otimes m}\otimes V^{\otimes n} \lra V^{\otimes m + n}
\]
acts as 
\be\label{eq:shuffle-prod}
\bigstar_R^{m,n} = \sum_{w\in \Sha_{m,n}} R_w.
\ee
Here 
\[
\Sha_{m,n} = \{w\in \Sigma_{m+n}|\ \text{if\ }i\leq m, j \geq m + 1\ \text{then\ }
w(i) < w(j)\}\subset \Sigma_{m+n}
\]
is the  subset of $(m,n)$-shuffles.  The braiding axioms for $R$
 imply the associativity of $\bigstar_R$. Sometimes $T_*(V)$ is referred to as the {\em big shuffle algebra}. 
 As a coalgebra, $T_*(V)$ is cofree. 
 The following is straightforward. 
 
 \begin{prop}
 (a) $T_!(V)$ and $T_*(V)$ are primitive (Definition \ref{def:prim-gen-bialg})  bialgebras in $\Vc$. In particular, they are Hopf algebras. 
 
 \vskip .2cm
 (b) Let  $V$ be a dualizable object of $\Vc$. Then $T_!(V)$ and $T_*(V)$ are dualizable as objects of the category of graded objects in $\Vc$
 and $T_*(V)$ is identified with the Hopf algebra dual to  $T_!(V)$.  \qed
 \end{prop} 
 
 \paragraph{The small shuffle (Nichols) algebra.}\label{par:nichols}
 Let us now assume, in addition,  that $\Vc$ is abelian and $\otimes$
is exact  in each variable. 

\vskip .2cm

 Note that we have a canonical morphism of algebras in $\Vc$
 \[
 \alpha: T_!(V) \lra T_*(V). 
 \]
 It sends the component $T_!^1(V)=V$ by the identity into the component $T_*^1(V)=V$ and is uniquely defined by this
 property. Indeed, $T_!(V)$ is freely generated by $V$, so to define its algebra morphism into anything, it is enough the
 specify the values on $V$ in an arbitrary way. 
 
 \begin{prop}
  $\alpha$ is a morphism of bialgebras.
 \end{prop}
 
 \noindent{\sl Proof:} 
 One has to check that $\alpha$ preserves comultiplication, i.e., that
\[
\alpha\otimes \alpha\circ \Delta_! = \Delta_*\circ \alpha:\ 
T_!(V) \lra T_*(V)  \otimes T_*(V).
\]
Both maps are morphisms of algebras and $T_!(V)$ is generated as an algebra by 
$V = T_!^1(V)$. But on this component both maps agree. In the element notation,  they both send any 
$x$ to $x\otimes 1 + 1\otimes x$. \qed

\begin{Defi}
The  {\em Nichols } (or {\em  small shuffle}) algebra  of an object $V\in \Vc$   is the bialgebra 
\[
T(V)_{!*} := \alpha(T_!(V)). 
\]
So it is a subalgebra of $T_*(V)$ and a quotient of $T_!(V)$. 
\end{Defi}

\begin{ex}[(abelian monodromy)]\label{ex:abel-mon}
Let $\Lambda \simeq \ZZ^r$ be a free abelian group of finite rank,  and 
$
c:\ \Lambda\otimes_\ZZ \Lambda \lra \k^*
$
be a  group homomorphism. 
If we pick  a $\ZZ$-basis $e_1, \ldots, e_r\in \Lambda$,  then $c$ will be determined by 
an  $r\times r$ matrix 
$
c_{ij} = c(e_i\otimes e_j) \in \k^*
$
Consider the category $\Vect^\Lambda$ of $\Lambda$-graded vector spaces whose 
objects we write as 
\[
V = \oplus_{\lambda\in \Lambda}\ V_\lambda. 
\]
We can introduce on it a braiding given by the $R$-matrices 
\[
R(x\otimes y) = c(\lambda,\mu) y \otimes x
\]
for $x\in V_\lambda, y\in V_\mu$  (``an abelian monodromy'', as a physicist would say).   

\vskip .2cm

In particular, let $A=\|a_{ij}\|$ be a Cartan matrix and put $c_{ij}=q^{a_{ij}}$. Let $V = \k^{r} = \Lambda\otimes_\ZZ\k$ with its standard 
$\Lambda$-grading. Let $\gen$ be the semi-simple Lie algebra with Cartan matrix $A$. Then:

\begin{itemize}

 \item
If $q$ is not a root of $1$, then
$
T_{!*}(V) \,=\, U_q^+(\gen)
$
is the positive part of  quantum universal enveloping algebra of $\gen$, see \cite{feigin-odesskii}  \cite{rosso}.

\item If $q$ is a root of $1$, then it was proved in \cite{rosso} that $T_{!*}(V) = {\bf u}_q^+(\gen)$ is the
positive part of Lusztig's {\em small} (or {\em restricted}, see \cite{lusztig:Rt1})   quantum enveloping algebra associated to $\gen$ and $q$.  In particular,
$T_{!*}(V)$ is finite-dimensional. 

\end{itemize}
  This can be seen as an interpretation of the Lusztig's construction of    $U_q^+(\gen)$ 
  (resp. ${\bf u}_q^+(\gen)$) as the  quotient
of a free Lie algebra by the kernel of the natural bilinear form $\langle -,-\rangle$, see  \cite{lusztig}, Ch. 1, \cite{rosso-bilin}.  
 In other words,  taking the image of $\alpha$ in forming the Nichols algebra has the effect of
  quotienting  by the kernel of $\langle -,-\rangle$,
which is the universal way to encode the quantum Serre relations in $U_q^+(\gen)$. 
\end{ex}

 
 \subsection {Factorizable perverse sheaves}\label{subsec:fact-perv}
 
 \paragraph{The $E_2$-operad and braided categories.}\label{par:E2-br}
 An interpretation going back to Deligne, views a braided monoidal structure in a category $\Vc$ as a datum, for each $n\geq 0$,  of
 $n$-fold tensor product functors in $\Vc$ labelled by $n$-tuples of distinct points in $\CC$. It is convenient to reformulate this
 approach as follows,  replacing points by disks.  
 
 \vskip .2cm
 
 Let $D=\{|z|\leq 1\}$ be the standard unit disk in $\CC$.
  We denote by  $E_2$  the operad of little $2$-disks in $\RR^2=\CC$, see \cite{may} \cite{CG}.
   Thus, for each $n\geq 0$, the topological space
  $E_2(n)$ consists of $n$-tuples of disjoint round open disks $(U_1, \cdots, U_n)$ in $D$. Equivalently, we can see $E_2(n)$
  as consisting of embeddings
  \[
  \phi = (\phi_1,\cdots, \phi_n): \, \bigsqcup_{i=1}^n D \, \lra D
  \]
such that each $\phi_i: D\to D$ is a composition of a dilation and a translation. This interpretation makes manifest the
{\em operadic compositions}
\[
E_2(k) \times E_2(n_1) \times \cdots \times E_2(n_k) \lra E_2(n_1+\cdots + n_k). 
\]
Thus $E_2$ is an operad in the category $\Top$ of topological spaces. Passing to fundamental groupoids, we get an
operad $\Pi_1(E_2)$ in the category of groupoids. We note that each $E_2(n)$ is a $K(\pi, 1)$-space, so
$\Pi_1(E_2(n))$ contains all the information about its homotopy type. 

\begin{prop}\label{prop:braid-e2}
Let $\Vc$ be a category. Giving a braided monoidal structure in $\Vc$ is equivalent to making $\Vc$ into a $\Pi_1(E_2)$-algebra
in the category of categories, that is, defining:
\begin{itemize}
\item [(1)] For any $n\geq 0$ and any disjoint open round disks $U_1,\cdots, U_n\subset D$, an $n$-variable functor
\[
\bigotimes^{(n)}\nolimits _{(U_1, \cdots, U_n)}: \,\,  \Vc^n \lra \Vc, \quad (V_1, \cdots, V_n) \mapsto
 \bigotimes\nolimits_{(U_1, \cdots, U_n)} V_i, 
\]
\item[(2)] For any path $\gamma:[0,1]\to E_2(n)$,  joining $(U_1, \cdots U_n)$ and $(U'_1, \cdots U'_n)$, a natural isomorphism of functors
\[
R_\gamma: \bigotimes^{(n)}\nolimits _{(U_1, \cdots, U_n)}\lra \bigotimes^{(n)}\nolimits _{(V_1, \cdots, V_n)}
\]
depending only on the homotopy class of $\gamma$. 
\item[(3)] These data are required to be compatible with the operadic compositions. \qed
\end{itemize}
\end{prop}

For any open $U\subset\CC$ let $U^n_\neq\subset U^n$ be the set of $(z_1,\cdots, z_n)$ where   $z_i\in U$ 
are such that $z_i \neq z_j$ for $i \neq j$. Given $n$ objects $V_1,\cdots, V_n\in\Vc$ and  $(z_1,\cdots, z_n)\in D^n_\neq$, we denote
\be\label{eq:braided-tensor}
\bigotimes\nolimits_{(z_1, \cdots, z_n)} V_i  \,=\,  \bigotimes\nolimits_{(U_1, \cdots, U_n)} V_i,
\ee
where the $U_i$ are sufficiently small open disks centered around $z_i$. This defines a local system on $D^n_\neq$
with values in $\Vc$. Since the embedding $D^n\neq \CC^n_\neq$ is a homotopy equivalence, this local system
extends canonically to a local system in $\CC^n_\neq$, and so we can and
 will use the notation $\bigotimes\nolimits_{(z_1, \cdots, z_n)} V_i $
for any $(z_1,\cdots, z_n)\in \CC^n_\neq$. 

\vskip .2cm

We will refer to the above point of view on braided categories (many tensor products operations, labelled by $(U_1,\cdots, U_n)$
or $(z_1,\cdots, z_n)$) as the {\em operadic point of view}, as opposed to the {\em classical point of view}
(one operation $\otimes$ together with braiding isomorphisms). 

\paragraph{Braidings on derived categories and perverse sheaves.}
Let $\Vc$ be a braided abelian $\k$-linear category, with $\otimes$ exact in each argument. We extend the braiding to
$C^b(\Vc)$, the category of bounded cochain complexes over $\Vc$ in a standard way. In  other words,  
  the operation $\otimes$
on $C(\Vc)$ is defined  to be   the usual
graded tensor product of comlexes (with differential defined by the Leibniz rule). The braiding on $C^b(\Vc)$ is
modified by the Koszul sign rule. That is, for objects $V,W\in\Vc$ and integers $m,n\in\ZZ$, the tensor product
and braiding of the shifted objects are defined by
\be\label{eq:koszul}
(V[m])\otimes (W[n])= (V\otimes W)[m+n], \quad R_{V[m], W[n]} = (-1)^{mn} R_{V,W}[m+n]. 
\ee
Because of bi-exactness of $\otimes$, the braided structure descends to $D^b(\Vc)$, the 
bounded derived category of $\Vc$. 

\vskip .2cm

For any stratified complex manifold $(X,S)$ this induces a braided monoidal structure on the derived category
$D^b(X,S,\Vc)$ of $S$-constructible complexes on $X$ with values in $\Vc$. The abelian subcategory of perverse sheaves
is typically not closed under $\otimes$. However, the {\em exterior} tensor product of perverse sheaves is perverse.

\vskip .2cm

That is, if $(X,S)$ and $(Y,T)$ are two stratified complex manifolds, we have the functor
\[
\boxtimes: \Perv(X,S,\Vc) \times \Perv(Y,T,\Vc) \lra\Perv(X\times Y, S\times T, \Vc).
\]
 Denoting by $\pi: X\times Y\to Y\times X$ the permutation, we have the external braiding isomorphism
 \[
 R_{\Fc, \Gc} : \Fc\boxtimes\Gc \lra \pi^*(\Gc\boxtimes\Fc). 
 \]
 Note that a ``locally constant perverse sheaf''  on $X$ is in fact a local system placed in degree $(-\dim X)$. 
 Therefore $R_{\Fc,\Gc}$ involves the Koszul sign rule even for ``locally constant'' $\Fc$ and $\Gc$. 
 
 \vskip .2cm
 
 Further, we extend the operadic point of view on the braiding to the external tensor products of perverse sheaves.
 That is, let $(X_i, S_i)$, $i=1,\cdots, m$, be several stratified complex manifolds and $\Fc_i\in\Perv(X_i, S_i, \Vc)$.
 Then for any disjoint open disks $U_1,\cdots, U_m$ in the unit disk $D$ we have the object 
 \[
 \bigboxtimes\nolimits_{(U_1,\cdots, U_m)} \Fc_i \,\, \in\,\, \Perv\left(\prod X_i, \prod S_i, \Vc\right). 
 \]

  
  \paragraph{ Factorizable perverse sheaves with values in a braided abelian category.} Let now $\Vc$ be a braided monoidal
  abelian category, as before.
  
      \vskip .2cm
    
    For any open set $U\subset \CC$ we denote $\Sym(U) = \bigsqcup_{n\geq 0} \Sym^n(U)$. 
   We denote by $S_{n,U}$ the diagonal stratification of $\Sym^n(U)$ and by $S_U$ the resulting stratification of $\Sym(U)$.
 In the case $U=\CC$ we abbreviate $S_{n,U}$ to $S_n$ and $S_U$ to $S$. 
 Denote by $\Perv(\Sym(U), S_U;\Vc )$ the product of the categories $\Perv(\Sym^n(U), S_{n,U};\Vc)$. 
 Thus an object of  $\Perv(\Sym(U), S_U;\Vc)$
 is a collection $\Fc=(\Fc_n)_{n\geq 0}$ of perverse sheaves $\Fc\in\Perv(\Sym^n(U), S_{n,U};\Vc)$. 
 
 \vskip .2cm
 
 Let $(U_1, \cdots U_m)$ be a point of $E_2(m)$, i.e., a tuple of disjoint round open disks in $D$. 
 Consider the Cartesian product $\prod_{i=1}^m \Sym(U_i)$ with its stratification  $\Pi S_{U_i}$
 formed by the products of the strata of the $S_{U_i}$. 
  Since $U_1, \cdots, U_m$ are
    disjoint,  we have a canonical identification (``addition map'')
    \be\label{eq:add-map}
  a:   \prod_{i=1}^m  \Sym(U_i) \buildrel\sim\over \lra  \Sym\biggl( \bigcup_{i=1}^m U_i\biggr). 
    \ee

 \begin{Defi}\label{def:FPS}
 A {\em factorizable perverse sheaf} on $\Sym(\CC)$ with values in $\Vc$ is an object $\Fc= (\Fc_n)$ of $\Perv(\Sym(\CC), S;\Vc)$ together with
 the  data of
 
 \begin{itemize}
 \item[(*)] For each $m\geq 0$ and each $(U_1, \cdots U_m)\in E_2(m)$, an isomorphism
 \[
\mu_{U_1, \cdots, U_m}:  \bigboxtimes \nolimits_{(U_1, \cdots, U_m)} \Fc|_{\Sym(U_i)} \lra    a^*\left(  \Fc|_{\Sym(\bigcup U_i)}\right),
 \]
these isomorphisms compatible with the operadic compositions.
 \end{itemize} 
 \end{Defi}
 
 \begin{exas}
 The condition ``compatible with operadic compositions'' includes, in particular, the following:
 
 \vskip .2cm
 
 (a) Associativity for binary tensor products. That is, suppose we have  three disjoint round disks $U_1, U_2, U_3\subset D$. 
 Then the two ways of decomposing the map  
 \[
 \mu_{U_1, U_2, U_3}: \boxtimes_{(U_1,U_2, U_3)} \Fc|_{\Sym(U_i)} \lra a^* \Fc|_{\Sym(U_1\cup U_2\cup U_3)}
 \]
 as a composition of two 2-variable $\mu$'s, are equal.
 
 \vskip .2cm 
 
 (b) Compatibility with restrictions. Suppose $U_1, \cdots, U_n$ are disjoint disks as above and $U'_i\subset U_i$
 are smaller disks. Then the restriction of $\mu_{U_1, \cdots, U_n}$ to $\bigboxtimes_{(U'_1, \cdots, U'_n)} \Fc|_{\Sym(U_i)}$
 is equal to $\mu_{U'_1, \cdots, U'_n}$
 
 \end{exas}
 
 We denote by $\FPS(\Vc)$ the category formed by factorizable perverse sheaves on $(\Sym(\CC),S)$ with
 values in $\Vc$. 
 
 \begin{rem}\label{rem:non-round}
 Note that the isomorphism \eqref{eq:add-map} holds for arbitrary, not necessarily round, disjoint topological disks
  $U_1, \cdots, U_m\subset \CC$ . One can show that for a factorizable perverse sheaf $\Fc$ 
  and any such disjoint $U_i$
  there is  a natural
  identification of $a^*\left(  \Fc|_{\Sym(\bigcup U_i)}\right)$ with ``the'' tensor product of the $\Fc|_{\Sym(U_i)}$.
  To formulate it canonically, we need a version of braided  tensor product labelled by arbitrary sequence
  of disjoint topological disks. This can be achieved by forming an extended, colored version of the $E_2$-operad
  where each topological disk $U\subset \CC$ is a separte color, cf. \cite{CG} \S 3.1.2.  Since this generalization is not
  necessary for us, we do not pursue it here. 
  \end{rem}
 
 
 \subsection{From bialgebras to factorizable sheaves: the localization theorem }\label{par:loc-fun}
 
 \paragraph{Statement of results.} 
 
 Let $\Br_n$ be the Artin braid group on $n$ strands. 
 Let 
 \[
 j_n: \Sym^n_\neq (\CC) = X_{1^n}^\CC \hookrightarrow  \Sym^n(\CC)
 \]
  be the embedding of the open stratum (polynomials with distinct roots). 
 Its fundamental group is $\Br_n$. 
 
 As before, let $\Vc$ be a $\k$-linear braided monoidal abelian category with $\otimes$ bi-exact.  Recall that
 we  have extended the braided monoidal structure to the bounded derived category $D^b(\Vc)$. 
 
 For any object $V\in D^b(\Vc)$ the tensor power $V^{\otimes n}$ a representation of $\Br_n$ and so gives a 
 $D^b(\Vc)$-valued local system
 on $\Sym^n_\neq (\CC)$ which we denote  $\Lc_{n}(V)$. In the language of Proposition \ref{prop:braid-e2}, various determinations of
 $V^{\otimes n}$ form, by the very definition, a $D^b(\Vc)$-valued
 local system on the space $E_2(n)$ which is homotopy equivalent to 
 $\Sym^n_\neq(\CC)$.  This is $\Lc_{n}(V)$. We will view $\Lc_n(V)$ as a complex of locally constant sheaves on
 $\Sym^n_\neq(\CC)$. 
  Recall also the notation $K_n = \Sym^n(\RR)$ for the real skeleton of $\Sym^n(\CC)$. 
 
 \vskip .2cm

 Let $\PB(\Vc)$ be the category of primitive bialgebras in $\Vc$. 
 
 \begin{thm}\label{thm:localize}
 There is an equivalence of categories (called the  {\em localization functor}) 
 \[
 L: \PB(\Vc) \lra \FPS(\Vc), \quad A \mapsto L(A) \,=\,\bigl(L_n(A)\in \Perv(\Sym^n(\CC), S_n; \Vc) \bigr)
 \]
 with the following properties:
 \begin{itemize}
 
 \item[(a)] The restriction $L_n(\Fc)|_{\Sym^n_\neq(\CC)}$ is identified with the local system
 $\Lc_{n}(A_1[1])$ of determinations of
 $(A_1[1])^{\otimes n}$.
 
 \item[(b)] The total space of vanishing cycles $\Phi_\tot(L_n(A))$, see \eqref {eq:total-van}, is identified with $A_n$.
 
 \item[(c)] More generally, the cubical diagram associated to the constructible sheaf $R_{K_n}(L_n(A))$
 is identified with $\BB_n^\dagger(A)$, the $n$th cobar-cube of $A$. 
 
 \item[(d)] The stalk if $L_n(A)$ at $0$ is identified with the $n$th bar-complex of $A$
 \[
 B_n^\bullet(A) \,=\,\biggl\{ A_1^{\otimes n} \to \bigoplus_{l(\alpha)=n-1} \bigotimes_{i=1}^{n-1} A_{\alpha_i} \to 
 \cdots \to  \bigoplus_{\alpha_1+\alpha_2=n} A_{\alpha_1}\otimes A_{\alpha_2} \to A_n\biggr\},
 \]
 with grading normalized so that $A_1^{\otimes n}$ is in degree $(-n)$.
 
 \item[(e)] Suppose $\Vc$ is rigid. Then for the primitive bialgebra $A^*$ dual to $A$ we have that  $L_n(A^*) = \DD(L_n(A))$
 is the Verdier dual to $L_n(A)$. 
  \end{itemize}
 \end{thm}
 
 \begin{rem}\label{rem:tw-sign}
 In part (a) of the theorem, $(A_1[1])^{\otimes n}$ is identified with $A_1^{\otimes n}[n]$, that is,
 the $n$th tensor power of $A_1$ put in degree $(-n)$. However, the braid group action on $(A_1[1])^{\otimes n}$
 differs from that on $A_1^{\otimes n}$ by twisting with the sign character $\Br_n\to\Sigma_n \buildrel\on{sgn}\over\to
 \{\pm 1\}$.  This is because we extend the braided structure on $\Vc$ to the derived category of $\Vc$
by  using the Koszul sign rule \eqref{eq:koszul}.

 \end{rem}
 
 We denote $\Sym_\neq(\CC) = \bigsqcup_n \Sym^n_\neq (\CC)$ and write $j: \Sym_\neq(\CC)\hookrightarrow \Sym(\CC)$
 for the embedding. For an object $V\in\Vc$ we denote $\Lc (V[1])$ the local system on $\Sym_\neq(\CC)$ which on $\Sym^n_\neq(\CC)$ is
 $\Lc_{n}(V[1])$. 
 
 \begin{thm}\label{thm:shuffle}
 Let $V$ be any object of $\Vc$. Then:
 
 (a) We have $ L(T_!(V))\,\simeq \, j_! \Lc  (V[1])$.
 
 \vskip .2cm
 
 (b) We also have $ L(T_*(V))\, \simeq \, j_* \Lc  (V[1])$
 
 \vskip.2cm
 
 (c) Further, 
 the canonical map $\alpha: T_!(V)\to T_*(V)$ is sent by $L$ to the canonical map $c: j_!\to j_*$, and
 the Nichols algebra $T_{!*}(V)$  corresponds, under $L$, to the intersection cohomology extension $j_{!*}\Lc (V[1])$.
 \end{thm}
 
 Theorems \ref{thm:localize} and \ref{thm:shuffle} will be proved in the next Ch. \ref{sec:cousin}.
 Here we note the following corollary which was proved, in an equivalent form, in \cite{westerland}, Th. 1.3. 
 
 \begin{cor}\label{cor:braid-cohom}
 We have an identification
 \[
 H^j(\Br_n, V^{\otimes n}\otimes_\k\on{sgn}) \,\simeq \, \Tor_{n-j,n}^{T_*(V)} (\1, \1). 
 \]
 \end{cor}
 
 \noindent {\sl Proof of the corollary:} Let $j_n: \Sym^n_\neq(\CC)
 \to\Sym^n(\CC)$ be the embedding. 
 By  \ref{thm:shuffle} (b)  
 $L_n(T_*(V)) \simeq  j_{n*} \Lc_n(V[1])$. Because of the quasi-homogeneity of the diagonal stratification $S$ of
 $\Sym^n(\CC)$ (the $\CC^*$-action coming from the standard dilation action on $\CC$),
 the stalk of 
 $ j_{n*} \Lc_n(V[1])$ at $0$ is identified with
  the complex 
  \[
  R\Gamma(\Sym^n(\CC),  j_{n*} \Lc_n(V[1]))\,\simeq \, R\Gamma(\Sym^n_\neq(\CC), \Lc_n(V[1]))
  \]
 and the  complex in the RHS calculates $H^\bullet(\Br_n, V^{\otimes n} \otimes_\k\sgn)[n]$.
 So our statement follows from Theorem  \ref{thm:localize}(d). \qed
 
 \paragraph{Example: Yetter-Drinfeld modules and Hurwitz spaces.}
 Let us explain the relation of our results with the work 
 \cite{westerland} of Ellenberg, Tran and Westerland. 
 
 \vskip .2cm
 
 Let $G$ be a finite group with multiplication $m: G\times G\to G$.
  A {\em Yetter-Drinfeld} $G$-module (YD $G$-module, for short)
 is a finite-dimensional $G$-graded $\k$-vector
 space $V=\bigoplus_{g\in G} V_g$ with a right $G$-action such that $V_g \cdot h  = V_{h^{-1}gh}$; this is the same as a module over the Drinfeld double of the group 
ring $\k[G]$. 
 Such objects form a  $\k$-linear abelian rigid  braided monoidal category $\Yc_G$ with $\otimes$
 being the usual graded tensor product and braiding given by
 \[
 R(v\otimes w) \,=\, w\otimes (v\cdot h), \quad v\in V_g, w\in V_h. 
 \]
 Alternatively, we can consider {\em character sheaves}, i.e., conjugation equivariant sheaves $F$ of
 finite-dimensional $\k$-vector spaces on $G$ as a discrete space.  
 Such a sheaf $F$ gives a YD $G$-module
 $V=\Gamma(G,F)$ of global sections, and the monoidal operation corresponds to the convolution of
 sheaves $F*F'=m_*(F\boxtimes F')$. 
 The braided category of character sheaves
 (i.e., the category $\Yc_G$) can be, therefore, seen as a categorification of the commutative ring of
 class functions on $G$ under convolution. See \cite{bakalov} \S 3.2, 
\cite{lusztig:exotic} for more background. 
 
 \vskip .2cm

 If $C\subset G$ is a conjugation invariant subset, then the constant sheaf $\ul\k_C$ is a character sheaf
 with corresponding YD $G$-module $V=\k C$ (the $\k$-span of $C$). 
 In this case the $\Br_n$-module $V^{\otimes n} = \k C^n$ is the permutation module corresponding to the $\Br_n$-action on
 the Cartesian power $C^n$ given on the standard generators  by
 \[
 \sigma_i(g_1,\cdots, g_n) \,=\,\bigl( g_1, \,\cdots, \,g_{i-1}, \,g_{i+1}, g_{i+1}^{-1} g_i g_{i+1}, \,g_{i+2}, \,\cdots, \,g_n\bigr),
 \quad i=1,\cdots, n-1. 
 \]
 (This action reflects the embedding of $\Br_n$ into the automorphism group of the free group.)
 The {\em Hurwitz space} associated to $G$, $C$ and $n$ is the unramified covering
 \[
 p_n: \Hur^C_{G,n}\lra \Sym^n_\neq(\CC)
 \]
 associated to the $\Br_n$-set $C^n$. It has the following relations to moduli of  branched covers 
 (see \cite{ellenberg-venkatesh} \S 2 for a systematic discussion):
 
 \begin{itemize}
 \item[(1)] Let $D=\{|z|<1\}$ be the unit disk in $\CC$ and $\ol D$ be its closure $|z|\leq 1$. The complex manifold 
 \[
 \Hur^C_{G,n}(D) = p_n^{-1}(\Sym^n_\neq(D))
 \]
parametrizes 
  $G$-coverings $\pi: Y\to \ol D$ ramified at some $n$ points in $D$, with monodromies around these
 points liying in $C$, together with a distinguished point $y\in \pi^{-1}(1)$.
 
 \item[(2)] The quotient Deligne-Mumford stack $\Hur^C_{G,n}/\!/G$ parametrizes
   $G$-coverings of the 
   affine line $\AAA^1$ ramified at $n$ points with monodromies around these points
 lying in $C$. 
  \end{itemize}
  
  Note that the intrepretation (1) gives a map $M:  \Hur^C_{G,n}(D)\to G$ given by the counterclockwise
  monodromy around the boundary of $\ol D$, and so
  we have a disjoint union decomposition ($G$-grading)
  \[
   \Hur^C_{G,n}(D)\,=\, \bigsqcup_{g\in G}  \Hur^C_{G,n}(D)_g, \quad  \Hur^C_{G,n}(D)_g :=M^{-1}(g).
  \]
  Together with the natural $G$-action on $  \Hur^C_{G,n}(D)$ (change of $y$) this makes  $\Hur^C_{G,n}(D)$
  into a ``YD $G$-module in the category of complex manifolds". 

 \vskip .2cm

 Thus the (untwisted) local system $\Lc_n(V)$ of determinations of $V^{\otimes n}$ is the direct image
 $p_{n*} \ul\k _{ \Hur^C_{G,n}}$. By introducing a minus sign into the braiding, the authors
 of \cite{westerland} realize $V^{\otimes n}\otimes \sgn$ as $p_{n*} \ul\k _{ \Hur^C_{G,n}}$,
 which, 
 by Corollary  \ref{cor:braid-cohom} (i.e., Thm. 1.3 of \cite{westerland}) identifies the topological cohomology  of Hurwitz spaces as
 \[
 H^j( \Hur^C_{G,n}, \k) \,\simeq \,  \Tor_{n-j, n}^{T_*(V)}(\k, \k). 
 \]

 We note that Theorem \ref{thm:shuffle}(c) 
  provides a similar geometric intnerpretation of the Nichols algebra
 $T_{!*}(V)$:

 \begin{cor}
 We have an identification
 \[
 \Tor_{n-j,n}^{T_{!*}(V)}(\k, \k) \,=\, \HH^j(\Sym^n(\CC), j_{!*} p_{n*} \k_{ \Hur^C_{G,n}}).  \qed
 \]
 \end{cor}
 It is interesting whether the RHS of this identification, involving the intersection
 cohomology extension of a local system, can be described in terms of   intersection cohomology of some
  natural compactification
 of the Hurwitz space. 
 
 \begin{rem}
 The formulation of \cite{westerland} proceeds in the dual form,
 using $\Ext$-{\em cohomology},  not $\Tor$-homology of $T_*(V)$,
 which is related to the {\em homology} rather than cohomology of $\Hur^C_{G,n}$.  This has the convenience
 that $\Ext^\bullet(\k,\k)$ forms an algebra, rather than a coalgebra. The fact that this algebra is (graded)
 braided-commutative
 reflects the phenomenon, mentioned in the Introduction,  that the Koszul dual of a Hopf algebra is an $E_2$-algebra.
 In fact, the appearance of  $E_2$-algebras in this context is more fundamental, since:
 
 \begin{itemize} 
 
 \item[(1)] The Hurwitz
 spaces $\Hur^C_{G,n}(D)$ themselves, taken for all $n\geq 0$, form, naturally, an $E_2$-algebra in the 
 (braided monoidal) category of YD-modules in complex manifolds. 
 
 \item[(2)] Similarly, 
 the quotients
 $\Hur^C_{G,n}(D)/\!/G$ form an $E_2$-algebra in the more usual (symmetric monoidal) category of  analytic stacks.
 
 \end{itemize}
 These algebra structures are given by  appropriate versions of the gluing maps of \cite{ellenberg-venkatesh} \S 2.3.

 \end{rem}
 
 \vfill\eject
 
 \section{ Proof of Localization Theorem: Cousin complexes}\label{sec:cousin}

  \subsection{ The Cousin  sheaves  of a coalgebra}\label{subsec:cousin-coalg}
  
  Our proof of Theorem \ref{thm:localize} consists in constructing the perverse sheaf $L_n(A)$
  as an explicit complex of sheaves on $\Sym^n(\CC)$, which we call the Cousin complex.
  In this section we describe the sheaves constituting this complex.

 \paragraph{ Objects $A_\bz$ and local systems $\Lc_\alpha(A)$.}
 Let $A = \bigoplus_{n=0}^\infty$ be a graded object in $\Vc$  with $A_0=\1$.
 Let  $\bz\in \Sym^n(\CC)$ be a point understood as an unordered divisor
 $\bz = \sum_{i=1}^d \lambda_i z_i$ with $\lambda_i\geq 1$ and $z_i\neq z_j$. The operadic point of view on
 the braided structure in $\Vc$ allows us to form the object
 \[
 A_\bz\,=\,\bigotimes\nolimits_{(z_1, \cdots, z_d)} A_{\lambda_i},
 \]
 which is canonically (up to a unique isomorphism) independent on the ordering of $z_1,\cdots, z_d$. 
 Given an unordered partition $\lambda = (\lambda_1, \cdots, \lambda_d)\in\Pc(n)$, the objects $A_\bz$, $\bz\in X_\lambda^\CC$,
 form a local system on the complex stratum $X_\lambda^\CC$. We denote this  local system $\Lc_\lambda(A)$. 
 For $\lambda=1^n$ we get the local system $\Lc_n(A_1)$ of determinations of $A_1^{\otimes n}$ on $\Sym^n_\neq(A)$,
 see \S \ref{subsec:fact-perv}\ref{par:E2-br}
 
 \paragraph{Cousin sheaves $\Ec^\alpha(A)$ of a coalgebra.}\label{par-cous-sh}
 
 Let us further assume
  that $A$ is a graded  coassociative coalgebra in $\Vc$ with counit given by the projection to $A_0=\1$,  satisfying the primitivity
   condition (2) of Definition 
 \ref{def:prim-gen-bialg}. 
 
 Since the objects $A_\bz$ are canonically determined by the points  $\bz\in\Sym^n(\CC)$,
 it is natural to consider the following.

\begin{prob}\label{prob:sheaf}
 Arrange the objects $A_\bz$ into a constructible sheaf on $\Sym^n(\CC)$  so that these objects
 are the stalks and the sheaf structure (generalization maps) is given by the comultiplication $\Delta$ in $A$. 
 \end{prob}

 However, this is impossible unless $\Delta$ is cocommutative. 
 
 \begin{ex}
 Let $n=2$. After quotienting by translations, $\Sym^2(\CC)$ becomes the complex plane, with the open stratum,
 the reduction of $X_{1^2}^\CC = \Sym^2_\neq(\CC)$ being $\CC^*$ and the closed stratum, the reduction
 of $X_{(2)}^\CC$, being the point $0$. So the local system $\Lc_{1^2}(A)$ on $\CC^*$ has monodromy
 $R: A_1\otimes A_1\to A_1\otimes A_1$, and the   stalk of $\Lc_{(2)}(A)$ at $0$ is $A_2$. 
 To combine them into a constructible sheaf on $\CC$, we need to specify a morphism 
 $\eps: A_2\to A_1\otimes A_1$
 such that $R\circ\eps = \eps$,
  so the comultiplication cannot do this in general. 
 \end{ex}
 
 Nevertheless we have interesting {\em partial} solutions to Problem \ref{prob:sheaf}, which produce sheaves
 not on all $\Sym^n(\CC)$ but on some subspaces.
 
 \begin {prop}
 Let $\alpha\in\2^{n-1}$. The comultiplication $\Delta$ in $A$ defines a sheaf $\wt\Ec^\alpha(A)$
 on $X_\alpha^\Ien$ whose stalk at any $\bz\in X_\alpha^\Ien$  is identified with $A_\bz$
  and the generalization maps are given by  the appropriate components of $\Delta$. 
 \end{prop}
 
  \begin{rem}
 The reason why it becomes possible to arrange the $A_\bz$ into a sheaf on $X_\alpha^\Ien$ is
 that the intersection of each $X_\alpha^\Ien$ with
 each $X_\lambda^\CC$ is   contractible, see Proposition \ref{prop:fine-cells}. 
    \end{rem}
 
 \noindent {\sl Proof of the proposition:} 
  Recall the ``real parts'' projection 
\[
\rho_\alpha: X_\alpha^\Ien \lra \prod_{\nu=1}^k \Sym^{\alpha_\nu}(\RR)  
\]
  For each $m\geq 0$ let $\BB^\dagger_m(A)$ be the $m$th cobar-cube
 associated to the coalgebra structure on $A$ and let 
 $R_m(A)$ be the constructible sheaf on
 the stratified space $(K_m=\Sym^m(\RR), K_{m,S})$ corresponding to  $\BB^\dagger_m(A)$. Define, using the
 ``classical'' point of view (one tensor product operation)  on the braided structure: 
 \be\label{eq:wt-Ea}
 \wt\Ec^\alpha(A) =  \rho_\alpha^* \biggl(\bigboxtimes_{\nu=1}^k  R_{\alpha_\nu}(A)\biggr). 
  \ee
  Let us show why, at the level of stalks, this gives the desired answer. 
   
   \vskip .2cm
  
Let first $m\geq 0$ and consider a point $\bf x$ of $\Sym^m(\RR)$ written  as an ordered divisor $m_1 x_1 + \cdots + m_p x_p$,
with $x_1<\cdots x_p$ and $m_1+\cdots +m _p=m$. The stalk of $R_m(A)$ at $\bf x$ is, in the ``classical'' notation
\be
R_m(A)_{\bf x} \,=\, A_{m_1}\otimes \cdots A_{m_p}.
\ee
Let now $\alpha=(\alpha_1, \cdots, \alpha_k)\in\2^{n-1}$ and  ${\bf z}  = \{z_1, \cdots, z_n\}  \in X_\alpha^\Ien$.
In other words, according to the increasing order of the  $\Im(z_j)$,
the $z_j$ are subvidided into $k$ groups, so that the elements  the $\nu$th group, of cardinality $\alpha_\nu$, have the same
 imaginary part. The $\nu$th component of the projection $\rho_\alpha$ is the map $X_\alpha^\Ien\to\Sym^{\alpha_\nu}(\RR)$
 which sends $\bf z$ into the points represented by the real parts of the $z_j$ from the $\nu$th group.
 Let us now write $\bf z$ itself as a divisor $\sum_{i=1}^d \lambda_i z_i$ with $z_1,\cdots, z_d$ being distinct and 
 $\lambda_i >0$. 
  From the above we conclude that in the ``classical''  notation,  the stalk of $\wt\Ec^\alpha(A)$ at $\bf z$  would be 
   simply the tensor product
 $\bigotimes_{i=1}^d  A_{\lambda_i}$, taken in the lexicographic order with respect first, to the imaginary and  then to the real parts
 of $z_1,\cdots, z_d$. This is  precisely $A_\bz$. 
 
  \qed

    Let us  form the   complexes  (recall that $j_{\alpha *}$ means the derived direct image)
  \be\label{eq:cousin-sh}
   \Ec^\alpha(A) = j_{\alpha *}\, \wt\Ec^\alpha(A), \quad  j_\alpha: X_\alpha^\Ien  \hookrightarrow 
\Sym^n(\CC). 
\ee
We call $\Ec^\alpha(A)$ the $\alpha$th {\em Cousin sheaf} of $A$; 
it will be proved later that is indeed quasi-isomorphic to a single sheaf in degree $0$.
In order to establish this and other properties of $\Ec^\alpha(A)$ we need some preparations. 


\paragraph{Contingency matrices and double cosets.}\label{par:stoch-mat}
  Let $\xi=(\xi_1,\cdots, \xi_r)$ and $\eta = (\eta_1,\cdots, \eta_s)$
be two ordered partitions of the same number $m = \sum \xi_i = \sum\eta_j$. 
By a $(\xi,\eta)$-{\em contingency matrix}\footnote{the statisticians call such matrices {\it contingency tables with margins $(\xi,\eta)$} 
(the name given by Karl Pearson at the beginning of the last century), cf. \cite{diaconis}.} we will mean an $r\times s$ matrix 
$\zeta = \|\zeta_{ij}\|_{i=1,\cdots, r}^{j=1,\cdots, s}$ 
  with entries in $\ZZ_{\geq 0}$ such that
\[
\sum_j \zeta_{ij} = \xi_i, \quad \forall i=1,\cdots, r,\quad\quad \sum_i \zeta_{ij} = \eta_j,  \quad\forall j=1,\cdots, s,
\]
and, in addition, each row and each column  has a nonzero entry. 
We denote the set of $(\xi,\eta)$-contingency matrices by $S(\xi,\eta)$.
 
\begin{ex}\label{ex:1n-stoch}
Let $\eta = 1^m$ be the maximal ordered partition. Then a $(\xi, 1^m)$-contingency matrix is the same as
a shuffle of type $\xi$, i.e., a sequence $(X_1, \cdots, X_r)$ of  subsets of $\{1,\cdots, m\}$ with
$|X_i|=\xi_i$, forming a disjoint decomposition. In particular,
\[
|S(\xi, 1^m)| \,=\, {m! \over (\xi_1)! \, \cdots \, (\xi_r)!}.  
\]
 \end{ex}
 
 Let $F$ be a field and 
 Let $P(\xi)\subset GL_m(F)$ be the parabolic subgroup of block upper triangular matrices with blocks of
 sizes $\xi_1,\cdots, \xi_r$, and similarly for $P(\eta)$. Recall that $\Sigma_m$ denotes the symmetric group,
 i.e., the Weyl group of $GL_m(F)$. 
 Let  $\Sigma(\xi)= \Sigma_{\xi_1}\times\cdots\times \Sigma_{\xi_r}\subset \Sigma_m$ be the Weyl group
 of the Levi subgroup of $P(\xi)$, and similarly for $\Sigma(\eta)$.

\begin{prop}
We have the identifications (parabolic Bruhat decomposition)
 \[
 P(\xi) \backslash GL_m(F) / P(\eta)  \,\simeq \,  S(\xi,\eta)
  \,\simeq \, \Sigma(\xi) \backslash \Sigma_m / \Sigma(\eta). 
 \]
 \end{prop}
 
 \noindent {\sl Proof:} This is well known.  For a discussion of the first identification
 (possible relative positions of a pair of partial flags in $F^m$)  see, e.g., \cite{BLM}.
 As for the second identification, see \cite{diaconis} Lemma 3.3,  which has been visualized as
  the ``balls in boxes" picture in \cite{petersen} 
 \S 6.  
 
  \qed

 \vskip .2cm
 
 Given $\zeta\in S(\xi,\eta)$, we have an ordered partition $\lex(\zeta)$ of $m$ obtained by reading the
{\em nonzero}  entries of $\zeta$ in the lexicographic order, starting from the first 
column $\zeta_{1,1}, \cdots, \zeta_{s,1}$,
 then the second column $\zeta_{1,2},\cdots, \zeta_{s,2}$ and so on. 
 Note that $\lex(\zeta)>\eta$ is a refinement
 of $\eta$ (but not necessrily a refinement of $\xi$).

 \vskip .2cm
 
 Let now $n\geq 0$ and $\alpha,\beta,\gamma\in \2^{n-1}$ be three ordered partitions of $n$
 such that $\alpha\geq\beta \leq\gamma$.  Writing $\beta= (\beta_1,\cdots, \beta_k)$, the inequality $\beta\leq\gamma$
 means that each $\beta_i$ is the sum of a string of consecutive $\gamma$'s:
 \be\label{eq:beta-gamma}
 \beta_i \,=\,\gamma_{r_i} + \gamma_{r_i+1} + \cdots + \gamma_{r_{i+1}-1}, \quad
 \text{for some}\quad  1=r_1< \cdots < r_{k+1}= l(\gamma)+1. 
 \ee
 Similarly, the inequality $\alpha\geq\beta$ means that
 \be\label{eq:alpha-beta}
  \beta_i \,=\,\alpha_{s_i} + \alpha_{s_i+1} + \cdots + \alpha_{s_{i+1}-1}, \quad
 \text{for some}\quad  1=s_1< \cdots < s_{k+1}= l(\alpha)+1. 
 \ee
This means that for each $i$ we have two ordered partitions of $\beta_i$ which we denote
\[
\alpha/\beta_i \,=\,(\alpha_{s_i},  \alpha_{s_i+1},  \cdots,  \alpha_{s_{i+1}-1}),\quad
\gamma/\beta_i  \,=\, (\gamma_{r_i},  \gamma_{r_i+1},  \cdots,  \gamma_{r_{i+1}-1}). 
\]
We put
\[
S(\alpha, \beta, \gamma) \,=\,\prod_{i=1}^k S(\alpha/\beta_i, \gamma/\beta_i). 
\]
Thus, an element of $S(\alpha, \beta, \gamma)$ is a sequence $\vec{\zeta} = (\zeta^{(1)}, \cdots, \zeta^{(k)})$
of contingency matrices $\zeta^{(i)}\in S(\alpha/\beta_i, \gamma/\beta_i)$. For such  a sequence $\vec\zeta$
we denote by $\Lex(\vec{\zeta})$ the ordered partition of $n$ formed by writing $\lex(\zeta^{(1)}), \cdots,
\lex(\zeta^{(k)})$ in a single sequence. We note that $\Lex(\vec{\zeta})\geq \beta$ is a refinement of 
$\beta$ (but not necessarily a refimenent of $\alpha$).


\paragraph{Properties of the Cousin sheaves.}
Recall that $A=\bigoplus_n A_n$ is a graded coalgebra in the braided category $\Vc$, as in
\S \ref{par-cous-sh}  We study the complexes $\Ec^\alpha(A)$ defined in \eqref{eq:cousin-sh}. 

\vskip .2cm

For any  ordered partition $\theta = (\theta_1,\cdots, \theta_p)\in\2^{n-1}$ we denote $A_\theta = A_{\theta_1}\otimes\cdots\otimes A_{\theta_p}$.  

\vskip .2cm

For two ordered partitions $\xi,\eta\in \2^{m-1}$ and a contingency matrix 
$\zeta = \|\zeta_{ij}\|\in S(\xi,\eta)$
we denote $A_\zeta = \bigotimes_{i,j} A_{\zeta_{ij}}$, the tensor product in the sense of
\eqref{eq:braided-tensor}, formed with respect to
positioning $A_{\zeta_{ij}}$ at the point $(i,j)\in\RR^2=\CC$, i.e., at the complex number $i+j \sqrt{-1}$. 

\vskip .2cm

For $\alpha\geq\beta = (\beta_1,\cdots, \beta_k) \leq\gamma \in \2^{n-1}$ and a sequence of matrices
 $\vzeta = (\zeta^{(1)}, \cdots, \zeta^{(k)})\in S(\alpha,\beta,\gamma)$ we put 
 $A_\vzeta = A_{\zeta^{(1)}}\otimes\cdots
 \otimes A_{\zeta^{(k)}}$.

\begin{prop}\label{prop:Ea-stalks}
(a) The restriction of $\Ec^\alpha(A)$ to any fine cell $X_{[\beta:\gamma]}$, $\beta\leq\gamma$,
vanishes unless $\alpha\geq\beta$. 

\vskip .2cm

(b) Suppose $\alpha\geq\beta\leq\gamma$. Then the stalk of $\Ec^\alpha(A)$ at any point
 $\bz\in X_{[\beta:\gamma]}$ is quasi-isomorphic to the object
 \[
 \bigoplus_{\vec{\zeta}\in S(\alpha,\beta,\gamma)} A_{\vec{\zeta}}. 
 \]
 in degree $0$. 
\end{prop}

This implies immediately:

\begin{cor}
$\Ec^\alpha(A)$ is quasi-isomorphic to a single sheaf in degree $0$ which is constructible with respect
to the stratification by fine cells. \qed
\end{cor}

To prove  Proposition \ref{prop:Ea-stalks},
let $\bz\in X_{[\beta:\gamma]}$ and $U\subset\Sym^n(\CC)$ be a small ball around $\bz$. The stalk of
$\Ec^\alpha(A)$ at $\bz$ is, by definition,  $R\Gamma(U\cap X_\alpha^\Ien, \wt\Ec^\alpha(A))$. 
So it is enough to show the following.

\begin{prop}\label{prop:U-intersect}
\begin{itemize}
\item[(1)] Unless $\alpha\geq\beta$, the intersection $U\cap X_\alpha^\Ien$ is empty.

\item[(2)] Assume $\alpha\geq\beta$. Then:

\begin{itemize}
\item[(2a)] $U\cap X_\alpha^\Ien$ is the union of contractible components $T_{\vec{\zeta}}$ labelled
by $\vec{\zeta}\in S(\alpha,\beta,\gamma)$. 

\item[(2b)] Further, consider the stratification of $X_\alpha^\Ien$ given by 
$X_\alpha^\Ien = \bigsqcup_{\eps\geq\alpha}
X_{[\alpha:\eps]}$. The intersection of $T_{\vec{\zeta}}$ with each  stratum $X_{[\alpha:\eps]}$ is contractible or empty,
and the minimal stratum meeting $T_{\vec{\zeta}}$ is $X_{[\alpha: \Lex(\vec{\zeta})]}$. 
\end{itemize}
\end{itemize}
\end{prop}

\noindent {\sl Proof:} 
Part (1) is clear since the closure of $X_{\alpha}^\Ien$ is the union of $X_{[\beta:\gamma]}$ for
 $\alpha\geq\beta\leq\gamma$. 
 
 \vskip .2cm
 
 Let us prove part (2). 
Since $\bz\in X_{[\beta:\gamma]}\subset X_\beta^\Ien$, we have a decomposition
$\bz=\bz_1+\cdots + \bz_k$ where  each $\bz_i\in\Sym^{\beta_i}(\CC)$ has the imaginary parts of
all its component points equal;   denoting this common value    $\Im(\bz_i)$, we further have
$\Im(\bz_1) < \cdots < \Im(\bz_k)$.  Further, the condition that $\bz\in X_{[\beta:\gamma]}$ means that
each $\bz_i$, positioned on a horizontal line $\Im(z)=\on{const}$, has multiplicities
$\gamma_{r_i}, \gamma_{r_i+1}, \cdots, \gamma_{r_{i+1}-1}$, see \eqref{eq:beta-gamma}, 
if read from left to right on that line. 

Now suppose we deform $\bz$ to a point $\bz'$ in a small ball $U$ around $\bz$, such that
$\bz'\in X_\alpha^\Ien$. Recall that $\alpha\geq\beta$, in particular, $\beta_i$ is represented as the sum
\eqref{eq:alpha-beta}. 

\vskip .2cm

Such a deformation  deforms each $\bz_i$ to some $\bz'_i$.
The deformation of $\bz_i$ to $\bz'_i$ subdivides each $\bz_i$ to a sum $\bz_{i, s_i} + \bz_{i, s_i+1} + \cdots + \bz_{i, s_{i+1}-1}$
so that $\bz_{i,\nu}$ moves to  a divisor $\bz'_{i,\nu}$ elying on the  $\nu$th horizontal line $\Im(z)=\on{const}$ featured in the definition of
 the stratum $X_\alpha^\Ien$.  
Such a separation defines a contingency matrix $\zeta^{(i)}\in S(\alpha/\beta_i, \gamma/\beta_i)$.
This matrix is just the matrix of multiplicities of the divisor $\bz'_i$. So the tensor product $A_{\bz'_i}$ is
identified with $A_{\zeta^{(i)}}$.

\vskip .2cm

Now,  deforming the whole of $\bz$ to $\bz'$  amounts to deforming each $\bz_i$ to a $\bz'_i$ as above, and this
defines a sequence of contingency matrices $\vzeta = (\zeta^{(1)}, \cdots, \zeta^{(k)})$.
The component $T_\vzeta$ of the intersection $U\cap X_\alpha^\Ien$ consists, by definition,
 of deformations $\bz'$ giving the sequence
$\vzeta$.  The space $\Gamma(T_\vzeta, \wt\Ec^\alpha(A))$ is then identified with the
 tensor product $A_{\bz'}$ for any $\bz'\in T_\vzeta$ , and this tensor product is  identified, by the above, with $A_\vzeta$.


    \subsection{The Cousin complexes of a bialgebra as  perverse sheaves}\label{subsec:cousin-perv} 
    
    \paragraph{Definition of  the Cousin complexes.}
 
  Assume now that $A$ is a primitive bialgebra in $\Vc$.   For $\alpha,\beta\in\2^{n-1}$  and $d>0$ we write $\alpha>_d\beta$ if $\alpha > \beta$ and 
 $l(\alpha)=l(\beta)+d$. 

  \vskip .2cm

 Let $\alpha = (\alpha_1, \cdots, \alpha_{k+1}) >_1\beta = (\beta_1,\cdots, \beta_k)$. Define the morphism of sheaves on $X_\beta^\Ien$
 \[
 \delta'_{\alpha,\beta}: j_\beta^*\Ec^\alpha(A) = j_\beta^*  j_{\alpha*}\wt\Ec^\alpha(A) \lra \wt\Ec^\beta(A)
 \]
 using the multiplication $\mu$ in $A$. More precisely, $\alpha>_1\beta$ means that $\beta$ is obtained
 from $\alpha$ by combining together two consecutive parts, say $\alpha_i$ and $\alpha_{i+1}$:
 \be\label{eq:alpha-and-beta}
 \beta \,=\, (\alpha_1, \cdots, \alpha_{i-1}, \alpha_i+\alpha_{i+1}, \alpha_{i+2}, \cdots, \alpha_{k+1}),
 \ee
 which is a particuilar case of \eqref{eq:alpha-beta}. 
 Accordingly, we apply a particular case of the reasoning of Proposition \ref{prop:U-intersect}(2a). 
 That is, 
  let $\bz$ be a point of  $X_\beta^\Ien$,
 considered as an unordered divisor. The fact that $\bz\in X_\beta^\Ien$ means that 
 $\bz$ is decomposed into a sum 
 $\bz_1+\cdots +\bz_{k}$,
 where $\bz_j$ is of degree $\beta_j$ and its points have the same imaginary part. 
 The stalk of $j_\beta^* j_{\alpha *}\wt\Ec^\alpha$(A) at $\bz$ is, by definition,
 the space of sections of $\wt\Ec_\alpha(A)$ on $U\cap X_\alpha^\Ien$, where $U$ is a small ball
 around $\bz$ in $\Sym^n(\CC)$. This intersection consists of several components, labelled by
 the ways of decomposing $\bz_i$ into a sum of effective divisors 
 \[
 \bz_i^1 + \bz_i^2,\quad \bz_i^{1}\in \Sym^{\alpha_i}(\CC), \,\, \bz_i^{2}\in\Sym^{\alpha_{i+1}}(\CC). 
 \]
 The component, denote it 
 $T = T(\bz_i^1, \bz_i^2)$,
   corresponding
 to such a decomposition,  is obtaind by increasing the imaginary parts of  all the points of  $\bz_i^1$ by the small
 increment $\eps$, getting a new divisor $\bw_i^1$ and similarly decreasing the imaginary parts of the points of $\bz_i^2$
 getting a divisor $\bw_i^2$. The divisor
 \[
\bw =  \bz_1 + \cdots + \bz_{i-1} + \bw_i^1 + \bw_i^2 + \bz_{i+1} +\cdots + \bz_{k-1}
 \]
represents a point in  $T$ mapped by $\rho_\alpha: X_\alpha^\Ien\to\prod \Sym^{\alpha_i}(\RR)$ into the stratum which is minimal possible for T and so
the space of sections of  the sheaf $\wt\Ec^\alpha(A)$  over $T$ is  canonically identified with $A_\bw$. 
The multiplication $\mu$ in $A$ gives a morphism   $A_\bw\to A_\bz$. 
We define the action of $\delta'_{\alpha,\beta}$ on the stalk of $ j_\beta^*  j_{\alpha*}\wt\Ec^\alpha(A)$
at $\bz$ to be the map  $A_\bw\to A_\bz$
 given by multiplication in $A$.
 
 \begin{prop}\label{prop:delta-is-morph}
  The maps of stalks defined above commute with the generalization maps and so define
 a morphism of sheaves  $ \delta'_{\alpha,\beta}:   j_\beta^*  j_{\alpha*}\wt\Ec^\alpha(A) \lra \wt\Ec^\beta(A)$. 
  \end{prop}
  
  \noindent {\sl Proof:} The generalization maps are given by the comultiplication in $A$, the maps of stalks we defined
  are given by the multiplication. The commutation of the two sets of maps
  follows from the compatibility of multiplication and comultiplication in $A$.
  
  More precisely, let us write Definition \ref{def:br-bialg} of a braided bialgebra diagrammatically, as the commutativity of
  \be\label{eq:br-bialg}
  \xymatrix{
  A\otimes A
  \ar[dd]_{\mu}
   \ar[r]^{\hskip -0.5cm \Delta\otimes\Delta} & A\otimes A\otimes A \otimes A
  \ar[d]^{R_{23}}
  \\
  & A\otimes A\otimes A \otimes A
  \ar[d]^{\mu\otimes\mu}
  \\
  A\ar[r]^{\Delta}& A\otimes A,
  }
  \ee
  and decompose this condition into  parts corresponding to various graded components of $A\otimes A$.
  That is, we fix $(p,q), (r,s)\in \ZZ_+^2$  with $p+q=r+s$ and look at the matrix element of $\Delta\circ \mu$ (the lower left path in the diagram)
  mapping from $A_p\otimes A_q$ to $A_r\otimes A_s$. 
  
  Let $S((p,q), (r,s))$ be the set of $2\times 2$ matrices of nonnegative integers
  \[
  \begin{pmatrix} p_1&p_2 \\q_1&q_2 \end{pmatrix} \text{ such that }  p_1+p_2=p, \, q_1+q_2=q, \quad p_1+q_1=r,\, p_2+q_2=s.
  \]
    (This is a particular case of the definition of contingency matrices in \S \ref{subsec:cousin-coalg}\ref{par:stoch-mat})
  Then \eqref{eq:br-bialg} gives
  \be\label{eq:br-bialg-grad}
  \Delta_{r,s}\circ \mu_{p,q} \,=\, \sum_{\tiny{
    \begin{pmatrix} p_1&p_2 \\q_1&q_2 \end{pmatrix} }\in S ((p,q), (r,s))}  
  \bigl(\mu_{p_1, q_1}\otimes\mu_{p_2, q_2} \bigr) \circ  \bigl(\Id \otimes R_{A_{p_2}, A_{q_1}}\otimes \Id\bigr) \circ
 \bigl (\Delta_{p_1, p_2} \otimes \Delta_{q_1, q_2}\bigr). 
   \ee
   We now reduce our question 
 (that $\delta'_{\alpha,\beta}$ commutes with generalization maps) to an instance of \eqref{eq:br-bialg-grad}.
 Let $\alpha>_1\beta$ and $\bz\in X_\beta^\Ien$ be as before. Let us write the action of  $\delta'_{\alpha,\beta}$ 
 on the stalks at $\bz$ more explicitly.  By the above assumption \eqref{eq:alpha-and-beta}, we have
 $\beta_i=\alpha_i+\alpha_{i+1}$.  Let us write the summand $\bz_i$ of $\bz$ as
 \[
 \begin{gathered}
 \bz_i = m_1 z_{i,1} +\cdots + m_l z_{i,l}, \quad \sum m_\nu=\beta_i,
 \\
 \Im(z_{i,1})=\cdots=\Im(z_{i,l}), \quad \Re(z_{i,1}) <\cdots < \Re(z_{i,l}). 
 \end{gathered}
 \]
A decomposition $\bz_i=\bz_i^1+\bz_i^2$ defining a compoment $T(\bz_i^1, \bz_i^2)$ above, i.e.,
a summand in $(j_\beta^*\j_{\alpha*} \wt\Ec^\alpha(A))_\bz$, is given by writing, in some way,
\[
m_\nu = m_\nu^1 + m_\nu^2, \quad \nu=1,\cdots, l. 
\]
The corresponding summand is simply
 \[
\bigl(  \bigotimes_{m=1}^{i-1} A_{\bz_m} \bigr) \otimes \bigotimes\limits _{\nu=1}^l (A_{m^1_\nu} \otimes A_{m^2_\nu}) 
 \otimes \bigl(\bigotimes_{m=i+1}^{k} A_{\bz_m}\bigr),
 \]
the stalk  $\wt\Ec^\beta(A)_\bz$ is 
\[
\bigl(  \bigotimes_{m=1}^{i-1} A_{\bz_m} \bigr) \otimes 
\bigotimes\limits_{\nu=1}^l A_{m_\nu}\otimes 
\bigl(\bigotimes_{m=i+1}^{k} A_{\bz_m}\bigr). 
\]
and the map on the stalks  given by  $\delta'_{\alpha,\beta}$ is 
 \[
\xymatrix{ \bigotimes\limits _{\nu=1}^l (A_{m^1_\nu} \otimes A_{m^2_\nu}) 
\ar[rr]^{\hskip 0.8cm \mu\otimes\cdots\otimes\mu} &&\bigotimes\limits_{\nu=1}^l A_{m_\nu}
}
\]
tensored on the left and right by the identities. 
We need to show that such maps commute with the generalization maps for points inside 
the stratum $X_\beta^\Ien$ (i.e., from a smaller to a bigger fine cell inside this tratum). For this,
it suffices to consider {\em elementary} generalizations, when the dimension of the fine cell
increases by $1$. If  we start from $\bz$ as before, then an elementary generalization corresponds to
splitting just one $m_j$, in some way,  as $m'_j+m''_j$ and moving
 from $\bz$ to a nearby $\wt\bz$ which has the form
\[
\wt\bz = \bz_1+\cdots + \bz_{i-1} + \wt\bz_i + \bz_{i+1}+\cdots + \bz_k,
\]
where $\bz_t, t\neq i$ are the same as for $\bz$, while $\wt\bz_i$
is obtained by splitting some the divisor  $m_j z_{i,j}$ into $m'_j z'_{i,j}+ m''_j z''_{i,j}$ so that 
 \[
\Im(z'_{i,j})=\Im(z''_{i,j}) = \Im (z_{i,j}), \quad \Re(z'_{i,j}) < \Re(z''_{i,j}).
 \]
 So the number $m_j$ gets split in  two ways:
 \[
 m_j = m_j^1+m_j^2 = m'_j+m''_j. 
 \]
 We notice that all the tensor products involved will contain the tensor product of the $A_{\bz_m}$,
 $m\neq i$, and all the maps involved will be given by identities on these tensor factors.
 So in the further discussion we will ignore these factors. What  will remain is the diagrams of the form 
  \begin{tikzcd}
\bigotimes\limits_{\nu=1}^l  (A_{m_\nu^1}\otimes A_{m_\nu^2})  \arrow[r]
  \arrow[d,"\mu\otimes\cdots\otimes\mu" ]
&  \bigl(\bigotimes\limits_{\nu=1}^{j-1} A_{m_\nu^1}\otimes A_{m_\nu^2}\bigr)\otimes \bigl( A_{p_1}\otimes A_{p_2}\otimes
A_{q_1}\otimes A_{q_2}\bigr) \otimes\bigl (\bigotimes\limits_{\nu=j+1}^l  A_{m_\nu^1}\otimes A_{m_\nu^2}\bigr)
 \arrow[d] \\
 \bigotimes\limits_{\nu=1}^l A_{m_\nu}
  \arrow[r, " \Id\otimes\Delta_{m'_j, m''_j}\otimes\Id" ]
&  \bigl(\bigotimes\limits_{nu=1}^{j-1} A_{m_\nu}\bigr) \otimes (A_{m'_j}\otimes A_{m''_j}) \otimes 
\bigl(\bigotimes\limits_{nu=j+1}^l A_{m_\nu}\bigr) 
\end{tikzcd}

\vskip .2cm

 with 
 \[
 p_1+p_2=m_j^1, ,\, q_1+q_2=m_j^2, \quad p_1+q_1=m'_j, \, p_2+q_2=m''_j. 
 \]
 In each such  diagram, the left vertical arrow  represents (after omitting the inessential tensor factors, as we agreed)
  the action of $\delta$ on the stalks at $\bz$. 
 The lower horizontal arrow represents the generalization map of $\wt\Ec^\beta(A)$ 
 from $\bz$ to $\wt\bz$. 
 The top horizontal arrow, summed over the admissible $p_1, p_2, q_1, q_2$, represents
 the generalization from $\bz$ to $\wt\bz$ in $j_{\beta}^* j_{\alpha *} \wt\Ec^\alpha(A)$.
 The right vertical arrow, similarly summed, represents the action of $\delta'_{\alpha,\beta}$
 on stalks at $\wt\bz$. So removing even more inessential factors (corresponding to $\nu\neq j$),
 we reduce our statement to  \eqref{eq:br-bialg-grad}.
  Proposition \ref{prop:delta-mor} is proved. \qed
  
  \vskip .3cm
  
  By adjunction, the morphism $\delta'_{\alpha,\beta}$ gives a morphism of sheaves on $\Sym^n(\CC)$
  \[
  \delta_{\alpha, \beta}:  \Ec^\alpha (A)\lra \Ec^\beta(A), \quad \alpha >_1\beta. 
  \]
  
    \begin{prop}\label{prop:delta-mor}
    The morphisms $ \ul \delta_{\alpha,\beta}$ define an anti-representation of $\2^{n-1}$ in the category of sheaves on
  $\Sym^n(\CC)$, and therefore we have a complex of sheaves
 \[
\Ec^\bullet_n(A) \,=\,\biggl\{ \Ec^{(1^n)} (A) \buildrel\ul \delta\over\to \bigoplus_{l(\alpha) = n-1} \Ec^\alpha(A) \buildrel
\ul \delta\over\to \cdots \buildrel \ul\delta\over\to \Ec^{(n)}(A)\biggr\}.
\]
 Set the  grading of this complex so that $\Ec^{1^n}(A)$ is situated in degree $(-n)$. 
  \end{prop}
  
   \noindent {\sl Proof:} Follows from the associativity of the multiplication. 
   
   \qed

We call $\Ec_n^\bullet(A)$ the $n$th {\em Cousin complex} of $A$. By definition, it is a complex of sheaves
on $\Sym^n(\CC)$ constructible with respect to the stratification by fine cells. 

\begin{ex}
 Let  $n = 2$. The complex
  \[
 \Ec_2^\bullet(A)\,=\,\biggl\{ \Ec^{(1,1)}(A)\buildrel\ul\delta \over \lra \Ec^{(2)}(A)\biggr\}
 \]
 consists of $2$ terms. After  ``reducing" (quotienting) by translations, it can be seen as a complex of sheaves on $\CC$
 stratified by $\CC\setminus \RR_{\geq 0}$, the reduction of $X_{(1,1)}^\Ien$ and $\RR_{\geq 0}$, the
 reduction of $X_{(2)}^\Ien$. The sheaf $\Ec^{(2)}(A)$ lives on $\RR_{\geq 0}$ with stalks  (at $0$ and elsewhere) and the generalization map   given by $A_2\buildrel\Delta\over\lra A_1\otimes A_1$. 

The sheaf $\Ec^{(1,1)}(A)$ is the direct image of the local system associated to $A_1\otimes A_1$ from
$\CC\setminus \RR_{\geq 0}$ to the whole of $\CC$. So its stalk at $0$ is $A\otimes A$
while the stalk at any point of $\RR_{>0}$ is the sum of two copies of $A\otimes A$ corresponding to
the two sides of the cut $\RR_{\geq 0}$. The generalization map is given by $(\Id, R)$.
Thus the diagram of stalks and generalization maps given by $\ul\delta$ over $\RR_{\geq 0}$
has the form ($\ul\delta$ acts vertically):
\[
\xymatrix{
A_1\otimes A_1 \ar[r]^{\hskip -0.5cm (\Id,R)}
\ar[d]_{\mu_{1,1}} & (A_1\otimes A_1)^{\oplus 2}
\ar[d]^{+}
\\
A_2 \ar[r]^{\Delta_{1,1}} & A_1\otimes A_1
}
\]
 The commutativity of this diagram (i,e., the fact that $\ul\delta$ is indeed a morphism of sheaves)
 is the identity $\Delta_{1,1}\circ\mu_{1,1}  = \Id+R$ which is a particular case of \eqref{eq:br-bialg-grad}, see also
   \eqref{uv:n=2} below.
\end{ex}
  
  
\paragraph{Vertical factorizability of Cousin complexes.}

\begin{Defi}\label{def:prec}
Let $U_1, U_2\subset \CC$ be open sets. We say that $U_1, U_2$ are {\em vertically disjoint} and write
$U_1\prec U_2$, if
$\Im(z_1) < \Im(z_2)$ for any $z_1\in U_1$, $z_2\in U_2$. 
\end{Defi}  

If $U_1, U_2$ are vertically disjoint, then $U_1\cap U_2=\emptyset$.
As in \eqref{eq:add-map}, we have the addition map
\[
a = (a_{pq}) : \bigsqcup_{p+q=n} \Sym^p(U_1)\times\Sym^q(U_2)\buildrel\simeq \over
 \lra \Sym^n(U_1\cup U_2). 
\]

\begin{prop}\label{prop:cous-vert-fact}
The collection of complexes $\Ec(A)=(\Ec_n^\bullet(A))_{n\geq 0}$ is {\em vertically factorizable}. That is,
for any vertically disjoint $U_1, U_2$ we have an isomorphism of complexes of sheaves
 \[
\mu_{U_1, U_2}^{pq}:   \Ec^\bullet_p(A)|_{\Sym^p(U_1)} \boxtimes\Ec^\bullet_q(A)|_{\Sym^q(U_2)}  \lra  a_{pq}^*
  \bigl(  \Ec_n^\bullet(A)|_{\Sym^n(U_1\cup U_2)}\bigr), 
    \]
    these isomorphisms being compatible with restrictions and associative for any triple $U_1 \prec  U_2\prec U_3$
    of vertically disjoint open sets. 
\end{prop}

\noindent{\sl Proof:} Let $p+q=n$. For any two ordered partitions
\[
\beta=(\beta_1,\cdots, \beta_k)\in\2^{p-1}, \quad \gamma = (\gamma_1,\cdots, \gamma_m)\in \2^{q-1}
\]
 we denote their {\em concatenation} by
 \[
 \beta\sqcup\gamma =(\beta_1,\cdots, \beta_k, \gamma_1,\cdots, \gamma_m)\in \2^{n-1}. 
 \]
 Our statement follows from the identifications of sheaves
 \[
\mu_{U_1, U_2}^{\beta, \gamma}:  \Ec^\beta(A)|_{\Sym^p(U_1)} \boxtimes  \Ec^\gamma(A)|_{\Sym^q(U_2 )}
\buildrel\simeq\over\lra a_{pq}^* \left( \Ec^{\beta\sqcup\gamma}(A)|_{\Sym^n(U_1)\cup U_2)}\right).
 \]
 To construct these, we recall that $\Ec^\beta(A)$ is defined as the direct image of the sheaf $\wt\Ec^\beta(A)$
 from the imaginary stratum $X_\beta^\Ien$, and similarly for $\Ec^\gamma(A)$. 
 More precisely, the construction of $\mu_{U_1, U_2}^{\beta, \gamma}$ follows  from:
 \begin{itemize} 
 \item [(1)] The equality of subsets
 \[
 \bigl(X_\beta^\Ien\cap \Sym^p(U_1)\bigr) \, \times \, \bigr(X_\gamma^\Ien\cap \Sym^q(U_2)\bigr) 
 \,\, = \,\, 
 a_{pq}^{-1}\biggl( X_{\beta\sqcup\gamma}^\Ien \cap \Sym^n(U_1\cup U_2)\biggr),
 \]
 and  
  \item[(2)] The identification of two sheaves on the space in (1):
 \begin{itemize}
\item[(2L)] The restriiction of  $\wt\Ec^\beta(A)\boxtimes\wt\Ec^\gamma(A)$ to the LHS in (1), and 

\item[(2R)] The pullback, under $a_{pq}$, of th restriction of $\wt\Ec^{\beta\sqcup\gamma}(A)$
to the RHS in (1). 
\end{itemize} 
 \end{itemize} 
 
 The equality (1)  follows from the very definition of the imaginary strata, while the identification (2)
 follows from the very definition of the sheaves $\wt\Ec^\alpha$. 
  \qed
 

 \paragraph{Stalks of Cousin complexes.}\label{par:En-stalks}
 Recall the notations $B_n(A)$ for the $n$th bar-complexe of $A$.
 For any point $\bz=\sum_{i=1}^m \lambda_i z_i\in\Sym^n(\CC)$ we denote
 \[
 B_\bz(A) \,=\bigotimes\nolimits_{(z_1,\cdots, z_m)} B_{\lambda_i}(A).
  \]

 \begin{prop}\label{prop:En-stalks}
 (a) The complex $\Ec_n^\bullet(A)$ is constructible with respect to the diagonal stratification $S$ of
 $\Sym^n(\CC)$ by complex strata $X_\lambda^\CC$. 
 
 \vskip .2cm
 
 (b) In fact, 
 the stalk of $\Ec^\bullet_n(A)$ at a point $\bz\in\Sym^n(\CC)$  is quasi-isomorphic to $B_\bz(A)$. 
 \end{prop}

 The proof of  Proposition \ref {prop:En-stalks} will occupy the rest of this \S \ref{par:En-stalks}
 Before giving the proof in general, we consider an example which illustrates the situation.
 
 \begin{ex}\label{ex:E-permu}
 Proposition \ref {prop:En-stalks} means, in particular, that the restriction of $\Ec_n^\bullet(A)$ to
 the generic stratum $\Sym^n_\neq (\CC) = X_{(1^n)}^\CC$ is quasi-isomorphic to $\Lc_{1^n}(A)[n]$, the local system of
 determinations of $A_1^{\otimes n}$, put in degree $(-n)$.  However, the stalk $\Ec_n^\bullet(A)_\bz$, $\bz\in \Sym^n_\neq(\CC)$,
 as a complex, depends on the Fox-Neuwirth cell $X_{[\beta:1^n]}$ in which $\bz$ lies.  Let us see why, up to quasi-isomorphism,
 the answer is always the same, namely  the stalk of $\Lc_{1^n}(A)[n]$. 
 
 \vskip .2cm
 
 The simplest case is when $\beta=(1^n)$, i.e., $\Im(z_1) <\cdots < \Im(z_n)$ after a reordering. In this case 
 only the sheaf $\Ec^{(1^n)}(A)$ has non-trivial stalk at $\bz$ (the other sheaves being supported on higher codimension strata),
 and this stalk is $A_1^{\otimes n}$. So in this case $\Ec^\bullet(A)_\bz = A_\bz[n] = B_\bz(A)$ as a complex.
 
 \vskip .2cm
 
 The most complicated case is when $\beta=(n)$, i.e., $\Im(z_1) = \cdots = \Im(z_n)$, and we can assume $\Re(z_1) <\cdots < \Re(z_n)$
 after a reordering. In this case all the  $\Ec^\alpha(A)$ will have nontrivial stalks at $\bz$. The stalk of a given  $\Ec^\alpha(A)$ 
 is, by Proposition \ref{prop:Ea-stalks},
  the sum of as many copies of $A_1^{\otimes n}$, as there are contingency matrices $\zeta$ of type $(\alpha, (1^n))$. 
  Such matrices are, by Example, 
 \ref{ex:1n-stoch}, in bijection with shuffles of type $\alpha$, i.e., with sequences $(I_1, \cdots, I_l)$ of subsets of $\{1,\cdots, n\}$
 forming a disjoint decomposition and such that $|I_\nu|=\alpha_\nu$. Such a sequence gives a face $[I_1,\cdots, I_l]$
 of the $n$th permutohedron $P_n$, see Proposition \ref{prop:Pn-faces}.  When $\alpha$ and $\zeta$ vary, 
 we get precisely all the faces.  Thus the stalk at $\bz$, as a complex, is:
 \[
 \Ec^\bullet(A)_\bz \,=\, A_1^{\otimes n}[n] \otimes_\k  C^\bullet(P_n,\k),
 \]
where $C^\bullet(P_n,\k)$ is the cellular cochain complex of $P_n$ (with respect to the CW-decomposition given by the faces)
with coefficients in $\k$.
Since $P_n$, being a convex polytope, is contractible, $C^\bullet(P_n,\k)$ is quasi-isomorphic to $\k$ in degree $0$,
which gives our desired quasi-isomorphism for $\Ec^\bullet(A)_\bz$. 
 \end{ex}
 
 \vskip .2cm

 \noindent \ul{\sl Proof of Proposition  \ref {prop:En-stalks}:} We will prove (a) and (b) simultaneously. 
To prove  (a) we need to show that for any $\lambda\in\Pc(n)$,  any $\bz\in X_\lambda^\CC$
 and any nearby $\bz'$ lying in the same $X_\lambda^\CC$, the generalization map
 $\gamma_{\bz,\bz'}: \Ec^\bullet_n(A)_\bz \to  \Ec^\bullet_n(A)_{\bz'}$ is a quasi-isomorphism. 
 So let $\bz$ be fixed. 
 
 \vskip  .2cm

 Let $\beta\leq\gamma$ be such that $\bz\in X_{[\beta:\gamma]}$.
 We notice first that because of vertical factorizability of $\Ec^\bullet(A)$,  both our statements
 (quasi-isomorphicity of $\gamma_{\b\bz'}$ for nearby $\bz'\in X_\lambda^\CC$ as well as the identification
 of $\Ec^\bullet_n(A)_\bz$)
 reduce to the case where $\beta=(n)$, i.e., $\Im(z_1)=\cdots=\Im(z_m)$. Ordering the $z_i$ 
 according to the increasing real part, we see that $z_i$ enters into $\bz$ with multiplicity $\gamma_i$,
 i.e., $\lambda =\ol\gamma$ is the unordered partition of $n$ associated to then ordered partition $\gamma$. 
 This means that
 \[
 B_\bz(A) \,=\, B_\gamma(A) \,=\, B_{\gamma_1}(A) \otimes \cdots \otimes B_{\gamma_m}(A)
 \]
 (isomorphism of complexes of objects of $\Vc$). Further, let us write each bar-complex $B_{\gamma_i}(A)$
 as a sum of explicit summands, i.e., identify, as a graded object:
 \[
 B_{\gamma_i}(A) \,= \bigoplus_{\theta^{(i)} = (\theta^{(i)}_1, \cdots, \theta^{(i)}_{m_i})\in\2^{\gamma_i}}
  A_{\theta^{(i)}}[m_i], \quad A_{\theta^{(i)}}= A_{\theta^{(i)}_1} \otimes\cdots \otimes A_{\theta^{(i)}_{m_i}},
  \quad m_i = l(\theta^{(i)}). 
 \]
  Thus, as a graded object,
  \be\label{eq:Bz-A}
  B_\bz(A) \,\,= \bigoplus_{(\theta^{(1)}, \cdots, \theta^{(m)})\,\in\, \prod_{j=1}^m \2^{\gamma_j-1}}
 A_{\theta^{(1)}}\otimes\cdots\otimes A_{\theta^{(m)}}\left[\sum l(\theta^{(i)})\right]. 
  \ee
  We compare this with the shape of the stalk  $\Ec^\bullet_n(A)_\bz$. By Proposition
  \ref {prop:Ea-stalks}, it is a certain direct sum: first, over $\alpha\in\2^{n-1}$ (parametrizing the individual
  $\Ec^\alpha(A)$)
  and, second, over contingency
  matrices $\zeta\in S(\alpha,\gamma)$ (parametrizing the summands in the stalk of  $\Ec^\alpha(A)$
 at $\bz$). 
  For such $\zeta=\|\zeta_{ij}\|$ we denote by $\zeta_{\bullet j}$ the
  $j$th column of $\zeta$. This is a vector in $\ZZ_{\geq 0}^{l(\alpha)}$ with entries summing up to $\gamma_j$,
  with some of these entries possibly being zero. We denote by $\c_j(\zeta)$ and call the 
  {\em reduced $j$th column} of $\zeta$ the ordered partition of $\gamma_j$ obtained by removing the
  zeroes from the entries of $\zeta_{\bullet j}$. The summand corresponding to $\alpha$ and $\zeta$
  is, by Proposition
  \ref {prop:Ea-stalks}, the tensor product $A_\zeta$ which we can canonically identify as
  \[
  A_\zeta \,=\, A_{\c_1(\zeta)} \otimes\cdots \otimes A_{\c_m(\zeta)}. 
  \]
So as a graded object,
\be\label{eq:EnA-z}
\Ec_n^\bullet(A)_\bz \,=\,\bigoplus_{\alpha\in\2^{n-1}} \bigoplus_{\zeta\in S(\alpha,\gamma)} 
A_{\c_1(\zeta)} \otimes\cdots \otimes A_{\c_m(\zeta)} [l(\alpha)].
\ee
This means that each summand $A_{\theta^{(1)}}\otimes\cdots\otimes A_{\theta^{(m)}}$ in
\eqref{eq:Bz-A}  will be found in \eqref{eq:EnA-z} possibly several times: once for each
occurrence of $\alpha\in\2^{n-1}$ and $\zeta\in S(\alpha,\gamma)$ such that $\c_j(\zeta) = \theta^{(j)}$ for
all $j=1,\cdots, m$. 

The possible $\zeta$ with this property are in bijection with XO-tables of type
$(l(\theta^{(1)}),\cdots, l(\theta^{(m)})$, see Definition \ref{def:XO-tab}. Indeed, all we need to do to recover $\zeta$ is,
for any $j$,  to mark by X the positions in the $j$th column which will be occupied by nonzero entries of $\theta^{(j)}$. 
There are $l(\theta^{(j)})$ such positions. 
Once we know $\zeta$, the partition $\alpha$ is found uniquely by summing the rows. This implies that at the level of graded
objects, we have an identification
\be\label{eq:En-XO}
\Ec_n^\bullet(A)_\bz \,=\, 
 \bigoplus_{(\theta^{(1)}, \cdots, \theta^{(m)})\,\in\, \prod_{j=1}^m \2^{\gamma_j-1}}
 \left( A_{\theta^{(1)}}\otimes\cdots\otimes A_{\theta^{(m)}}\left[\sum l(\theta^{(i)})\right]
 \right) \otimes_\k C^\bullet(\XO_{ l(\theta^{(1)}),\cdots, l(\theta^{(m)})}),
\ee
where $C^\bullet(\XO_{  l(\theta^{(1)}),\cdots, l(\theta^{(m)})})$ is the XO-cochain complex, see 
Definition \ref{def:xo-cochain},
and $\otimes_\k$ means the tensor product of a (graded)  object of $\Vc$ and  a  (graded) $\k$-vector space  
 given by the $\k$-linear structure on $\Vc$. 
 
 \vskip .2cm
 
 We next define a morphism of graded objects
 \[
 \gamma: \Ec_n^\bullet(A)_\bz\lra B_\bz(A)
 \]
 in terms of decompositions \eqref {eq:En-XO} and \eqref {eq:Bz-A}. For every 
 $(\theta^{(1)}, \cdots, \theta^{(m)})\,\in\, \prod_{j=1}^m \2^{\gamma_j-1}$
 we consider the morphism of complexes
 \be\label{eq:q-theta}
 \begin{gathered}
 q_{\theta^{(1)}, \cdots, \theta^{(m)}}\,=\, \Sigma\circ p, 
 \\
 C^\bullet(\XO_{  l(\theta^{(1)}),\cdots, l(\theta^{(m)})})\buildrel p \over\lra
 C^0(\XO_{  l(\theta^{(1)}),\cdots, l(\theta^{(m)})}) \buildrel\Sigma\over\lra \k, 
 \end{gathered} 
 \ee
 where $p$ is the projection to the degree $0$ component and $\Sigma$ is given by summing the values of a
 $0$-cochain over all the vertices of the XO-complex. We then define $\gamma$ to act om the summand in 
 \eqref{eq:En-XO} corresponding to $(\theta^{(1)}, \cdots, \theta^{(m)})$ as 
 $\Id \otimes_\k  q_{\theta^{(1)}, \cdots, \theta^{(m)}}$. 
 
 \vskip .2cm
 
 We claim that $\gamma$ is in fact a morphism of complexes. This can either be seen directly by analyzing the
 differential in $\Ec_n^\bullet(A)_\bz$ (described below) or, more conveniently for us, deduced from the following
 observation:
 \vskip .2cm
 {\em 
 $\gamma$ is a particular instance of the generalization maps $\gamma_{\bz,\bz'}$ where
 $\bz'$ is a small displacement of $\bz$ in the same stratum $X_\lambda^\CC$.
 }
 \vskip .2cm
 
 More precisely, $\gamma$ corresponds to
   the case when $\bz'$ is a {\em maximally generic displacement} of $\bz$ inside $X_\lambda^\CC$,
 that is $\bz'\in X_{[\gamma:\gamma]}$, while $\bz\in X_{[(n):\gamma]}$. In other words,  
 assume that the
  displacement is of the form
 \[
 \bz = \sum_{i=1}^m  \gamma_i z_i \quad \rightsquigarrow \quad \bz' = \sum_{i=1}^m  \gamma_i z'_i,
 \quad \Im(z'_1) < \cdots < \Im (z'_m) 
 \]
 (we move each $z_i$ in the vertical direction a bit to make the imaginary parts all different). 
 In this case, similarly to Example \ref{ex:E-permu},
 \[
 \Ec_n^\bullet(A)_{\bz'} \,=\, B_{\bz'}(A)
 \]
 as a complex. If we identify $B_{\bz'}(A)$ with $B_{\bz}(A)$ by continuity of the braided tensor prouct 
 $\bigotimes^{(m)}_{(z_1,\cdots, z_m)}$
 under displacement of the points, then the map $\gamma_{\bz,\bz'}$ becomes
 precisely $\gamma$.  Since we known that each $\gamma_{\bz, \bz'}$ is a morphism of complexes,
 so is $\gamma$.

\vskip .2cm
 
 We now analyze the differential in 
 $\Ec^\bullet_n(A)_\bz$. Fix $\alpha\in\2^{n-1}$ and $\zeta\in S(\alpha,\gamma)$, so
 that the corresponding summand in  $\Ec^\bullet_n(A)_\bz$ is $A_\zeta$.  Let $\zeta_{\nu,\bullet}$ be the $\nu$th row of
 $\zeta$.  (After omitting the zeroes,   $\zeta_{\nu,\bullet}$ becomes  an ordered partition of $\alpha_\nu$).
 As before, let also $\theta^{(j)}= \c_j(\zeta)$ be the ordered partition of $\gamma_j$ obtained by
 omitting zeroes of the column $\zeta_{\bullet,j}$.

 Further let $\beta\in\2^{n-1}$ be such that
 $\alpha>_1\beta$ so $\beta$ is obtained by combining some $\alpha_i$ and $\alpha_{i+1}$ together, as in 
 \eqref{eq:alpha-and-beta}. 
 Let us look at the action of $\delta_{\alpha,\beta}: \Ec^\alpha(A)_\bz\to\Ec^\beta(A)_\bz$ on
 the summand $A_\zeta$. This action consists of a morphism $\delta_{\zeta, \zeta'}: A_\zeta\to A_{\zeta'}$,
 where $\zeta'\in S(\beta,\gamma)$ is obtained by adding together the rows $\zeta_{i,\bullet}$ and $\zeta_{i+1,\bullet}$
 and keeping the other rows of $\zeta$ intact. Now, there are the following possibilities:
 
 \begin{itemize}
 \item[($0$)] The rows  $\zeta_{i,\bullet}$ and $\zeta_{i+1,\bullet}$ have disjoint sets of X's (i.e., of positions filled
 with non-zero entries).  
 \item[($1$)]   $\zeta_{i,\bullet}$ and $\zeta_{i+1,\bullet}$ have  exactly one X-position in common.

 \item[($2^+$)] $\zeta_{i,\bullet}$ and $\zeta_{i+1,\bullet}$ have  some number $k\geq 2$
  X-positions in common.
  \end{itemize}
  
  In the case ($0$) we have $A_{\zeta'}=A_\zeta$ and $\delta_{\zeta,\zeta'}=\pm\Id$. This is 
a matrix element of the differential $d_\XO$ in 
$C^\bullet(\XO_{ l(\theta^{(1)}, \cdots, l(\theta^{(m)}})$ multiplied by the identity of the first factor in the 
corresponding summand in \eqref{eq:En-XO}. 

\vskip .2cm

  In the cases ($1$) and ($2^+$) the differential $\delta_{\zeta,\zeta'}$ is given by the multiplication in $A$.
  But the difference between them is as follows. 
   In the case ($1$) we have a matrix element of the differential $d_B$ in 
 $B_\bz(A) = \bigotimes B_{\gamma_\nu}(A)$ (tensored by the identity of $C^\bullet(\XO_{ l(\theta^{(1)},
 \cdots, l(\theta^{(m)})})$).
That last differential  is the tensor products of the differentials in the $B_{\gamma_\nu}(A)$ and so  each
   matrix element of it acts by increasing the degree by $1$ in just a single tensor factor $B_{\gamma_\nu}(A)$,
   leaving the other tensor factors intact.
 In contrast, in the case ($2^+$) the matrix element, while having a similar general nature (multiplying some {\em groups of
 components} of $A$ together)
 acts simultaneously on more than one tensor factor constituting $B_\bz(A)$. 
 
 \vskip .2cm
 
  To account for this difference, we note that  the decomposition \eqref{eq:En-XO} gives a bigrading of 
  $\Ec^\bullet_n(A)_\bz$,
  and $\delta_{\zeta,\zeta'}$ has, with respect to this bigrading, the  bidegree $(1-k,k)$  in all cases
  ($k=0,1, 2,\cdots$). 
 So considering the filtration of  \eqref{eq:En-XO} whose $d$th layer is the sum of the
 \[
 \left( A_{\theta^{(1)}}\otimes\cdots\otimes A_{\theta^{(m)}}\left[\sum l(\theta^{(i)})\right]
 \right)^{\geq d} \otimes_\k C^\bullet(\XO_{ l(\theta^{(1)}),\cdots, l(\theta^{(m)})}),
 \]
 we see that the quotients of the filtration are the complexes $C^\bullet(\XO_{l(\theta^{(1)}, \cdots, l(\theta^{(m)})})$
 tensored with various homogeneous parts of $A_{\theta^{(1)}}\otimes\cdots\otimes A_{\theta^{(m)}}$
 (the latter considered as graded objects with no differential). 
  Since $C^\bullet(\XO_{ l(\theta^{(1)}), \cdots, l(\theta^{(m)}) })$ is, by Proposition \ref{cor:CXO-contr},
 quasi-isomorphic to $\k$ in degree $0$, the morphism $q_{\theta^{(1)}, \cdots, \theta^{(m)}}$
 in \eqref{eq:q-theta} is a quasi-isomorphism (here we use the assumptions $\on{char}(\k)=0$). 
 Therefore 
  the associated graded object, i.e., the term $E_1$
 of the spectral sequence of our filtered complex, is identified with $B_\bz(A)$,
 and the differential $d_1$ is the differential in $B_\bz(A)$ (coming from Case ($1$)). 
 We note further that $\gamma$ is a morphism of filtered complexes if we consider the
 ``stupid filtration''  of the target by the $B_\bz(A)^{\geq j}$ and so gives a morphism of the
 spectral sequences associated to the source and target.  By what we shown, this  morphism
 of spectral sequences
 gives an isomorphism of the $E_1$-terms so it is an isomorphism, and therefore $\gamma$
 is a quasi-isomorphism. In particular,  
  $\Ec^\bullet_n(A)_\bz$  is quasi-isomorphic to $B_\bz(A)$. 
 This proves part (a) of   Proposition  \ref {prop:En-stalks}.
 
 \vskip .2cm
 
 Let us prove part (b). This part means that each generalization map $\gamma_{\bz, \bz'}$ where $\bz, \bz'$
 lie in the same complex stratum $X_\lambda$, is a quasi-isomorphism. Now, by the transitive property of the
 generalization maps, it is enough to establish quasi-isomorphicity for the case when $\bz'$ is the maximally
 generic displacement within $X_\lambda$, considered above. In this case, as we have seen, $\gamma_{\bz, \bz'}$
 is identified with the map $\gamma$ which we have just shown  to be a quasi-isomorphism. 
  Proposition  \ref {prop:En-stalks} is proved.

  
  \paragraph {Costalks of Cousin complexes.} 
  We now identify the {\em costalks} of $\Ec_n^\bullet(A)$, i.e., the complexes of vector spaces $i_\bz^! \Ec_n^\bullet(A)$,
  $\bz=\sum_{i=1}^m \lambda_i z_i\in\Sym^n(\CC)$.   Here $i_\bz: \{\bz\}\to\Sym^n(\CC)$ is the embedding of $\bz$. 
  
  Recall that  $B_n^\dagger(A)$ is the $n$th cobar-complex of $A$.
  As before, we write 
  \[
 B^\dagger_\bz(A) \,=\bigotimes\nolimits_{(z_1,\cdots, z_m)} B^\dagger_{\lambda_i}(A). 
\]

  \begin{prop}\label{prop:En-costalk}
   The costalk of $\Ec^\bullet_n(A)$ at a point $\bz\in\Sym^n(\CC)$ is quasi-isomorphic to $B^\dagger_\bz(A)$. 
  \end{prop}
  
  \noindent{\sl Proof:} As with Proposition  \ref {prop:En-stalks},  we assume that $\bz\in X_{[\beta:\gamma]}$. 
  Vertical factorizability of the $\Ec_n^\bullet(A)$ reduces our statement to the case $\beta=(n)$. That is,
  $\bz\in X_{(n)}^\Ien$ lies in the minimal imaginary stratum.  Recall that
 $X_{(n)}^\Ien$ consists of $\bz$ such that   $\Im(z_1)=\cdots=\Im(z_m)$, or, put differently,
  $\bz$ is obtained by a translation, by a purely imaginary number,  from a divisor in $K_n=\Sym^n(\RR)$. 
  This means that we have
  an identification $X_{(n)}^\Ien \simeq K_n\times i\RR$.    Let 
    $ \rho: X_{(n)}^\Ien \lra K_n$
 be the projection.

   We look at the costalks at $\bz$ of individual
  terms   $\Ec^\alpha(A)$ of the complex $\Ec_n^\bullet(A)$. Recall that we denote by $j_\alpha: X_\alpha^\Ien\to \Sym^n(\CC)$
  the embedding of the imaginary stratum corresponding to $\alpha$.

   If $\alpha\neq (n)$, then $\bz$ lies in the closure
  of $X_\alpha^\Ien$ but not in $X_\alpha^\Ien$ itself. This implies that the costalk at $\bz$ of 
  $\Ec^\alpha(A)=j_{\alpha_*}\wt\Ec^\alpha(A)$ is $0$ by the base change theorem. 
  
  If $\alpha=(n)$, then $\Ec^{(n)}(A)$  is already supported on the closed stratum $X_{(n)}^\Ien$.
  Denoting, by a minor abuse of notation, the embedding $\{\bz\}\to X_{(n)}^\Ien$ by the same symbol $i_\bz$,
  we have 
  \[
  i_\bz^! \Ec^{(n)}(A) \,=\, i_\bz^! \bigl(j_{(n)}^!\Ec^{(n)}(A) \bigr)
  \]

  The   sheaf  $j_{(n)}^!\Ec^{(n)}(A)$ on $X_{(n)}^\Ien$ coincides wth $j_{(n)}^*\Ec^{(n)}(A)$ and
  is what we denoted earlier
  \eqref{eq:wt-Ea}  by $\wt\Ec^{(n)}(A)$. So
  it has the form
  $\rho^* R_n(A)$, where $R_n(A)$ is  the constructible sheaf on $K_n$ corresponding
   to $\BB_n^\dagger(A)$, the $n$th
  cobar-cube of $A$. 
  
  Now, $K_n$ with its stratification $S_\RR$ and $X_{(n)}^\Ien$ with the pullback stratification are quasi-regular
  cell complexes (see Definition \ref {def:quasi-reg}), and $ R_n(A)$ and, therefore, $\Ec_n^{(n)}(A)$ are cellular sheaves on them.
  We now recall the following well known fact.

  \begin{lem}
  Let $Y$ be a quasi-regular cell complex and $\Fc$ be a cellular sheaf on $Y$.
  For a cell $\sigma$ of $Y$ let $j_\sigma:\sigma\to Y$ be the embedding. 
  Denote the stalk of $\Fc$
  at  $\sigma$ by $F_\sigma = \Gamma(\sigma, j_\sigma^*\Fc)$.
  Then the costalk $ \Gamma(\sigma, j_\sigma^!\Fc)$  at $\sigma$ is quasi-isomorphic to the complex
  \[
  F_\sigma \buildrel\gamma\over\lra \bigoplus_{\ol\tau\supset\sigma
  \atop \dim(\tau=\dim(\sigma)+1} F_\tau\otimes\OR(\tau/\sigma)  \buildrel\gamma\over\lra
   \bigoplus_{\ol\tau\supset\sigma
  \atop \dim(\tau=\dim(\sigma)+2} F_\tau\otimes\OR(\tau/\sigma) \buildrel\gamma\over\lra\cdots
  \]
  where $\OR(\tau/\sigma)$ is the 1-dimensional relative orientation $\k$-vector space,
  and the differential $\gamma$ is induced by the generalization maps of $\Fc$.
  The grading of this complex is normalized so that $F_\sigma$ is in degree $0$. 
   \qed
  
  \end{lem}  
  
 We apply this lemma to $\Fc= \wt\Ec^{(n)}(A)$. The cells in $X_{(n)}^\Ien$  are the fine cells $X_{[(n):\gamma]}$. 
 If $\bz \in X_{[(n):\gamma]}$ as we assumed, then the complex of the lemma will 
 be precisely $B_\bz(A)$ but with  a shift: the grading  will start from degree $0$ on the left. 
 This is because it describes the stalk {\em at the cell} (i.e., $j_\sigma^!$ instead of $j_\bz^!$).
  Passing from that to the stalk at $\bz$ itself
 removes the shift (by introducing the cohomology of a constant sheaf on the cell with support at $\bz$) 
 and we get $B_\bz(A)$ with correct grading. \qed

 
 \paragraph{Perversity of Cousin complexes.}
 
 \begin{prop}
 The complex $\Ec_n^\bullet(A)$ is a perverse sheaf constructible with respect to the stratification $S$,
 i.e., an object of $\Perv(\Sym^n(\CC), S,\Vc)$. 
 \end{prop}
 
 \noindent{\sl Proof:} The $S$-constructibility is shown in  Proposition   \ref {prop:En-stalks}. 
So it remains to prove the perversity conditions  $(P^+)$ and $(P^-)$, see \S \ref{subsec:skel-pure}\ref{par:skel-pure}
 
 \vskip .2cm
 
 We start with $(P^+)$. As the stalk  $i_\bz^*\Ec_n^\bullet(A)$ is quasi-isomorphic to $B_\bz(A)$ (Proposition \ref{prop:En-stalks}),
  we need to show that the locus of $\bz$
 for which $H^{-j} (B_\bz(A))\neq 0$, has complex dimension $\leq j$. Indeed, suppose $\bz\in X_\lambda^\CC$,
 so $\bz = \sum_{i=1}^m \lambda_i z_i$ with $z_i$ distinct. By definition, $B_\bz (A)$ is the tensor
 product of the $B_{\lambda_i}(A)$, and each $B_{\lambda_i}(A)$ is a complex concentrated in degrees fron
 $(-\lambda_i)$ to $(-1)$. So $B_\bz(A)$ is concentrated in degrees from $(-n)$ to $(-m)$. 
 Now, the complex dimension of $X_\lambda^\CC$ is $m$, and $(P^+)$ follows.
 
  \vskip .2cm
  
 The condition $(P^-)$ follows similarly from Proposition \ref{prop:En-costalk}. \qed
 
 
 \subsection{Factorizability of Cousin complexes:  proof of Theorem \ref{thm:localize}} 
 
 \paragraph{ From vertical factorizability to full factorizability.} Recall (\S \ref{subsec:fact-perv}) that the  space $E_2(m)$ of the $E_2$-operad consists of $m$-tuples $(U_1,\cdots, U_m)$ of disjoint round disks in the unit disk $D\subset\CC$. 
 Let $E_2^\prec(n)\subset E_2(m)$ be the subspace formed by $(U_1,\cdots, U_m)$ such that $U_1\prec\cdots\prec U_m$,
 see  Definition \ref{def:prec}. The space $E_2^\prec(n)$ is contractible and together they form a (non-symmetric)
 suboperad  $E_2^\prec$ in $E_2$. Thus $E_2^\prec$ is a non-symmetric verstion of the homotopy associative operad $E_1$. 
 
\vskip .2cm

Let $A$ be a primitive bialgebra in $\Vc$, as before. We want to equip the collection
of perverse sheaves  $\Ec(A) = (\Ec_n^\bullet(A))$ with a structure of a factorizable perverse sheaf in the sense of
Definition \ref{def:FPS}. That is, 
 for each $m\geq 0$ and each $(U_1, \cdots U_m)\in E_2(m)$,  we want to construct an isomorphism of perverse
 sheaves (i.e, an isomorphism of the objects of the derived category)
 \[
\mu_{U_1, \cdots, U_m}:  \bigboxtimes \nolimits_{(U_1, \cdots, U_m)} \Fc|_{\Sym(U_i)} \lra    a^*\left(  \Fc|_{\Sym(\bigcup U_i)}\right),
 \]
so that these isomorphisms are compatible with the operadic compositions.

\vskip .2cm

For $(U_1,\cdots, U_m)\in E_2^\prec(m)$,  such an isomorphism is  already given by Proposition \ref {prop:cous-vert-fact}.
In this case $\mu_{U_1,\cdots, U_m}$ is in fact an isomorphism of complexes of sheaves. These isomorphisms
are compatible with the operadic compositions in $E_2^\prec$. 

\vskip .2cm

Next, we want to define  the isomorphisms
$\mu_{U'_1,\cdots, U'_m}$ for any $(U'_1,\cdots, U'_m)\in E_2(m)$ by extending them
from $E_2^\prec(m)$ to $E_2(m)$  ``by continuity''.  That is, let 
\[
\gamma: [0,1] \lra E_2(m), \quad t\mapsto \gamma(t) = (U_1(t), \cdots, U_m(t))
\]
be a continuous path in $E_2(m)$. We denote $U_i=U_i(0)$ and $U'_i=U_i(1)$. Let also $U(t)=\bigcup U_i(t)$ and  
$U=U(0)$, $U'=U(1)$. 

Because perverse sheaves are a purely topological concept, the 1-parameter family of disks $U_i(t)$ defines a 
1-parameter family of equivalences of categories (``isomonodromic deformations
of perverse sheaves'')
\[
h_t: \Perv(\Sym(U_i(0), S) \lra \Perv(\Sym(U_i(t)), S).
\]
Moreover, if $\Fc\in\Perv(\Sym(D), S)$, then $h_t$ takes the restriction  $\Fc|_{\Sym(U_0)}$ to $\Fc|_{\Sym(U_t)}$.
We have a similar family of equivalences of categories for the unions
\[
k_t: \Perv(\Sym(U(0), S) \lra \Perv(\Sym(U(t)), S).
\]
Taking $t=1$, we obtain
the isomorphism  (``monodromy map'')
\[
M_\gamma: \Hom\biggl( \bigboxtimes \nolimits_{(U_1, \cdots, U_m)} \Fc|_{\Sym(U_i)} , a^*\left(  \Fc|_{\Sym(U)}\right)
\biggr) \to  \Hom\biggl( \bigboxtimes \nolimits_{(U'_1, \cdots, U'_m)} \Fc|_{\Sym(U'_i)} , a^*\left(  \Fc|_{\Sym(U')}\right)
\biggr).
\]
These isomorphisms are compatible with concatenation of paths. 

\vskip .2cm

 We now fix $(U_1,\cdots, U_m)\in E_2^\prec(m)$,  consider some $(U'_1,\cdots, U'_m)\in E_2(m)$
 and choose a path $\gamma$ joining them as above. We define
 \be\label{eq:mu-continued}
 \mu_{U'_1,\cdots, U'_m} \,=\, M_\gamma(\mu_{U_1,\cdots, U_m}).
 \ee
For this  to be well-defined, we need to establish the following.

\begin{lem}
Suppose $\gamma$ is a closed path in $E_2(m)$ with $\gamma(0)=\gamma(1) = (U_1,\cdots, U_m)\in E_2^\prec(m)$.
Then $M_\gamma(\mu_{U_1,\cdots, U_m}) = \mu_{U_1,\cdots, U_m}$.
\end{lem}

\noindent{\sl Proof:} This follows from the identification, up to quasi-isomorphism, of the stalks of the $\Ec_n^\bullet(A)$
in Proposition \ref{prop:En-stalks}. That is, the stalks at $\bz=\sum_{i=1}^m \lambda_i z_i$, $z_i\neq z_j$, are,
up to quasi-isomorphism (in particular, at the level of their cohomology), manifestly factorizable
whether $z_1 \prec\cdots \prec z_m$ or not. This  means that the two
morphisms of perverse sheaves,  $M_\gamma(\mu_{U_1,\cdots, U_m})$
and $ \mu_{U_1,\cdots, U_m}$ ,  induce the same maps on the cohomology of the stalks of the perverse sheaves in
question. This, in its turn, implies that these two morphisms of perverse sheaves coincide. \qed

Well-definedness being established, we easily conclude that:

\begin{prop}\label{prop:E-fact}
The  isomorphisms \eqref{eq:mu-continued} make $\Ec(A)=(\Ec_n^\bullet(A))$ into a factorizable perverse sheaf
on $\Sym(\CC)$ with values in $\Vc$. \qed
\end{prop}


\paragraph{Proof of Theorem  \ref{thm:localize}.}
Proposition \ref{prop:E-fact} defines a functor
\[
L: \PB(\Vc) \lra \FPS(\Vc), \quad A \mapsto L(A) = (L_n(A))_{n\geq 0}, \,\,\, L_n(A):=\Ec_n^\bullet(A). 
\]
To prove that $L$ is an equivalence, we construct a functor in the opposite direction
\[
\Phi:  \FPS(\Vc)\lra \PB(\Vc), \quad \Fc= (\Fc_n)_{n\geq 0} \mapsto \Phi(\Fc) = \bigoplus_{n\geq 0} \Phi_n(\Fc),\,\,\,
\Phi_n(\Fc) := \Phi_\tot(\Fc_n). 
\] 
We recall  \eqref{eq:total-van} that $\Phi_\tot(\Fc_n)$ is the stalk at $0$ of the sheaf $R_{K_n}(\Fc_n) = \ul\HH^0_{K_n}(\Fc)$,
where $K_n=\Sym^n(\RR)$. As $\Fc$ is factorizable, we see that the stalk of $R_{K_n}(\Fc_n)$ at the stratum
$K_\alpha$, $\alpha=(\alpha_1,\cdots, \alpha_m)\in\2^{n-1}$, is canonically identified with
$\Phi_{\alpha_1}(\Fc) \otimes\cdots\otimes\Phi_{\alpha_m}(\Fc)$. 
Consider the generalization maps from  strata corresponding to $\alpha$ with $l(\alpha)=1$ to 
strata corresponding to $\alpha$ with $l(\alpha)=2$,
for the sheaf structures of the $R_{K_n}(\Fc_n)$. These maps have the form $\Delta_{p,q}:\Phi_{n}(\Fc)\to\Phi_p(\Fc)\otimes\Phi_q(\Fc)$,
$p+q=n$. So they 
give a 
comultiplication $\Delta$ on $\Phi(\Fc)$. This comultiplication is coassociative, as one can see by considering
strata corresponding to $\alpha$ with $l(\alpha)=3$. It also
satisfies the primitivity condition because
$\Fc_0 = \1$.  

\vskip .2cm
 
To  define a multiplication on $\Phi(\Fc)$ we use the Cousin resolution of $\Fc_n$ from Proposition
\ref{prop:cousin}. That is, we denote $j_\alpha: X_\alpha^\Ien\to \Sym^n(\CC)$ the embedding of the imaginary stratum
and put $\wt\Ec^\alpha(\Fc_n)=j_\alpha^!\Fc_n$, which is a single sheaf on $X_\alpha^\Ien$. Then we put
$\Ec^\alpha(\Fc_n)=j_{\alpha*}\wt\Ec^\alpha(\Fc_n)$, which is a single sheaf of $\Sym^n(\CC)$, and these
sheaves, taken  for all $\alpha\in\2^{n-1}$, form a complex quasi-isomorphic to $\Fc_n$. 

Now, fix a point $\bz\in X_\alpha^\Ien$. We can write $\bz = \sum \bz_i$ where $\bz_i$ is a divisor of degree $\alpha_i$
with all the points on the same horizontal line $\Im(z)=\on{const}$. Therefore 
the pair $(X_\alpha^Ien, \Sym^n(\CC))$  near $\bz$ is homeomorphic to the product of
the pairs $(\Sym^{\alpha_i}(\RR), \Sym^{\alpha_i}(\CC))$  near $\bz_i$. Together with factorizability of $\Fc$ this implies
that the stalk of $\wt\Ec^\alpha(\Fc_n)$ at $\bz$ is the tensor product of the stalks of $\ul \HH^0_{K_{\alpha_i}}(\Fc_{\alpha_i})$
at $\bz_i$. This further implies that the stalk at $0$  of the direct image is identified as
\[
\Ec^\alpha(\Fc_n)_0 \,=\, \Phi_{\alpha_1}(\Fc)\otimes\cdots\otimes\Phi_{\alpha_m}(\Fc).
\]
Let now $\delta$ be the differential in the Cousin complex of $\Fc_n$. As a graded sheaf, this complex is the sum of
summands
\[
\bigoplus_{\alpha\in\2^{n-1}} \Ec^\alpha(\Fc_n)[l(\alpha)]
\]
 Consider the action of the  components of $\delta$ on the stalks of these summands at $0$.
 The components acting from the summands with $l(\alpha)=2$ to the summands with $l(\alpha)=1$
 have the form $\mu_{p,q}: \Phi_p(\Fc)\otimes\Phi_q(\Fc)\to\Phi_n(\Fc)$, $p+q=n$. So they define
 a multiplication $\mu$ in $\Phi(\Fc)$. The condition $\delta^2=0$ implies the associativity of $\mu$.
 The compatibility of $\delta$ and $\mu$ follows from the fact that $\delta$ is a morphism of sheaves,
 i.e., commutes with the generalization maps, by the same analysis as in the proof of Proposition
 \ref{prop:delta-is-morph}. 
 
 \vskip .2cm
 
This defines the functor $\Phi$. If $\Fc$ is a factorizable perverse sheaf, then the above reasoning identifies
the complex $\Ec^\bullet(\Fc_n)$ with $\Ec_n^\bullet(\Phi(A))$, and this shows that  $L(\Phi(\Fc))\simeq\Fc$. 
  If $A$ is a primitive bialgebra, then $\Phi(L(A))\simeq A$ by Propositions
\ref{prop:En-stalks} and \ref{prop:En-costalk}.  This shows that $L$ is an equivalence of categories.

\vskip .2cm

Further, Proposition \ref{prop:En-stalks} implies properties (a) and (d) of the functor $L$ claimed in Theorem 
 \ref{thm:localize}. Proposition  \ref{prop:En-costalk} implies properties (b) and (c).
 Finally, to show the property (e) it is enough to notice that the Verdier duality interchanges stalks and costalks
 so the above identification of the multiplication and comultiplication in $A = \Phi(\Fc)$ identifies comultiplication
 and multiplication in $\Phi(\DD \Fc)$, thus showing that $\Phi(\DD\Fc)\simeq A^*$. Theorem 
  \ref{thm:localize} is proved.

  
  \subsection {Perverse sheaves corresponding to shuffle algebras: proof of Theorem \ref{thm:shuffle} }
  Here we prove Theorem  \ref{thm:shuffle}, with the three paragraphs A,B,C below corresponding to
  parts (a), (b), (c) of the theorem. 
  
  \paragraph{ Localization of the tensor algebra $T_!(V)$.} Let $A=T_!(V)$, so $A_n=V^{\otimes n}$. 
  We know that
  \[
  \Ec_n^\bullet(A)|_{\Sym^n_\neq(\CC)} \,=\,\Lc_{(1^n)}(A)\,=\,\Lc_n(V)[n]
  \]
  is the local system of determinations of $V^{\otimes n}$, pu in degree $(-n)$. To show that 
  $\Ec_n^\bullet(A)  \simeq  j_{n!} \Lc_n(V)[n]$
  it is enough to show that the stalk of $\Ec_n^\bullet(A)$ at any point $\bz \in \Sym^n (\CC)\,\setminus \, \Sym^n_\neq(\CC)$
  is an exact complex. By Proposition \ref{prop:En-stalks}, this stalk is the tensor product of the $B_{\lambda_i}(A)$,
  where $\bz=\sum \lambda_i z_i$ with $z_i$ distinct. The condition $\bz \in \Sym^n (\CC)\,\setminus \, \Sym^n_\neq(\CC)$
  means that some $\lambda_i\geq 2$. So our statement follows from the following.
  
  \begin{lem}\label{lem:T!-stalks}
  If $A=T_!(V)$, then for any $m\geq 2$ the complex $B_m(A)$ is exact. 
  \end{lem}
  
  For example, for $m=2$ the complex $B_2(A)=\{A_1\otimes A_1\buildrel\mu\over\to A_2\}$ is identified with
  $\{V\otimes V\buildrel\Id\over\to V\otimes V\}$ and so is exact. 
  
  \vskip .2cm
  
  \noindent{\sl Proof of the lemma:}  For any graded algebra $A$ we can write $B_m(A)$ as the graded object as follows:
  \[
  B_m(A) \,=\bigoplus_{\alpha=(\alpha_1,\cdots, \alpha_l)\in\2^{m-1}} A_{\alpha_1}\otimes\cdots\otimes A_{\alpha_l}[l],
  \]
  with the differential given by the multiplication in $A$. In our case $A_{\alpha_i}=V^{\otimes \alpha_i}$
  and each component of the multiplication is the identity map.  Therefore
  we have an {\em isomorphism of complexes}
  \[
  B_m(A) \,=\bigoplus_{\alpha=(\alpha_1,\cdots, \alpha_l)\in\2^{m-1}} V^{\otimes m}[l] \,\,\simeq \,\, V^{\otimes m} \otimes_\k
  C^\bullet(\2^{m-1},\k), 
    \]
    where 
    \[
      C^\bullet(\2^{m-1},\k)\,=\ \bigoplus_{\alpha=(\alpha_1,\cdots, \alpha_l)\in\2^{m-1}} \k[l], 
      \]
      with the differential being  the alternating sum of elementary contractions of ordered
      partitions.  Now notice that $ C^\bullet(\2^{m-1},\k)$ is (up to shift) nothing but the augmented cellular chain complex
      of the $(n-2)$-dimensional simplex and so is exact. \qed

  \paragraph{Localization of the cotensor algebra $T_*(V)$.} 
  Let $A=T_*(V)$, so $A_n=V^{\otimes n}$ as before but the multiplication is the shuffe product while
  the comultiplication $\Delta_{pq}: A_{p+q}\to A_p\otimes A_q$ is the identity. 
  As before, $\Ec_n^\bullet(A)|_{\Sym^n_\neq(\CC)} = \Lc_n(A)[n]$. To prove that $\Ec_n^\bullet(A)\simeq j_* \Lc_n(A)[n]$,
  it is enough to prove that the {\em  costalks}  $i_\bz^! \Ec_n^\bullet(A)$ are exact for any 
   $\bz \in \Sym^n (\CC)\,\setminus \, \Sym^n_\neq(\CC)$. This follows from Proposition
   \ref {prop:En-costalk} and from
   the statement dual to 
   Lemma \ref{lem:T!-stalks}:
   
   \begin{lem}
   Let $A=T_*(V)$. Then for any $m \geq 2$ the $m$th cobar-complex $B_m^\dagger(A)$ is exact. 
   \end{lem}
   
   \noindent {\sl Proof:} completely parallel to that of   Lemma \ref{lem:T!-stalks}.\qed

  \paragraph{Localization of the Nichols algebra $T_{!*}(V)$.}
  We notice the following:
  
  \begin{itemize}
  \item[(1)] The canonical map $\alpha: T_!(V)\to T_*(V)$ is the unique morphism of algebras extending
  the identity on $T_!(V)_1 = V = T_*(V)_1$. 
  
  \item[(2)] The canonical map $c_n: j_{n!} \Lc_n(V)\to j_{n*} \Lc_n(V)$ is the unique morphism of perverse sheaves
  which is the identity on $j_n^*  j_{n!} \Lc_n(V) = \Lc_n(V) = j_n^*  j_{n*} \Lc_n(V)$. 
  \end{itemize}
  So Theorem   \ref{thm:localize} on equivalence between factorizable perverse sheaves and primitive
  bialgebras implies that $\alpha$ corresponds, under this equivalence,  to the collection of the $c_n$. 
  In other words, the functor $\Phi_\tot$, see Theorem  \ref{thm:localize}(b),
   takes $j_{n!}\Lc_n(V)$ as well as $j_{n*}\Lc_n(V)$ to $V^{\otimes n}$,
  and takes the map $c_n$ to $\alpha_n$,  the $n$th component of $\alpha$. 
  
  Now, $\Phi_\tot$ is an exact functor from $\Perv(\Sym^n(\CC), S, \Vc)$ to $\Vc$ and so takes images to images.
  Therefore it takes $(j_n)_{!*} \Lc_n(V)$ to $T_{!*}(V)_n$, the $n$th component of the Nichols algebra. 
  This finishes the proof of Theorem  \ref{thm:shuffle}.

  
\section{Examples: Granger-Maisonobe quivers and bialgebras (Cases $n=2,3$)}\label{sec:examples}

\subsection {Generalities: reducing translation invariance}
In this chapter we reivew  the known elementary classification of perverse sheaves on $\Sym^n(\CC)$
for $n=2,3$, smooth with respect to the diagonal stratification $S$.  We also show how a primitive braided
bialgebra $A$ gives, in an elementary way, a perverse sheaf 
(denoted by $L_n(A)$ in Theorem \ref{thm:localize}) 
in terms of these descriptions. 

Let $\Vc$ be a $\k$-linear braided monoidal abelian category with $\otimes$ biexact. 

\vskip .2cm

As before, let  $\Sym^n_0(\CC)\subset \Sym_n(\CC)$ be the subspace of monic polynomials with zero subprincipal part,
see \eqref{eq:sum-divisor} ff.
Let $S^0$ be the stratification induced by $S$ on $\Sym^n_0(\CC)$. 
 Since $S$ is translation invariant, we have an identification
 \be
\Perv(\Sym^n(\CC), S,\Vc) \,\simeq \, \Perv(\Sym^n_0(\CC), S^0, \Vc). 
\ee
Note, however, that this identification involves a shift of degree:
perverse sheaves on $(\Sym^n(\CC), S)$ are situated, as complexes, in degrees $[-n,-1]$, while
perverse sheaves on $\Sym^n_0(\CC), S^0)$ are situated in degrees $[-n+1, 0]$. In
particular, it does
  affect the concept of factorizable
perverse sheaves because of the Koszul sign rule \eqref{eq:koszul} involved in the braiding of shifted objects.

In this chapter we will  work only with $\Sym^n_0(\CC)$, but keep the same notation that we used for various strata 
in $\Sym^n(\CC)$, assuming that they are intersected with $\Sym^n_0(\CC)$.

\subsection {Case $n=2$: the $(\Phi,\Psi)$ description.} 

Note that $\Sym^2_0(\CC)=\CC$, with the complex strata $X_{(2)}^\CC$ and $X_{1,1}^\CC$
 being  $\{0\}$ and $\CC-\{0\}$. Thus $\Perv(\Sym^2_0(\CC), S^0, \Vc)$ is identified with the category
$\Perv(\CC,0, \Vc)$ studied in Example \ref{ex:phi-psi} (for the case of $\Vc=\Vect_\k$, which generalizes
verbatim to  arbitrary $\Vc$). 

\vskip .2cm 

In other words,  $\Perv(\Sym^2_0(\CC), S^0, \Vc)$ is equivalent to the category of diagrams formed by
two objects and two morphisms of $\Vc$
\[
\xymatrix{
E_2 = \Phi  \ar@<.5ex>[r]^v &\Psi = E_{1,1}
\ar@<.5ex>[l]^u
}
\]
 such that  the morphisms $T_\Psi=\Id_\Psi-vu$  and $T_\Phi = \Id_\Phi-uv$  are invertible. 

\vskip .2cm

Let $A$ be a  primitive bialgebra in $\Vc$. We associate to it the diagram
\[
\xymatrix{
\Phi=A_2  \ar@<.5ex>[rr]^{v=\Delta_{1,1}} &&A_1\otimes A_1=\Psi.
\ar@<.5ex>[ll]^{u=\mu_{1,1}}
}
\]

\vskip .2cm
In the calculations below and in \S  \ref {subsec:GM} later we will assume that  we have a fully faithful embedding of $\Vc$, as a monoidal 
 (but not as braided monoidal) category into $\Vect_\k$, in other words, that objects $V$ of $\Vc$ have 
 consist of elements, or vectors  $x\in V$, so that the tensor product of two objects $V\otimes W$
 consists of sums of formal tensors $x_i\otimes y_i$. This assumption is for notational simplification only:
 it is straightforward (but more lengthy)  to restate our arguments  entirely 
 in terms of  diagrams involving various tensor products in $\Vc$.
 
 \vskip .2cm

Then for $a\otimes b\in A_1\otimes A_1$ we get, using the assumption that $A_1$ consists of
primitive elements:
\[
\begin{gathered}
vu(a\otimes b) = \Delta_{1,1}(ab) = \bigl(\Delta(a)\cdot \Delta(b)\bigr)_{1,1} = 
\\
= \bigl( (a\otimes 1 + 1\otimes a) (b\otimes 1 + 1\otimes b)\bigr)_{1,1} = (a\otimes 1)(1\otimes b)  +
 (1\otimes a)(b\otimes 1)=
 \\  = a\otimes b + R_{A_1,A_1}(a\otimes b),
\end{gathered}
\]
the last equality coming from the fact that the multiplication in $A\otimes A$ is defined using $R$.
Therefore
\be\label{uv:n=2}
T_\Psi= \Id-vu  =- R_{A_1, A_1}: A_1\otimes A_1\lra A_1\otimes A_1
\ee
  is the negative of braiding operator and it is invertible by axioms of braided monoidal categories. 
  This minus sign is an instance of the twist of  the braid group action by the sign character in 
  Remark \ref{rem:tw-sign}. 
  
  \vskip .2cm
  
Further,  the invertibility of  
\[
T_\Phi = \Id- \mu_{1,1}\circ\Delta_{1,1} : A_2\lra A_2
\]
follows from that of $T_\Psi$ by \eqref{eq:(1-uv)}. 
 Note that $T_\Phi=\sigma_2$ is   just the antipode
  of $A$ acting on $A_2$.

   
   \subsection  {Case n=3: the Granger-Maisonobe quiver. }\label{subsec:GM}

\paragraph{Imaginary strata for $n=3$.} 
 We have that $\Sym^3_0(\CC)=\CC^2$ is the space of
cubic polynomials  of the form $x^3+ax + b$.  
The open stratum $X_{(1,1,1)}^\CC$
consists of polynomials with distinct roots, i.e., with non-vanishing discriminant.
It is therefore the complement of the complex semi-cubic parabola
\[
Z=  \bigl\{ (a,b)\in\CC^2 | \,\, 4a^3 + 27 b^2 = 0\bigr\}.
\]
 The other two strata $X^\CC_{(2,1)}$ and $X^\CC_{(3)}$  are equal to
 $Z\setminus \{0\}$ (the smooth locus of $Z$) and $\{0\}$.

  The category of perverse sheaves on $\CC^2$ smooth with respect to a
 semi-cubic parabola was described in the papers of Granger-Maisonobe  \cite{GM-cusp} and
 MacPherson-Vilonen \cite{MV} (the latter authors considering the more general case of a curve
 of the form $a^m=b^n$).
 
 \vskip .2cm
 
 We want to emphasize that {\em our method (in the main body of the paper) is a direct generalization of
 the Granger-Maisonobe approach}.  In particular, the cuts (``coupures'') $H$ and $M_\RR$ of 
  \cite{GM-cusp} (whose idea is attributed by the authors to Kyoji Saito) are precisely the unions of  imaginary
  strata $X_{(2,1)}^\Ien \cup X_{(1,2)}^\Ien$ and $X_{(3)}^\Ien = \Sym^3_0(\RR)$ in our sense. 
  So we can say, together with the authors of   \cite{GM-cusp}, that K. Saito
  ``est a l'origine de ce travail''. 
  
  \vskip .2cm
  
  More precisely, for $n=3$ we have $4$ imaginary strata:
  \[
  X_{(3)}^\Ien = \Sym^3_0(\RR), \,\,  X_{(1,2)}^\Ien, \,\, X_{(2,1)}^\Ien, \,\, X_{(1,1,1)}^\Ien
  \]
  whose real dimensions are equal to $2,3,3,4$ respectively.  Thus $X_{(1,1,1)}^\Ien\subset \Sym^3_0(\CC)=\CC^2$
  is open, while $X_{(1,2)}^\Ien$ and $X_{(2,1)}^\Ien$ are {\em locally closed} 3-dimensional subanalytic hypersurfaces
  {\em with boundary}. The boundary of each of them is formed by $Z\setminus \{0\}$ (this boundary
  actually lies in the stratum). Further, each of them contains   {\em in their closure} (but  does not   actually contain)
    the stratum $X_{(3)}^\Ien$. 
    
      \vskip .2cm
      
      For further visualization it is instructive to take, as proposed in   \cite{GM-cusp}, the intersection
      of the entire picture with the unit sphere $S^3 = \{ |a|^2 + |b^2|=1\} \subset \CC^2$. The intersection
      $T= Z\cap S^3$ is, as well known, a trefoil knot in $S^3$. The intersection with $S^3$ of $ X_{(3)}^\Ien $,
      i.e., of the interior of the {\em real} semi-cubic parabola, see Fig. \ref{fig:semicubic},  is the same as the
      intersection of this interior with the unit circle in $\RR^2$. So it is  an interval (arc)  in $S^3$ beginning and ending
      on $T$ and not meeting $T$ at any other points, i.e., a {\em chord} $C$ of $T$.   
      The two points of $C\cap T$, denote them $a$ and $b$, subdivide $T$ into two halves $T_+$ and $T_-$. 
               \vskip .2cm

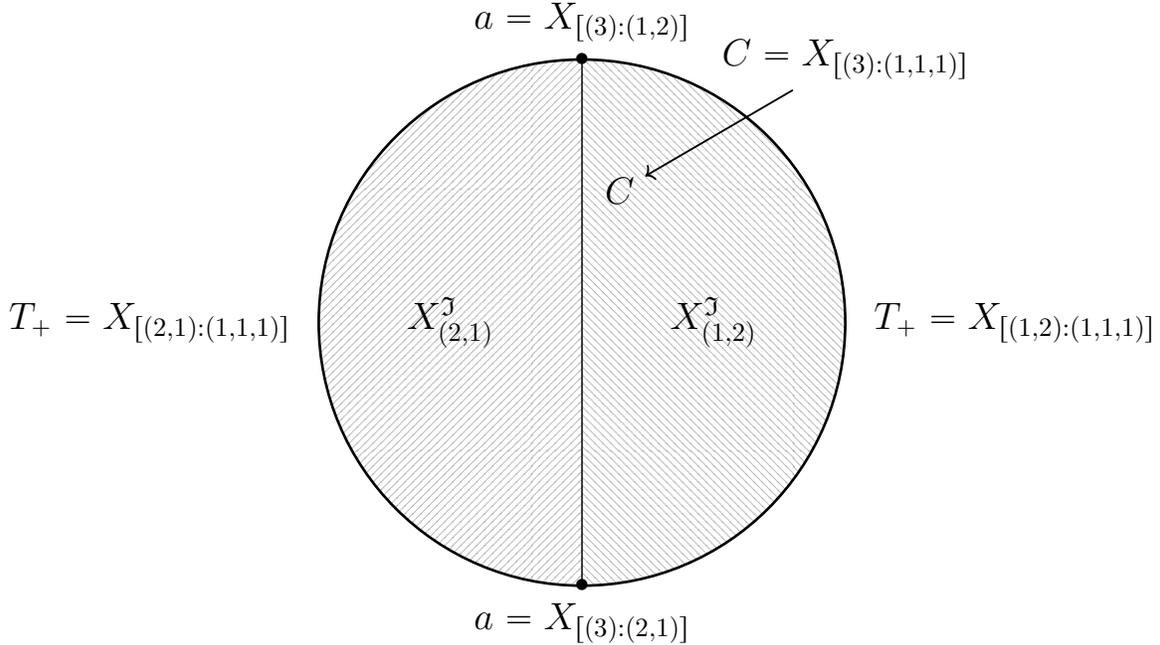
\begin{figure}[H]
\centering
\begin{tikzpicture}[scale=0.5]
\centering

\draw[line width=1]  (0,0) circle (7cm) ; 
\node at (0,7){$\bullet$}; 
\node at (0,-7){$\bullet$}; 
\draw (0,7) -- (0,-7); 
\node at (0,8) {\large$a=X_{[(3):(1,2)]}$}; 
\node at (0,-8) {\large$a=X_{[(3):(2,1)]}$}; 
\node at (-3.5,0) {\large $X_{(2,1)}^\Ien$}; 
\node at (3.5,0) {\large$X_{(1,2)}^\Ien$}; 
\node(C1) at (1,3.5){\large $C$}; 
\node at (-11.5,0){\large$T_+=X_{[(2,1):(1,1,1)]}$};
\node at (11.5,0){\large$T_+=X_{[(1,2):(1,1,1)]}$};
\filldraw[ opacity = 0.5, pattern=north east lines ]  (90:7cm)  arc (90:270:7cm)  (270:7cm) -- (90:7cm); 
\filldraw[ opacity=0.5, pattern=north west lines ]  (270:7cm)  arc (270:450:7cm)  (450:7cm) -- (270:7cm); 
\node (C2) at (7,7) {\large$C=X_{[(3):(1,1,1)]}$}; 

\draw[->, line width =0.7] (C2) -- (C1);

\end{tikzpicture}
\caption{  Imaginary strata and fine cells intersected with $S^3$. Knottedness not depicted. }\label{fig:imstrata}
\end{figure}

    The intersections $X_{(1,2)}^\Ien\cap S^3$
      and $X_{(2,1)}^\Ien\cap S^3$  are $2$-dimensional locally closed surfaces in $S^3$, with boundary.
      They have $T_+$ resp. $T_-$  as the boundary and approach (but do not actually contain)
      the chord $C$.  All parts of this picture are identified with (the intersection with $S^3$ of)
      various fine cells. This is shown symbolically in Fig. \ref{fig:imstrata}, where $T$ is represented as a circle
      but  knottedness is not depicted.

 \paragraph {The Granger-Maisonobe quiver. }
 The description of $\Perv(\Sym^3_0(\CC), S^0, \Vect)$ was given in \cite{GM-cusp} is in terms of diagrams of  $4$ vector spaces
  denoted there by $E, F_1, F_2, G$.
 Our method being  a direct generalization of the method of the Granger-Maisonobe approach,
 their spaces have transparent interpretation in our general terms. 
 We reproduce
 this type of diagram below, keeping their notation on the right and adding our more
 systematic notation on the left:
 \be\label{eq:diag-P3}
 \xymatrix{
 & E_{2,1}=F_1
 \ar@<-.5ex>[dd]_{\alpha}
 \ar@<-.5ex>[dr]_{v_1}
 \ar@<-.5ex>[dl]_{u'_1} &
 \\
 E_{3}=G 
\ar@<-.5ex>[ur]_{v'_1}
 \ar@<-.5ex>[dr]_{v'_2}
 && E_{1,1,1}=E
 \ar@<-.5ex>[dl]_{u_2}
 \ar@<-.5ex>[ul]_{u_1}
 \\&E_{1,2}=F_2. 
 \ar@<-.5ex>[uu]_{\beta}
\ar@<-.5ex>[ur]_{v_2}
 \ar@<-.5ex>[ul]_{u'_2}&
 }
 \ee
 We note that the arguments of  \cite{GM-cusp} generalize verbatim to perverse sheaves with
 values in any abelian category, so we have:

 \begin{thm}[\cite{GM-cusp}]\label{thm:GM}
 Let $\Vc$ be any $\k$-linear abelian category (not necessarily monoidal). 
 $\Perv(\Sym^3_0(\CC), S^0, \Vc)$ is equivalent to the category of diagrams \eqref{eq:diag-P3} 
 formed by objects and morphisms of $\Vc$ and  satisfying the following relations:
 \begin{itemize}
 \item[(1)] $M_1= v_1u_1+\Id$ and $M_2= v_2u_2+\Id$ are isomorphisms.
 
 \item [(2)] $u'_1u_1+u'_2u_2=0$ and $v_1v'_2+v_2v'_2=0$.
 
 \item[(3)] $v'_2 u'_1 = -\alpha + u_2 v_1$ and $v'_1 u'_2 = -\beta + u_1 v_2$. 
 
 \item[(4)] $\alpha$ and $\beta$ are isomorphisms.
 
 \item[(5)]  $v'_1 u'_1 = \Id + \beta u_2 M_1^{-1} v_1$ and $v'_2 u'_2 = \Id + \alpha u_1 M_2^{-1} v_2$.
 
 \item[(6)]  $v_2\alpha = M_1 M_2 v_1$ and $v_1\beta = M_2 M_1v_2$. \qed
 
 \end{itemize}
 \end{thm}
 
 Let us analyze this result in a series of remarks.
 
 \begin{rems}
 (a)  The relations (3) imply that $\alpha$ and $\beta$ are expressed through the rest of
 the generators. 
 
 \vskip .2cm
 
 (b) The relations (2) mean that the remaining generators form, up to replacing some arrows with their negatives, a double representation
 of $\2^{3-1}$ in $\Vc$, i.e.,  we get a  commutative square
  \[
\xymatrix{
&E_{2,1}\ar^{-v_1}[dr]&
\\
E_3\ar^{v'_1}[ur] \ar_{v'_2} [dr] && E_{1,1,1}\\
&E_{1,2}. \ar_{v_2}[ur]&
}
\]
and another commutative square obtained from  the ``$u$-part'' of the diagram. 
 
 \vskip .2cm
 
 (c) Further,  let $\Fc\in\Perv(\Sym^3_0(\CC), S^0, \Vc)$. 
 Then  the commutative square   in (b) obtained from the diagram 
 \eqref{eq:diag-P3} corresponding to $\Fc$, precisely describes the constructible sheaf $R_{K^0}(\Fc)$ on the
 real skeleton $K^0=\Sym^3_0(\RR)$. This means that $E_\alpha$ is the stalk of $R_{K^0}(\Fc)$ at the stratum $K^0_\alpha = K^0\cap K_\alpha$, $\alpha\in\2^{3-1}$,
 while the $v_i, v'_i$ are the generalization maps, 
 see  \cite{GM-cusp} \S II.4. 
 
 \vskip .2cm
 
 (d) From this point of view,  the rest  of relations in Theorem \ref{thm:GM} can be  viewed as giving a precise meaning
 to the concept of a bi-sheaf on $K^0$, by specifying  the compatibility conditions 
 between the sheaf structure ($v$-part) and the
 cosheaf structure ($u$-part). 
 
 \end{rems}

 
 \paragraph{Granger-Maisonobe relations from a braided bialgebra.}\label{par:GM-braided}
 
 \begin{thm}
 Let $A$ be a primitive  bialgebra in a braided monoidal category $(\Vc, \otimes, R)$.
  Let us associate to it the diagram \eqref{eq:diag-P3} by:
  \[
 \xymatrix{
&A_2\otimes A_1\ar@<.5ex>[dr]^{\Delta_{1,1}\otimes\Id}
\ar@<.5ex>[dl]^{-\mu_{2,1}}
 \ar@<-.5ex>[dd]_{R}
&
\\
A_3\ar@<.5ex>[ur]^{-\Delta_{2,1}} \ar@<-.5ex>[dr]_{\Delta_{1,2}} && A_1\otimes A_1\otimes A_1
\ar@<-.5ex>[dl]_{-\Id\otimes\mu_{1,1}}
\ar@<.5ex>[ul]^{-\mu_{1,1}\otimes\Id}
 \\
&A_1\otimes A_2. \ar@<-.5ex>[ur]_{\Id\otimes\Delta_{1,1}}
\ar@<-.5ex>[ul]_{\mu_{1,2}}
 \ar@<-.5ex>[uu]_{R}
&
}
\]
That is, up to signs, the maps $u_i, u'_i$ are given by the multiplication, $v_i, v'_i$ by the comultiplication 
 (so they form the bibar-cube $\BB_3(A)$ in $\Vc$) while
  $\alpha =  R_{A_2, A_1}$ and $\beta =  R_{A_1, A_2}$ are  the
   braiding isomorphisms.
 Then this diagram
  satisfies the relations of Theorem \ref{thm:GM}. 
 \end{thm}

 \noindent{\sl Proof:} The {\bf  condition (1)} is fulfilled, in virtue of \eqref{uv:n=2},  with
\[
M_1 = - R_{A_1,A_1}\otimes \Id,\quad  M_2 = -\Id \otimes  R_{A_1,A_1}. 
\]
As in the case $n=2$, these minus signs are an instance of the braid group action being twisted by the sign character
(Remark \ref{rem:tw-sign}). 

 \bigskip

The {\bf conditions (4)} are obvious. 
The {\bf relations (2)} follow from associativity and coassociativity. 
Let us check the remaining relations
 directly.

 \bigskip
 
\noindent {\bf Relations (5). } Let us check the relation
\[
v_1'u_1' = \Id + \beta u_2 M_1^{-1} v_1: F_1\lra F_1,
\]
that is,
\[
(-\Delta_{2,1})(-\mu_{2,1}) \,=\, \Id + R_{A_1, A_2} (-\Id\otimes\mu_{1,1}) (-R^{-1}_{A_1, A_1}\otimes \Id) (\Delta_{1,1}\otimes\Id):
A_2\otimes A_1\lra A_2\otimes A_1. 
\]
Let $x = a_2\otimes a\in A_2\otimes A_1 = F_1$. 
As traditional in the theory of Hopf algebras (``Sweedler notation''), we write symbolically 
\[
\Delta(a_2)_{11} = b\otimes c,
\]
 meaning by $b\otimes c$ a sum of the form $\sum b_i\otimes c_i$. 
Let us also write, in the same way,
\[
R(c\otimes a) = c'\otimes a', \quad R=R_{A_1, A_1}. 
\] 
Then 
\[
\begin{gathered}
v_1'u_1' (x) = [\Delta( a_2a)]_{2,1} = 
[(a_2\otimes 1 + \Delta(a_2)_{11} + 1 \otimes a_2) (a\otimes 1 + 1 \otimes a)]_{2,1} =
\\
= a_2 \otimes a + \Delta(a_2)_{11}(a\otimes 1)  = a_2 \otimes a + ba'\otimes c',
\end{gathered}
\]
i.e.,
\[
(v'_1u'_1 -\Id) (x) = ba'\otimes c'.
\]
 On the other hand
\[
\beta u_2 M_1^{-1} v_1(x) = \beta u_2 M_1^{-1}(  b \otimes c\otimes a) = 
R_{A_1, A_2} (\tc\otimes \tb a)
\]
where 
\[
R^{-1}(b\otimes c) = \tc\otimes \tb
\]
We have to show that 
\be
\label{eq:(1.2.1)}
ba'\otimes c' = R_{A_1, A_2} (\tc\otimes \tb a). 
  \ee
We have a commutative diagram (naturality of the braiding)
\[
\xymatrix{
A_1\otimes (A_1\otimes A_1)
\ar[d]_{1\otimes \mu}
\ar[rr]^{R_{A_1, A_1\otimes A_1}}&&(A_1\otimes A_1)\otimes A_1
\ar[d]^{\mu\otimes 1}
\\
A_1\otimes A_2 \ar[rr]_{R_{A_1, A_2}}&& A_2\otimes A_1. 
} 
\]
 On the other hand, $R_{A_1, A_1\otimes A_1 }$ equals the composition
\be\label{eq:(1.2.2)}
A_1\otimes A_1 \otimes A_1 \overset{R_{12}}\lra A_1\otimes A_1 \otimes A_1 \overset{R_{23}}\lra A_1\otimes A_1 \otimes A_1,
 \ee
whence 
\[
R_{A_1\otimes A_1, A_1}R_{12}^{-1} = R_{23}.
\]
It follows:
\[
R_{A_1, A_2}(1\otimes \mu)R_{12}^{-1} = 
(\mu\otimes 1)R_{A_1, A_1\otimes A_1}R_{12}^{-1} = (\mu\otimes 1)R_{23},
\]
thus
$$
R_{A_1, A_2}(1\otimes \mu)R_{12}^{-1}( b\otimes c\otimes a) = 
(\mu\otimes 1)R_{23}( b\otimes c\otimes a), 
$$
and this is the sought for relation \eqref{eq:(1.2.1)}.  The second relation in (5) is proved similarly.

\bigskip

\noindent {\bf Relations (3). }
 Let us check that 
\be\label{eq:(1.3.1)}
v_2'u_1' = - \alpha + u_2v_1: F_1\lra F_2,
 \ee
 i.e.,
 \[
 \Delta_{1,2} (-\mu_{2,1}) = -R_{A_2, A_1} + (-\Id\otimes\mu_{1,1}) (\Delta_{1,1}\otimes\Id): A_2\otimes A_1\lra
 A_1\otimes A_2. 
 \]
 Dropping the minus signs, we write this as
 \be\label{eq:1.3.1-short}
 \Delta_{1,2}\mu_{2,1} = R_{A_2, A_1} + (\Id\otimes\mu_{1,1})(\Delta_{1,1}\otimes\Id). 
 \ee
 As before, let $a_2\otimes a\in A_2\otimes A_1$ and let $\Delta(a_2)_{11}=b\otimes c$. Then
 \[
 \begin{gathered}
 \Delta_{1,2}\mu_{2,1}(a_2\otimes a) = [(a_2\otimes 1 + b\otimes c + 1 \otimes a_2)(a\otimes 1 + 1 \otimes a)]_{12}=
 \\
 = (1\otimes a_2)(a\otimes 1)  + b\otimes ca.
 \end{gathered}
 \]
 But, by definition of the product in $A\otimes A$,
 \[
  (1\otimes a_2)(a\otimes 1)  = (\mu\otimes\mu) (1\otimes R_{A_2, A_1}(a_2\otimes a)\otimes 1) =  R_{A_2, A_1}(a_2\otimes a). 
 \]
 So the LHS of \eqref{eq:1.3.1-short}, applied to $a_2\otimes a$, is equal to
 \be\label{eq:Rbca}
 R_{A_2, A_1}(a_2\otimes a) + b\otimes ca, 
 \ee
 with the first summand  in \eqref{eq:Rbca} being the same as the  value, at $a_2\otimes a$, of 
 the first summand in the RHS of  \eqref{eq:1.3.1-short}. 
 Now, the second summand in  the RHS of  \eqref{eq:1.3.1-short}, applied to $a_2\otimes a$, is
 \[
 (\Id\otimes\mu_{1,1})(b\otimes c\otimes a) = b\otimes ca,
 \]
 so it coincides with the  second summand of 
 \eqref{eq:Rbca} as well, thus proving \eqref {eq:1.3.1-short} and \eqref{eq:(1.3.1)}. The second relation in (3) is proved
 similarly.

\bigskip 

\noindent{\bf Relations (6):}  Let us check the identity

\be\label{eq:(1.3.2)}
v_2\alpha = M_1M_2 v_1: F_1\lra E,
 \ee
 i.e.,
 \[
 (\Id\otimes\Delta_{1,1}) R_{A_2, A_1} = (-\Id\otimes R_{A_1, A_1}) (-R_{A_1, A_1}\otimes\Id) (\Delta_{1,1}\otimes \Id):
 A_2\otimes A_1 \lra A_1\otimes A_1\otimes A_1. 
 \]
 This follows from the commutative diagram (naturality of the braiding)
 \[
 \xymatrix{
 A_2\otimes A_1
 \ar[d]_{R_{A_2, A_1}}
  \ar[rr]^{\hskip -0.7cm \Delta_{1,1}\otimes\Id} && (A_1\otimes A_1)\otimes A_1
  \ar[d]^{R_{A_1\otimes A_1, A_1}}
 \\
 A_1\otimes A_2 \ar[rr]_{\hskip -0.4cm \Id\otimes\Delta_{1,1}}&& A_1\otimes(A_1\otimes A_1)
 }
 \]
 and the fact  that 
 \[
 R_{A_1\otimes A_1, A_1} = (\Id\otimes R_{A_1, A_1}) (R_{A_1, A_1}\otimes\Id). 
 \]
 The second relation in (6) is proved similarly. 
 

  \appendix
  
  \renewcommand\theparagraph{(\Alph{section}.\arabic{paragraph})}
  \numberwithin{equation}{section}

  \section {Shuffle complexes}\label{app:shuff-com}
  
  \paragraph{The permutohedron and its faces.} Let $n\geq 0$. By $\Sigma_n$ we denote the symmetric group
  $\Aut\{1,\cdots, n\}$. The $n$th {\em permutohedron} is the convex polytope in $\RR^n$ defined as
  \[
  P_n = \Conv\bigl\{ [s], s\in\Sigma_n\bigr\}, \quad [s] := (s(1), \cdots, s(n))\in\RR^n. 
  \]
The following is well known.

\begin{prop}\label{prop:Pn-faces}
(a) $P_n$ is a convex polytope of dimension $n-1$ with $n!$ vertices $[s]$, $s\in\Sigma_n$. 

\vskip .2cm

(b) Faces of $P_n$ (of all dimensions) are labelled by sequences $(I_1,\cdots, I_m)$ of subsets $I_\nu\subset \{1,\cdots,n\}$
which form a disjoint decomposition of $ \{1,\cdots,n\}$. The face $[I_1,\cdots, I_m]$ corresponding to $(I_1,\cdots, I_m)$
has vertices $[s]$ for all permutations $s$ obtained by ordering elements in the $I_\nu$ in all possible ways.
Thus 
\[
[I_1,\cdots, I_m] \,\simeq \, P_{|I_1|} \times\cdots\times P_{|I_m|}
\]
is a product of permutohedra and has dimension $n-m$. \qed
\end{prop}

\paragraph{The shuffle polytope and the shuffle complex.} Let $(d_1,\cdots, d_m)$ be a sequence of non-negative
integers and $d=d_1+\cdots + d_m$. Recall that a $(d_1,\cdots, d_m)$-{\em shuffle} is a permutation $s\in\Sigma_d$
whiose restriction onto  every interval
\be\label{eq:intervals}
d_1+\cdots + d_j, \,\,d_1+\cdots + d_j+1, \,\,\cdots, \,\,  d_1+\cdots + d_{j+1}-1, \quad j=1, \cdots, m-1, 
\ee
is monotone (order preserving). The set of $(d_1,\cdots, d_m)$-shuffles will be denoted 
$\Sha_{d_1,\cdots, d_m}\subset \Sigma_d$. 

\begin{ex}
In examples it is convenient to label integers from the intervals \eqref{eq:intervals} by numbered letters
such as
\[
a_1,\cdots, a_{d_1}, \, b_1, \cdots, b_{d_2}, \, c_1, \cdots, c_{d_3}, \text{ etc.}
\]
Then a shuffle is a just reordering of these letters which preserves the order of $a$'s, the order of $b$'s etc. For example,
$a_1b_1a_2b_2$ and $b_1b_2 a_1 a_2$ are $(2,2)$-shuffles. 
\end{ex}

\begin{Defi}
The {\em shuffle polytope} $N_{d_1,\cdots, d_m}\subset \RR^d$ is the convex hull of the points 
$[s]$, $s\in\Sha_{d_1,\cdots, d_m}$. 
\end{Defi}

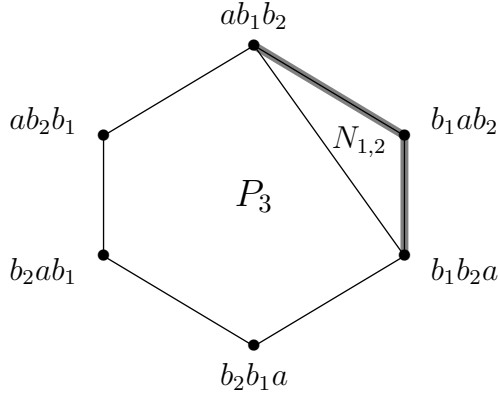
\begin{figure}[H]
 \centering
 \begin{tikzpicture}[scale=.4, baseline=(current bounding box.center)]
 
 \node(ab1b2) at (0,5){$\bullet$}; 
 \node (b1ab2) at (5,2) {$\bullet$}; 
 \node (b1b2a) at (5,-2){$\bullet$}; 
 \node (b2b1a) at (0,-5){$\bullet$}; 
 \node (b2ab1) at (-5,-2) {$\bullet$}; 
 \node (ab2b1) at (-5,2) {$\bullet$}; 
 
 \draw [line width=0.5] (ab1b2.center) -- (b1ab2.center) -- (b1b2a.center)
 -- (b2b1a.center) -- (b2ab1.center) -- (ab2b1.center) -- (ab1b2.center)
 ;
 
 \draw [line width=3,  opacity=0.5] (ab1b2.center) -- (b1ab2.center) -- (b1b2a.center); 
  \draw [line width=0.5] (ab1b2.center)  -- (b1b2a.center); 
  
  \node at (0,6) {$ab_1b_2$};
 \node at (7,2.5) {$b_1 a b_2$};
 \node at (7,-2.5) {$b_1  b_2a$};  
 \node at (0,-6) {$b_2b_1 a$};     
  \node at (-7,-2.5) {$ b_2a b_1$};  
  \node at (-7,2.5) {$ a b_2 b_1$};
  
  \node at (0,0) {\large$P_3$}; 
  
  \node at (3.5, 1.8) {$N_{1,2}$};

 \end{tikzpicture}
 \caption{The  hexagon $P_3$ and the triangle $N_{1,2}$.}\label{fig:N12}
 \end{figure}

This  $N_{d_1,\cdots, d_m}\subset P_d$, see Fig. \ref{fig:N12} for an illustration. 

\begin{rem}
As pointed out in \cite{GKZ}, Ch. 12,  \S 2, see also \cite{KV}, the 2-index shuffle polytope $N_{d_1, d_2}$ is 
the same as the {\em resultohedron},
i. e., the Newton polytope of the resultant $R(f_1, f_2)$ of two indeterminate polynomials of degrees $d_1, d_2$:
\[
f_1(x) = u_0 x^{d_1} + u_1 x^{d_1-1} + \cdots + u_{d_1}, \quad f_2(x) = v_0 x^{d_2} + v_1 x^{d_2-1} + \cdots + v_{d_2}.
\]
 An arbitrary $N_{d_1,\cdots, d_m}$ can also be interpreted in a similar way. For this we introduce
$m$ indeterminate polynomials $f_1(x), \cdots, f_m(x)$ of degrees $d_1,\cdots, d_m$ and form the expression
\[
\Rc (f_1, \cdots, f_m) \,=\,\prod_{i<j} R(f_i, f_j). 
\]
Then $N_{d_1,\cdots, d_m}$ can be identified with the Newton polytope of $\Rc$, i.e., with the Minkowski
sum of the resultohedra $N_{d_i, d_j}$. We do not prove this fact here since we will not need it. 
\end{rem}

\begin{Defi}
The {\em shuffle complex} $\ShD_{d_1,\cdots. d_m}$ is the union of all the faces of $N_{d_1,\cdots, d_m}\subset P_d$,
which are also faces of $P_d$.
\end{Defi} 

\begin{exas}
(a) If $(d_1,\cdots, d_m)=(1,\cdots, 1)$, then
\[
N_{1,\cdots, 1} \,=\, \ShD_{1,\cdots, 1}\, =\, P_m. 
\]

(b) $\ShD_{1,2}$ is the union of the two edges of the hexagon $P_3$ represented by thick lines in Fig. \ref{fig:N12}.

\vskip .2cm

(c) $\ShD_{2,2}$ is the complex in Fig. \ref{fig:Sh-22}. 

 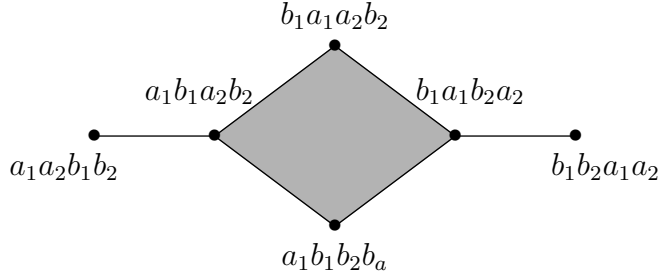
\begin{figure}[H]
 \centering
 \begin{tikzpicture}[scale=.4, baseline=(current bounding box.center)]
 
 \node (b1a1b2a2) at (0,3){$\bullet$}; 
 \node (a1b1a2b2) at (0,-3){$\bullet$}; 
 \node (a1a2b1b2) at (-8,0) {$\bullet$}; 
 \node (a1b1a2b2) at (-4,0) {$\bullet$}; 
 \node (b1a1b2a1) at (4,0) {$\bullet$}; 
 \node (b1b2a1a2) at (8,0) {$\bullet$}; 
  
  \draw[line width = 0.5] (-8,0) -- (-4,0); 
    \draw[line width = 0.5] (8,0) -- (4,0); 
    \draw[line width = 0.5] (-4,0) -- (0,3);   
    \draw[line width = 0.5] (0,3) -- (4,0);    
    \draw[line width = 0.5] (-4,0) -- (0,-3);  
    \draw[line width = 0.5] (0,-3) -- (4,0);  

 \filldraw[opacity = 0.3] (-4,0) -- (0,3) -- (4,0) -- (0,-3) -- (-4,0); 
 
\node at (-9, -1) {$a_1a_2b_1b_2$};  
 \node at (9, -1) {$b_1b_2a_1a_2$};  
 \node at (0,4) {$b_1 a_1a_2b_2$};
 \node at (0,-4) {$a_1 b_1b_2b_a$};
 
 \node at (-4.5,  1.5){$a_1b_1a_2b_2$}; 
 
  \node at (4.5,  1.5){$b_1a_1b_2a_2$};

 \end{tikzpicture}
 \caption{The complex $\ShD_{2,2}$.}\label{fig:Sh-22}
 \end{figure}
\end{exas}

\begin{prop}\label{prop:ShD-contr}
Each $\ShD_{d_1,\cdots, d_m}$ is contractible.
\end{prop}

 The proof will be given in the next \S \ref{par:tri-prod}. 

\paragraph{Relation to the triangulated product of simplices.}\label{par:tri-prod}. 
For any $p\geq 0$ we consider the standard coordinate space $\RR^{p+1}$ with basis  $e_0,\cdots, e_p$.
Let  $\Delta^p = \Conv\{e_0,\cdots. e_p\}\subset \RR^{p+1}$ be the standard $p$-dimensional simplex witn vertices $e_0,\cdots, e_p$.
It is convenient to identify the set of vertices of $\Delta^p$ with the standard ordinal $[p]=\{0,1,\cdots, p\}$. 

\vskip .2cm

Given $d_1,\cdots, d_m$ as before,  we consider the product $\Pi=\Delta^{d_1}\times \cdots\times \Delta^{d_m}$.  
It is a convex polytope of dimension $d=d_1+\cdots+d_m$. 
We
denote its vertices by $e_\sigma$ where $\sigma$ runs in the set $[d_1] \times \cdots\times [d_m]$
which we consider with the partial order given by  the product of the standard orders on the $[d_j]$. 
We denote by $T$ the {\em standard triangulation} of $\Pi$ well known in the theory of simplicial sets,
see, e.g., \cite{gelfand-manin} \S I.1.5. By definition, simplices of $T$  (of any dimension $m$) have the form
\[
\Delta_\tau \,=\, \Conv \{ e_{\tau_0},  e_{\tau_1}, \cdots, e_{\tau_m}\}, 
\quad  \tau = (\tau_0 <\cdots < \tau_m), \quad \tau_i\in \prod_{j=1}^m  [d_j],
\]
i.e., they are parametrized by chains in the poset $\prod_{j=1}^m  [d_j]$. 
Proposition \ref{prop:ShD-contr} follows from the next more precise fact.

\begin {prop}\label{prop:ShD=T*}
$\ShD_{d_1,\cdots, d_m}$ is identified with $T^\vee$, the polyhedral complex Poincar\'e dual to $T$.
\end{prop}

Explicitly, to construct $T^\vee$, we take one vertex for each maximal simplex of $T$, these vertices
being joined by an edge when the simplices are adjacent and so on.  That is,
cells of $T^\vee$ correspond to those  simplices of $T$ that {\em lie inside} $\Pi$ (i.e., do not lie on its boundary),
with the opposite order of inclusion.

\vskip .2cm

\begin{figure}[H] 
   
   \begin{tikzpicture}
   
   \node at (0,0){$\bullet$};
   \draw  (0,0) -- (5,0); 
   
   \draw (0,0) -- (0,3); 
  
  
   
  \node at (1,0){\tiny$\bullet$}; 
    \node at (2,0){\tiny$\bullet$}; 
      \node at (3,0){\tiny$\bullet$}; 
        \node at (4,0){\tiny$\bullet$}; 
          \node at (5,0){\tiny$\bullet$};

            \node at (0,1){\tiny$\bullet$}; 
              \node at (0,2){\tiny$\bullet$}; 
                \node at (0,3){\tiny$\bullet$}; 
  
  \node at (5,3){$\bullet$};  
  \draw [dashed] (0,3) -- (5,3); 
    \draw [dashed] (0,2) -- (5,2); 
      \draw [dashed] (0,1) -- (5,1); 
 
   \draw [dashed] (1,0) -- (1,3); 
     \draw [dashed] (2,0) -- (2,3); 
       \draw [dashed] (3,0) -- (3,3); 
         \draw [dashed] (4,0) -- (4,3); 
           \draw [dashed] (5,0) -- (5,3); 
           
  \draw[ line width = 1] (0,0) -- (0,1) --(3,1) -- (3,2) -- (5,2) -- (5,3);   
  
  
  \node at (-0.3, 0.5){$b_1$}; 
\node at (0.5, 1.2){$a_1$};   
 \node at (1.5, 1.2){$a_2$};    
  \node at (2.5, 1.2){$a_3$};   

\node at (3.2, 1.5){$b_2$};   

  \node at (3.5, 2.2){$a_4$};         
  \node at (4.5, 2.2){$a_5$};    
  
  \node at (5.2, 2.5){$b_3$};   
  
  \node at (-0.2, -0.2) {\tiny 0}; 
  
   \node at (1, -0.2){\tiny 1}; 
    \node at (2, -0.2){\tiny 2}; 
  \node at (3, -0.2){\tiny 3}; 
    \node at (4, -0.2){\tiny 4}; 
      \node at (5, -0.2){\tiny 5};  

\node at (-0.2, 1) {\tiny 1}; 
\node at (-0.2, 2) {\tiny 2}; 
\node at (-0.2, 3) {\tiny 3}; 

\node at (2.5,4) {\large $b_1 a_1a_2a_3b_2a_4a_5b_3$}; 
   
   \end{tikzpicture}\quad\quad\quad\quad
   \begin{tikzpicture}
   
   \node at (0,0){$\bullet$};
   \draw  (0,0) -- (5,0); 
   
   \draw (0,0) -- (0,3); 
  
  
   
  \node at (1,0){\tiny$\bullet$}; 
    \node at (2,0){\tiny$\bullet$}; 
      \node at (3,0){\tiny$\bullet$}; 
        \node at (4,0){\tiny$\bullet$}; 
          \node at (5,0){\tiny$\bullet$};

            \node at (0,1){\tiny$\bullet$}; 
              \node at (0,2){\tiny$\bullet$}; 
                \node at (0,3){\tiny$\bullet$}; 
  
  \node at (5,3){$\bullet$};  
  \draw [dashed] (0,3) -- (5,3); 
    \draw [dashed] (0,2) -- (5,2); 
      \draw [dashed] (0,1) -- (5,1); 
 
   \draw [dashed] (1,0) -- (1,3); 
     \draw [dashed] (2,0) -- (2,3); 
       \draw [dashed] (3,0) -- (3,3); 
         \draw [dashed] (4,0) -- (4,3); 
           \draw [dashed] (5,0) -- (5,3); 
           
  \draw[ line width = 1] (0,0) -- (0,1) --(3,1) -- (3,2) -- (5,2) -- (5,3);   
  
  \draw[line width=1] (0,0) -- (1,0) --(1,1) --(1,2) --(3,2); 
  
  \draw[line width=1] (4,2) -- (4,3) -- (5,3);

 \filldraw [opacity=0.3] (0,0) -- (1,0) -- (1,1) --(0,1) -- (0,0);  
 
  \filldraw [opacity=0.3]  (1,1) -- (3,1) -- (3,2) -- (1,2) -- (1,1); 
  
    \filldraw [opacity=0.3] (4,2) -- (5,2) -- (5,3) -- (4,3) -- (4,2); 
    
    \node at (1,1){$\bullet$}; 
    
     \node at (3,2){$\bullet$}; 
     
          \node at (4,2){$\bullet$}; 
          
    \node at (-0.2, -0.2) {\tiny 0}; 
  
   \node at (1, -0.2){\tiny 1}; 
    \node at (2, -0.2){\tiny 2}; 
  \node at (3, -0.2){\tiny 3}; 
    \node at (4, -0.2){\tiny 4}; 
      \node at (5, -0.2){\tiny 5};  

\node at (-0.2, 1) {\tiny 1}; 
\node at (-0.2, 2) {\tiny 2}; 
\node at (-0.2, 3) {\tiny 3};     

\node at (-0.5,0){$\tau_0$};   

\node at (0.7, 1.3){$\tau_1$};  
 
 \node at (3.1, 2.2){$\tau_2$};   
 
 \node at (3.8, 2.2){$\tau_3$};    

\node at (5.2, 3.2){$\tau_4$};  

\node at (0.5, 0.5) {$\Pen_1$}; 

\node at (2,1.5) {$\Pen_2$}; 

\node at (3.5, 1.7) {$\Pen_3$}; 

\node at (4.5, 2.5) {$\Pen_4$}; 
 
   \end{tikzpicture}
   
   \caption{A  taxicab path $\pen$ representing a shuffle and a block path $\Pen$.}\label{fig:taxicab}
   \end{figure}
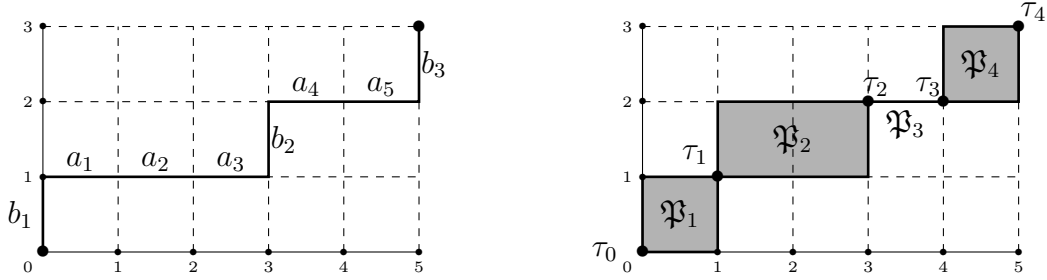

\noindent{\sl Proof of Proposition \ref{prop:ShD=T*}:}
Note first that maximal  chains $\tau$ in $\prod [d_j]$
(with $m=d$) are in bijection with $\Sha_{d_1,\cdots, d_m}$. Indeed, such chains can be seen as 
``lattice paths'' or ``taxicab paths'' $\pen$
in the $m$-dimensional parallelotope $\prod_{j=1}^m [0,d_j]$, going from the minimal vertex $(0,...,0)$ to the maximal one
$(d_1,\cdots. d_m)$. 
This is because, $\tau$ being maximal, each $\tau_{\nu+1}$ is obtained from $\tau_\nu$ by going exactly
one step in one of the $l$ directions.  See  the left part of Fig. \ref{fig:taxicab} for $m=2$.

   Further,  consider simplices $\Delta_\tau$ corresponding to arbitrary chains $ \tau = (\tau_0 <\cdots < \tau_m)$.
    The condition for  $\Delta_\tau$ to lie inside $\Pi$
    is that $\tau_0$ is the minimal element $(0,\cdots, 0)$ of $\prod [d_j]$
   and $\tau_m$ is the maximal element $(d_1, \cdots, d_m)$.

 Notice that such $\tau$    are in bijection with {\em block paths}  in the parallelotope  $\prod_{j=1}^m [0,d_j]$,
 i.e., sequences $\Pen = (\Pen_1, \cdots, \Pen_m)$ of 
   lattice parallelotopes (possibly degenerate, i.e., of smaller dimension) 
inside  $\prod_{j=1}^m [0,d_j]$ with the following properties:

\begin{itemize}
\item[(1)] The minimal (with respect to the product partial order on $\RR^m$) vertex of $\Pen_1$ is $(0,\cdots, 0)$.

\item[(2)] For $i=1, \cdots, m-1$, the minimal vertex of $\Pen_{i+1}$ coincides with the maximal vertex of $\Pen_i$.

\item[(3)] The maximal vertex of $\Pen_m$ is $(d_1,\cdots, d_m)$. 
\end{itemize}

 More precisely, since $\tau_{\nu+1} > \tau_\nu$, the 
difference $\tau_{i+1}-\tau_\nu$, as a vector in $\RR^m$, has the form $(q_1,\cdots, q_m)$ with $q_i\in\ZZ_{\geq 0}$,
and we have a parallelotope of these dimensions with minimal vertex $\tau_\nu$ and maximal vertex $\tau_{\nu+1}$.
 See the right part of  Fig. \ref{fig:taxicab} for $m=2$. 
 
 It remains to notice that the poset of such block paths $\Pen$ is identified with the poset of faces of
 $\ShD_{d_1,\cdots, d_m}$. For $m=2$ the reader can  compare with the description of all faces of the shuffle polytope
 $N_{d_1, d_2}$ in \cite{GKZ} Ch. 12, \S 2 in terms of  ``labirinths'':  block paths are precisely the
 labirinths corresponding to faces of  $N_{d_1, d_2}$ which are also faces of $P_{d_1+d_2}$.  
 
  This proves Propositions \ref{prop:ShD=T*} and  \ref{prop:ShD-contr}.

\paragraph{Relation to XO-tables.} Let $(d_1,\cdots, d_m)$ be as before. 

\begin{Defi}\label{def:XO-tab}
 By an {\em XO-table} of type
$(d_1,\cdots, d_m)$ we will mean an $l\times m$ matrix $\Xi$ for some $l\geq 1$ formed by X's and O's
such that:
\begin{itemize}
\item[(1)] The number of X's in the $j$th column of $\Xi$ is precisely $d_j$.

\item[(2)] Each row of $\Xi$ contains at least one X. 
\end{itemize}
The set of such tables will be denored $\XO(d_1,\cdots, d_m)$. 

\end{Defi}

The number $l=l(\Xi)$ will be called the {\em height} of $\Xi$. It can vary between $\max(d_1,\cdots, d_m)$
and $d_1+\cdots+ d_m$. 

\begin{figure}[H]
 \centering
 \begin{tikzpicture}[scale=.4, baseline=(current bounding box.center)]
 
 \node at (-10,0){$\bullet$}; 
  \node at (0,0){$\bullet$}; 
   \node at (10,0){$\bullet$}; 
   
   \draw (-10,0) -- (0,0) -- (10,0); 
   
 \node at (0,3) {
 $\begin{matrix}
 O&X\\
 X&O
 \\O&X
 \end{matrix}$  
 } ; 
 
  \node at (12, 1) {
 $\begin{matrix}
 X&O\\
 O&X
 \\
 O&X
 \end{matrix}$  
 } ; 
 
   \node at (-12, 1) {
 $\begin{matrix}
 O&X\\
 O&X
 \\
 X&O
 \end{matrix}$  
 } ; 
 
   \node at (5, -2) {
 $\begin{matrix}
 X&X\\
 O&X
  \end{matrix}$  
 } ; 
 
    \node at (-5, -2) {
 $\begin{matrix}
 O&X\\
 X&X
  \end{matrix}$  
 } ; 
   
   \end{tikzpicture}
 \caption{ XO-tables of type $(1,2)$, depicted on vertices and edges to indicate contractions.}\label{fig:XO}
 \end{figure}
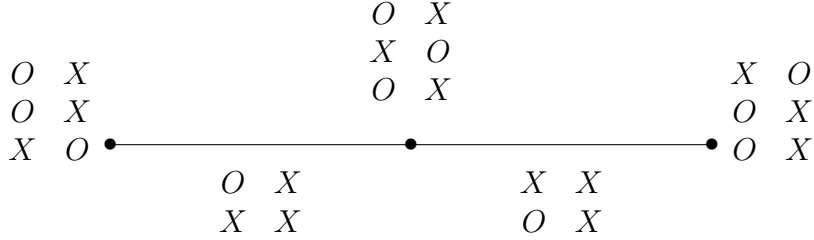

\begin{exas}
(a) XO-tables $\Xi$ of maximal height $d_1+\cdots + d_m$ are in bijection with $\Sha_{d_1,\cdots, d_m}$.
In this case there is exactly one X in every row, so  the columns give a shuffle. 

\vskip .2cm

(b) XO-tables of type $(1,2)$ are depicted in Fig. \ref{fig:XO}. 

\end{exas}

 To an XO-table $\Xi$ we can associate an integer matrix $\wt\Xi$ of the same size by replacing each X with $1$
 and each O with $0$. For any $i=1,\cdots, l=l(\Xi)$ we denote by $\wt\Xi_{i,\bullet}$ the $i$th row of $\Xi$,
 and similarly for $\Xi_{i\,\bullet}$. 
 
 \begin{Defi}
 Let $\Xi, \Xi'$ be two XO-table of the same type $(d_1,\cdots, d_m)$ of heights $l,l'$ respectively. We say that
 $\Xi'$ is a {\em contraction} of $\Xi$ and write $\Xi\leq\Xi'$, if there following holds:
 \begin{itemize}
 \item[(C)] 
 There is a surjective monotone map
 $\varphi: \{1,\cdots, l\} \to \{1,\cdots, l'\}$ such that for any $i'\in\{1,\cdots, l'\}$
 \[
 \wt\Xi'_{i',\bullet} = \sum_{i\in\varphi^{-1}(i')} \wt\Xi_{i, \bullet}
 \]
 (summation of integer vectors). 
 \end{itemize}
 We say that $\Xi$ is an {\em elementary contraction} of $\Xi$ if $\Xi\leq\Xi'$ and $l'=l-1$. 
 \end{Defi}

Since all the entries of $\wt\Xi$, $\wt\Xi'$ are $0$'s and $1$'s, the condition  (C) means that the rows $\Xi_{i,\bullet}$,
$i\in\varphi^{-1}(i')$, have mutually disjoint subsets of positions filled by X's, and the row $\Xi'_{i',\bullet}$
is obtained by putting X's in all of these positions. 

\begin{prop}
(a) $\leq$ is a partial order on $\XO_{d_1,\cdots, d_m}$ with minimal elements corresponding to shuffles
(XO-tables of height $d_1+\cdots + d_m$).

\vskip .2cm

(b) $(\XO_{d_1,\cdots, d_m},\leq) $ is isomorphic to the poset of faces of $\ShD_{d_1,\cdots, d_m}$. 
\end{prop}

\noindent{\sl Proof:} (a) is obvious. To prove (b), let us number the X's in an XO-table $\Xi$ from $1$ to $d$ lexicographically,
first from top to bottom in the first column, then similarly in the second column etc.  Then each row will pick up
a subset in $\{1,\cdots, d\}$. Denoting by $I_\nu$ the subset picked up by the $\nu$th row, $\nu=1,\cdots, l=l(\Xi)$,
we get a sequence $(I_1,\cdots, I_l)$ of subsets defining a face  $[I_1,\cdots, I_l]$ of $P_d$ by Proposition \ref{prop:Pn-faces}.
All the vertices of  $[I_1,\cdots, I_l]$  are $(d_1,\cdots, d_m)$-shuffles, so  $[I_1,\cdots, I_l]$ is a face of
$\ShD_{d_1,\cdots, d_m}$. This establishes the claimed isomorphism of posets.\qed

 \begin{Defi}\label{def:xo-cochain}
 The {\em XO-cochain complex}  of type $(d_1,\cdots, d_m)$ is  the cochain complex of $\k$-vector spaces defined as the graded space
\[
C^\bullet(\XO_{d_1,\cdots, d_m})\,=\, \bigoplus_{\Xi\in\XO_{d_1,\cdots, d_m}} \k[l(\Xi) - d],\quad d=d_1+\cdots + d_m, 
\]
with the differential  given by the alternating sum of elementary contractions. 

\end{Defi}

\begin{cor}\label{cor:CXO-contr}
The complex $C^\bullet(\XO_{d_1,\cdots, d_m})$ is quasi-isomorphic to $\k$ in degree $0$.
\end{cor}

\noindent{\sl Proof:} Indeed, it is identified with the cellular cochain complex of $\ShD_{d_1,\cdots, d_m}$
which is contractible by Proposition \ref {prop:ShD-contr}.\qed


\vskip 1cm

\small{

M.K.: Kavli IPMU, 5-1-5 Kashiwanoha, Kashiwa, Chiba, 277-8583 Japan, 
{\tt mikhail.kapranov@ipmu.jp}

\smallskip

 V.S.: Institut de Math\'ematiques de Toulouse, Universit\'e Paul Sabatier, 118 route de Narbonne, 
31062 Toulouse, France, 
 {\tt schechtman@math.ups-tlse.fr }

}

\ed